\pdfoutput=1
\documentclass[a4paper,11pt]{amsart}
\usepackage[english]{babel}
\usepackage{amsfonts}
\usepackage{amsthm}
\usepackage{amsmath}
\usepackage{amssymb}
\usepackage{enumerate}
\usepackage{mathtools}
\usepackage{tikz-cd}
\usepackage[all]{xy}
\usepackage{graphicx}
\usepackage{amsmath}
\usepackage{extpfeil}
\usepackage{comment}
\usepackage{float}
\usepackage[labelformat=empty]{caption}
\usepackage[toc]{appendix}
\usepackage[left=3cm,top=2.5cm,right=3cm,bottom=2.5cm]{geometry}
\usepackage{hyperref}
\usepackage{comment}
\usepackage{resizegather}
\usepackage[activate={true, nocompatibility}, final, tracking=true, kerning=true, spacing=true]{microtype}
\usepackage{xcolor}
\usepackage{cleveref}

\theoremstyle{plain}
\newtheorem{thm}{Theorem}[section]
\newtheorem{coroll}[thm]{Corollary}
\newtheorem{lemma}[thm]{Lemma}
\newtheorem*{claim}{Claim}

\newtheorem{example}[thm]{Example}
\newtheorem{prop}[thm]{Proposition}

\newtheorem{defn}[thm]{Definition}
\newtheorem{remark}[thm]{Remark}
\newtheorem{notn}[thm]{Notation}

\def\makeCal#1{%
\expandafter\newcommand\csname c#1\endcsname{\mathcal{#1}}}
\def\makeBB#1{%
\expandafter\newcommand\csname b#1\endcsname{\mathbb{#1}}}
\def\makeFrak#1{%
\expandafter\newcommand\csname f#1\endcsname{\mathfrak{#1}}}

\newcounter{int}
\loop
\addtocounter{int}{1}
\edef\y{\Alph{int}}%
\expandafter\makeCal\y
\expandafter\makeBB\y
\expandafter\makeFrak\y
\ifnum\value{int}<26%
\repeat
\makeatletter
\newcommand{\colim@}[2]{%
  \vtop{\m@th\ialign{##\cr
    \hfil$#1\operator@font colim$\hfil\cr
    \noalign{\nointerlineskip\kern1.5\ex@}#2\cr
    \noalign{\nointerlineskip\kern-\ex@}\cr}}%
}
\newcommand{\colim}{%
  \mathop{\mathpalette\colim@{\rightarrowfill@\textstyle}}\nmlimits@
}
\makeatother

\tikzset{
  symbol/.style={
    draw=none,
    every to/.append style={
      edge node={node [sloped, allow upside down, auto=false]{$#1$}}}
  }
}

\newcommand{\dash}{{\operatorname{-}}}

\microtypecontext{spacing=nonfrench}

\usepackage[normalem]{ulem}
\usepackage{pbox}

\DeclareMathOperator{\Coh}{Coh}

\newcommand{\dual}{\vee}

\DeclareMathOperator{\Filt}{Filt}

\DeclareMathOperator{\GL}{GL}

\DeclareMathOperator{\Map}{Map}

\DeclareMathOperator{\Pic}{Pic}

\DeclareMathOperator{\Spec}{Spec}

\DeclareMathOperator{\Gr}{Gr}

\DeclareMathOperator{\Pair}{Pair}
\DeclareMathOperator{\SL}{SL}
\DeclareMathOperator{\Quot}{Quot}
\DeclareMathOperator{\rk}{rk}
\DeclareMathOperator{\wt}{wt}

\begin{document}
\thispagestyle{empty}

\title{Moduli spaces of sheaves via affine Grassmannians}
\author{Daniel Halpern-Leistner, Andres Fernandez Herrero, and Trevor Jones}
\address{Department of Mathematics, Cornell University,
    310 Malott Hall, Ithaca, New York 14853, USA.}
\email{daniel.hl@cornell.edu}
\email{ajf254@cornell.edu}
\email{ttj9@cornell.edu}
\date{}
\maketitle
\begin{abstract}
We develop a new method for analyzing moduli problems related to the stack of pure coherent sheaves on a polarized family of projective schemes. It is an infinite dimensional analogue of Geometric Invariant Theory. We apply this to two familiar moduli problems: the stack of $\Lambda$-modules and the stack of pairs. In both examples, we construct a $\Theta$-stratification of the stack, defined in terms of a polynomial numerical invariant, and we construct good moduli spaces for the open substacks of semistable points. One of the essential ingredients is the construction of higher dimensional analogues of the affine Grassmannian for the moduli problems considered.
\end{abstract}

\tableofcontents

\begin{section}{Introduction}

The moduli stack $\Coh(X)$ of coherent sheaves on a projective scheme $X$ is of central interest in moduli theory. It has been used to formulate the non-abelian Hodge correspondence \cite{simpson-non-abelian-hodge}, it has been used to define and study algebraic Donaldson invariants of surfaces \cite{mochizuki-invariants-surfaces-book}, and several flavors of enumerative invariants for 3-folds \cite{thomas-3-fold, pt-stable-pairs}. It has also been used to study the moduli of varieties \cite{kollar2009hulls}, and has served as a particularly interesting testing ground for the minimal model program \cite{bayer-macri-mmp} and hyperk\"ahler geometry \cite{mukai-symplectic}.

As the examples illustrate, a thorough understanding of the structure of the stack $\Coh(X)$ is a fundamental building block for many other theories. We will focus on the open substack $\Coh^d(X)_P \subset \Coh(X)$ parameterizing sheaves that are \emph{pure} of dimension $d$ and have a fixed Hilbert polynomial $P$. In this paper, we develop a new method for studying the structure of $\Coh^d(X)_P$ and related stacks. We call it \emph{infinite dimensional Geometric Invariant Theory (GIT)}.

Let us recall the main structural results about $\Coh^d(X)$: Using an ample line bundle $\cO_X(1)$ on $X$, one defines a sheaf $\cF$ on $X$ to be Gieseker $\cO_X(1)$-semistable if and only if there is no proper subsheaf $\cF' \subset \cF$ whose reduced Hilbert polynomial is larger than that of $\cF$ (see Subsection \ref{subsection: hilbert polynomials background} below). The substack of semistable sheaves $\Coh^d(X)^{\rm ss}_P \subset \Coh^d(X)$ is open and admits a projective moduli space. Furthermore, $\Coh^d(X)^{\rm ss}$ is the open piece of a stratification of $\Coh^d(X)$ by locally closed substacks, where the strata parameterize unstable sheaves along with a canonical filtration, called the Harder-Narasimhan filtration \cite{harder_narasimhan, simpson-repnfundamental-I, nitsure-scheamtichn}. The strata are indexed by the Hilbert polynomials of the associated graded pieces of the Harder-Narasimhan filtration, which are themselves semistable.

This moduli problem has been previously studied using Geometric Invariant Theory (GIT) \cite{seshadri-unitary, gieseker_torsion_free, maruyama_moduli, simpson-repnfundamental-I}, which involves constructing and analyzing explicit families of sheaves over quot-schemes. Over the past few years, a new approach to these structures has developed that is more intrinsic to the moduli problem: the theory of good moduli spaces \cite{alper-good-moduli} and the theory of $\Theta$-stratifications \cite{halpernleistner2018structure}. We will review the key concepts in Subsection \ref{subsection: background theta stratifications}. From the intrinsic perspective, the main structural results for $\Coh^d(X)_P$ are: 
\begin{enumerate}[(S1)]
    \item To any relatively ample bundle on $X$, one can associate a polynomial-valued numerical invariant $\nu$ (\Cref{defn: poly numerical invariant on coh}) that defines a $\Theta$-stratification of $\Coh^d(X)_P$.
    \item If $S$ is defined over $\bQ$, then the semistable locus admits a good moduli space that is proper over $S$.
\end{enumerate}

% From the intrinsic perspective, the main structural result is:

% \begin{thm} \label{thm:A}
% Let $X$ be a projective scheme of finite presentation over a scheme $S$. To any relatively ample line bundle $\mathcal{O}_X(1)$ on $X$, one can associate a polynomial-valued numerical invariant $\nu$ (Definition \ref{defn: poly numerical invariant on coh}) that defines a $\Theta$-stratification of $\Coh^d(X)_P$. If $S$ is defined over $\bQ$, then the $\nu$-semistable locus $\Coh^d(X)_P^{\nu \dash \rm{ss}}$ admits a good moduli space that is proper over $S$.
% \end{thm}

Our methods give a new proof of these structural results for $\Coh^d(X)_P$ which more readily extends to other contexts. In the applications mentioned above, one often considers the moduli problem of a coherent sheaf $\mathcal{F}$ along with some auxiliary data. We will focus on two such variants: 
\begin{enumerate}[1)] 
\item The moduli stack of pairs $\Pair_\cA(X)$ parameterizing a pure coherent sheaf $\mathcal{F}$ of dimension $d$ on $X$ along with a homomorphism $\cA \to \mathcal{F}$ from a fixed coherent sheaf $\cA$ (often with $\cA=\cO_X$ or $\cA = \cO_X^{\oplus n}$ for some $n$ \cite{bradlow-pairs, thaddeus-verlinde, pt-stable-pairs, mochizuki-invariants-surfaces-book}).
\item The moduli stack $\Lambda\Coh^d(X)$ of pure coherent modules of dimension $d$ over a sheaf $\Lambda$ of rings of differential operators on $X$. 
\end{enumerate}
See Subsection \ref{section: related moduli problems} for precise definitions. Both moduli problems admit a morphism to $\Coh^d(X)$ that forgets the additional structure.
\begin{thm}[= Theorems \ref{thm: existence of weak theta stratification for pairs}, \ref{thm: moduli space pairs}] \label{thm:B}
Let $X$ be a projective scheme of finite presentation over a scheme $S$. To any relatively ample line bundle $\mathcal{O}_X(1)$ on $X$, one can associate a polynomial-valued numerical invariant $\nu$ on $\Coh^{d}(X)$ (Definition \ref{defn: poly numerical invariant on coh}). 

The restriction of the numerical invariant $\nu$ from $\Coh^d(X)$ defines a weak $\Theta$-stratification of the moduli stack $\Pair^d_\cA(X)_{P}$. If $S$ is defined over $\mathbb{Q}$, then the $\nu$-semistable locus admits a good moduli space that is proper over $S$. 

Furthermore, the same holds for a larger family $\nu^{(\delta)}$ of deformations of the numerical invariant on $\Pair^d_\cA(X)$ described in Definition \ref{defn: numerical invariant for pairs}.
\end{thm}

\begin{remark}
If we set $\mathcal{A}=0$, then  the forgetful morphism induces an identification $\Pair^{d}_{\mathcal{A}}(X) \xrightarrow{\sim} \Coh^{d}(X)$. In particular, \Cref{thm:B} recovers the classical structural results (S1) and (S2) above for the stack $\Coh^{d}(X)_{P}$.
\end{remark}

The main contribution of this paper is to develop the method of infinite dimensional GIT to prove \Cref{thm:B}. It is inspired by the theory of infinite dimensional symplectic reduction \cite{atiyah-bott-yangmills, donaldson-surfaces}. In the case of vector bundles on a curve, the differential geometric picture uses the description of the stack as a quotient of an infinite dimensional space of connections by an infinite dimensional gauge group. The analogous uniformization theorem in algebraic geometry describes this stack as the quotient of a Beilinson-Drinfeld Grassmannian, an ind-projective ind-scheme which is a relative version of the affine Grassmannian. We prove the theorem above by constructing a ``higher dimensional" analog of the Beilinson-Drinfeld Grassmannian for each moduli problem and applying ``Geometric Invariant Theory" to the action of the infinite-dimensional group of rational maps $X \dashrightarrow \GL_N$ on these ind-projective ind-schemes.

To our knowledge, the level of generality in which we construct moduli spaces of semistable pairs is new, although the question has been investigated in several more specific settings. The good moduli space of the semistable locus $\Pair^d(X)^{\nu \dash \mathrm{ss}}_{P}$ recovers the moduli spaces of Bradlow pairs considered in \cite{bradlow-pairs, thaddeus-verlinde} for curves, \cite{pt-stable-pairs} for Calabi-Yau threefolds, \cite{wandel-moduli-pairs} in the case of a complex projective variety, and \cite{kollar2009hulls} \cite[\S11]{bayer2019stability} \cite{lin2021decorated} when the stability parameter $\delta \gg 0$ (see also the related moduli spaces of coherent systems considered in \cite{le-potier-coherent-systems} for smooth projective varieties, and more recently in \cite{schmitt-coherent-systems}).

On the other hand, the stratification of $\Pair^d_{\cA}(X)_P$ is new, even in the more specific contexts mentioned above. One interesting feature is that the canonical filtration of an unstable point in $\Pair^d_\cA(X)_P$ is not necessarily ``convex" in the way that the usual Harder-Narasimhan filtration of an unstable sheaf or $\Lambda$-module is (see Example \ref{example: filtration nonconvex pairs}). This suggests that the stratification of $\Pair_\cA^d(X)_P$ cannot be constructed in the same way as the stratification of $\Coh^d(X)_P$.

Even for $\Coh^d(X)_P$, we believe that the method of infinite dimensional GIT is of intrinsic interest, and is conceptually cleaner than the classical approach. We hope that by developing these methods in the context of familiar moduli problems, this paper lays the foundation for a broad range of applications.

The method of infinite dimensional GIT originated in the forthcoming paper \cite{gauged_theta_stratifications} developing one such application, to the computation of gauged Gromov-Witten invariants. Another application, to the moduli of singular principal $G$-bundles on higher dimensional varieties, appears in \cite{rho-sheaves-paper}. Finally, the original version of this paper included applications to the moduli of $\Lambda$-modules, recovering results of \cite{simpson-repnfundamental-I} in the case when the base $S$ is Noetherian.
\begin{thm}[= Theorem \ref{thm: theta stratification pure sheaves}] \label{thm:A}
In the context of Theorem \ref{thm:B}, the restriction of the numerical invariant $\nu$ from $\Coh^d(X)$ defines a weak $\Theta$-stratification of the moduli stack $\Lambda\Coh^d(X)_{P}$. If $S$ is defined over $\mathbb{Q}$, then the $\nu$-semistable locus admits a good moduli space that is separated over $S$. 
\end{thm}
We summarize these results in Section \ref{section: applications lambda modules}, but the details will appear in the forthcoming paper \cite{companion-paper}.

%Our good moduli space for the semistable locus $\Lambda\Coh^d(X)^{\nu \dash \mathrm{ss}}_{P}$ recovers the moduli space of $p$-semistable $\Lambda$-modules constructed in \cite{simpson-repnfundamental-I} in the Noetherian case.

\subsection{Infinite dimensional GIT}
The theory of $\Theta$-stability allows one to formulate notions of semistability and canonical filtrations for points in an algebraic stack $\cX$. We will give a brief summary here, and refer the reader to Subsection \ref{subsection: background theta stratifications} for a more detailed discussion of the key ideas and results.

The input of the theory is a structure on $\cX$ called a \emph{numerical invariant} (Definition \ref{defn: numerical poly invariant}). Geometric invariant theory (GIT) is the special case where $\cX = Z/G$ for a reductive group $G$ acting on a projective scheme $Z$. The notion of semistability can be defined using the Hilbert-Mumford criterion with respect to a $G$-equivariant ample line bundle $L$ on $Z$. Combined with some additional data (a norm on cocharacters of $G$), $L$ determines a numerical invariant on $Z/G$, which encodes the Hilbert-Mumford criterion and also defines the Hesselink-Kempf-Kirwan-Ness stratification of the unstable locus. In our examples, semistability is defined via an intrinsic version of the Hilbert-Mumford criterion, but instead of a single line bundle $L$, we use an infinite sequence of line bundles $L_n \in \Pic(\cX)$ indexed by $n \in \bZ$. The resulting numerical invariant is a function in $n$ (in fact polynomial), and semistability and Harder-Narasimhan filtrations are characterized by the asymptotic behavior of this function as $n \to \infty$.

The main theorem of \cite{halpernleistner2018structure} identifies two properties of a numerical invariant $\nu$ on $\cX$ that imply that: 1) $\nu$ defines a $\Theta$-stratification of $\cX$; and 2) the substack $\cX^{\rm ss} \subset \cX$ of semistable points is open and admits a separated good moduli space \cite{alper-good-moduli}. The first condition, \emph{strict monotonicity}, is a condition about extending families over points of codimension $2$ (Definition \ref{defn: strictly theta monotone and STR monotone}). Under some mild hypotheses, satisfied in our examples, strict monotonicity implies that the semistable locus admits a separated good moduli space if it is open and bounded \cite[Thm.5.5.10]{halpernleistner2018structure}. The second condition is referred to as \emph{HN-boundedness} (Definition \ref{defn: HN boundedness}), and combined with strict monotonicity it implies the existence of a $\Theta$-stratification, which includes the openness of the semistable locus.\footnote{Technically, in general you only get a weak $\Theta$-stratification if $S$ does not have characteristic $0$.}

The main technical insight of this paper is that strict monotonicity can be guaranteed by choosing a moduli problem that fits a geometric template that we describe below. In particular, this reduces the construction of moduli spaces to the problem of verifying openness and boundedness of the semistable locus.

Let us explain how this works in the example of GIT. First, we regard the setup as an algebraic stack $\cX = Z/G$ along with a morphism $Z/G \to BG$ that is relatively representable by projective schemes. If $Y$ is a regular $2$-dimensional Noetherian scheme and $0 \in Y$ is a closed point, then one can consider whether a given a morphism $\xi : Y \setminus 0 \to \cX$ extends to $Y$. The composition $Y\setminus \{0\} \to BG$ extends uniquely to a morphism $\rho : Y \to BG$. In addition, if $Y$ is equipped with a $\mathbb{G}_m$-action fixing $0$ and the morphism $Y \setminus \{0\} \to BG$ is $\mathbb{G}_m$-equivariant, then the unique extension acquires a unique $\mathbb{G}_m$-equivariant structure. This filling property for any morphism $Y \setminus \{0\} \to BG$ is equivalent to $G$ being reductive \cite[Thm. 1.3]{alper2019cartaniwahorimatsumoto}. 

The morphism $\rho$ does not necessarily lift to $\cX$, but the original morphism $\xi : Y \setminus \{0\} \to \cX$ defines a section of the projective morphism $Y \times_{BG} \cX \to Y$ over $Y \setminus \{0\}$, and we let $\Sigma$ be the closure of this section. By construction, the morphism $\xi$ extends to a morphism $\tilde{\xi} : \Sigma \to \cX$. In other words, given the following diagram of solid arrows, one can always fill in the dotted arrows so that the diagram commutes
\begin{equation} \label{eqn: git filling condition}
\xymatrix{Y \setminus \{0\} \ar@{-->}[r] \ar[d]^\xi \ar@/^10pt/[rr] & \Sigma \ar@{-->}[r] \ar@{-->}[dl]^{\tilde{\xi}} & Y \ar@{-->}[d]^\rho \ar@{-->}[d] \\
\cX \ar[rr] & & BG }.
\end{equation}
If $\xi$ is equivariant with respect to a $\bG_m$-action on $Y$, then this construction equips $\Sigma$ with a canonical $\bG_m$ action such that $\tilde{\xi}$ is $\bG_m$-equivariant.

We postpone the formal definition of strict monotonicity until after our precise definition of a numerical invariant in Subsection \ref{subsection: background theta stratifications}. But roughly, strict monotonicity for the Hilbert-Mumford numerical invariant in GIT follows from two facts: i) if $Y$ has a $\bG_m$-action fixing $0$ and the morphism $\xi$ is $\bG_m$-equivariant, then there exists a $\bG_m$-equivariant filling of the dotted arrows in \eqref{eqn: git filling condition}, as we have established above; and ii) the line bundle $\tilde{\xi}^\ast(L)$ will be relatively ample for the map $\Sigma \to Y$ if $L$ is relatively ample for the map $\cX \to BG$.

The stack $\Coh^d(X)_P$ is not a quotient stack. For simplicity, we will assume in this introduction that $X$ is flat and has geometrically integral fibers of dimension $d$ over the base, and we restrict to the open substack of sheaves of fixed rank $r$. In this case, a family of pure sheaves of dimension $d$ is the same thing as a family of torsion-free sheaves. To verify strict monotonicity of the numerical invariant, we construct a morphism $\Coh^d(X)_{P} \to \Coh^d(X)_{rat}$ to a \emph{non-algebraic} stack $\Coh^d(X)_{rat}$ with the following properties analogous to the properties of the morphism $X/G \to BG$ in GIT:
\begin{enumerate}
    \item Any $\mathbb{G}_m$-equivariant morphism $Y \setminus \{0\} \to \Coh^d(X)_{rat}$ extends to an $\mathbb{G}_m$-equivariant morphism $Y \to \Coh^d(X)_{rat}$.\\
    \item The fibers of the morphism $\Coh^d(X)_{P} \to \Coh^d(X)_{rat}$ are ind-projective ind-schemes, and the duals of the line bundles $L_n$ are asymptotically positive in the sense that on any quasi-compact closed subscheme $L_n^\dual$ is ample for $n \gg 0$.
\end{enumerate}
The term ``infinite dimensional GIT" refers to the infinite dimensional ind-schemes that arise in (2). These conditions imply that any $\bG_m$-equivariant morphism $\xi : Y \setminus \{0\} \to \Coh^d(X)_{P}$ can be extended to a $\bG_m$-equivariant morphism from a birational cover of $Y$, $\tilde{\xi} : \Sigma \to \Coh^d(X)_{P}$, such that $\tilde{\xi}^\ast(L_n^\dual)$ is ample on $\Sigma$ for all $n\gg0$. This allows us to imitate the proof of strict monotonicity in GIT outlined above.

Objects of the groupoid of $S$-points of the stack $\Coh^d(X)_{rat}$ are pairs $(\mathcal{E}, D)$ consisting of a torsion-free sheaf $\cE$ on $X$ and an effective Cartier divisor $D \hookrightarrow X$. A morphism in $\Coh^d(X)_{rat}$ from $(\cE_1,D_1)$ to $(\cE_2,D_2)$ can be defined when $D_2-D_1$ is effective, in which case it is a morphism $\psi : \cE_2 \to \cE_1$ that induces an isomorphism after restriction to $X \setminus D_2$. We think of $\Coh^d(X)_{rat}$ as the stack of vector bundles defined on the complement of a divisor in $X$, i.e., as the stack of ``rational maps" $X \dashrightarrow B\GL_r$. Note that $\Coh^d(X)_{rat}$ is not a stack in groupoids (although all morphisms are monic). For any field $k$ over the base, if one formally inverts all morphisms in $\Coh^d(X)_{rat}(k)$, then the resulting groupoid is canonically equivalent to the groupoid of maps from the generic point of $X_k$ to $B\GL_r$.

Let $T \to \Coh^d(X)_{\rm rat}$ be a morphism corresponding to a $T$-flat family of torsion-free sheaves $\cE$ on $X_T$ along with a $T$-flat family of effective Cartier divisors $D \hookrightarrow X_T$. We define the affine Grassmannian $\Gr_{X,D,\cE}$ to be the fiber product:
\[
\Gr_{X,D,\cE} := \Coh^d(X)_P \times_{\Coh^d(X)_{\rm rat}} T.
\]
For any $T'/T$, a $T'$-point of $\Gr_{X,D,\cE}$ is a $T'$-flat family of torsion-free sheaves $\cF$ on $X_{T'}$ along with a morphism $\cE \to \cF$ that becomes an isomorphism after restriction to $X_{T'} \setminus D_{T'}$. This is a presheaf of sets on the category of $T$-schemes that is representable by an ind-projective ind-scheme over $T$ (Proposition \ref{prop: ind-representability of step 1 grassmannian}). Furthermore, for any quasi-compact closed subscheme $Y \subset \Gr_{X,D,\cE}$, the dual of the line bundle $L_n$ on $\Coh^d(X)_P$, which we use to define semistability, is $T$-ample on $Y$ for $n \gg 0$.

When $X$ is a smooth curve, the ind-schemes $\Gr_{X,D,\cE}$ are known as Beilinson-Drinfeld Grassmannians \cite[\S 5.3.10.(i)]{beilinson-drinfeld} \cite[\S3]{mircovic-vilonen-perverse} \cite[\S3]{zhu_affine_grassmannians}, which are fundamental objects in geometric representation theory. So we regard $\Gr_{X,D,\cE}$ as a higher dimensional analogue of the affine Grassmannian.

\subsection{Comparison with the classical approach}

Since the development of GIT, a standard approach to analyzing a moduli problem that is representable by a locally finite type algebraic stack $\cX$ is to attempt to find a projective quotient stack $Z/G$ with a morphism $Z/G \to \cX$. Hopefully, and often with great effort, one can then show that the restriction to the GIT semistable locus $Z^{\rm ss} / G \to \cX$ is an open immersion whose image admits a simple intrinsic description.

In the example of $\Coh^d(X)_P$, the basic idea in \cite{simpson-repnfundamental-I} is to observe that any sheaf $\mathcal{F}$ admits a surjection $\cO_X(-n)^{\oplus P(n)} \to \mathcal{F}$ for $n$ sufficiently large. The data of $\cF$ along with such a surjection can be regarded as a point in a quot-scheme $Q_n = \Quot(\cO_X(-n)^{\oplus P(n)})$, which is projective. One then studies the action of $\SL_{P(n)}$ on $Q_n$. After a somewhat involved analysis of the Hilbert-Mumford criterion, one finds that the notion of GIT semistability is independent of $n$, and it agrees with the intrinsic notion of Gieseker semistability of the underlying sheaf $\cF$. Finally, one shows that for semistable sheaves the choice of surjection $\cO_X(-n)^{\oplus P(n)} \twoheadrightarrow \mathcal{F}$ is unique up to the action of $\SL_{P(n)}$, and that the GIT quotient $Q_n^{\rm ss} /\!/ \SL_{P(n)}$ is a moduli space for semistable sheaves.

Constructing the stratification of the unstable locus is even more subtle, because unlike in the case of the semistable locus, no single quotient stack will suffice to study all of the strata. In the example of $\Coh^d(X)_P$, the existence of Harder-Narasimhan filtrations and the constructibility of the stratification of $\Coh^d(X)_P$ by Harder-Narasimhan type has been known for some time \cite{shatz77}. However, the canonical structure of locally closed substacks on the strata is a more recent observation \cite{nitsure-scheamtichn, gurjar-nitsure-schematic-hn-lambda}. As in the case of semistability, one can perform a careful analysis of the Hilbert-Mumford criterion to identify each Harder-Narasimhan stratum of $\Coh^d(X)$ with a Hesselink-Kempf-Kirwan-Ness stratum of the quotient stack $Q_n / \SL_{P(n)}$ coming from Geometric Invariant Theory, for $n$ sufficiently large and depending on the stratum \cite{gomez-sols-zamora, hoskins-hn-strat}.

One advantage of the infinite dimensional GIT approach is that it avoids auxiliary choices, and thus avoids the difficult task of showing that the ultimate result is independent of those choices.

Another advantage is the relationship between the analysis for $\Coh^d(X)$ and the analysis for related moduli stacks such as $\cM = \Pair^d(X)$ or $\cM = \Lambda \Coh^d(X)$. In the GIT approach, one starts with the morphism $Q_n / \SL_{P(n)} \to \Coh^d(X)$ described above. The fiber product is then
\[
(Q_n/\SL_{P(n)}) \times_{\Coh^d(X)} \cM \cong Y/\SL_{P(n)}
\]
for some scheme $Y$ that is \emph{affine} over $Q_n$, so one needs a trick in order to apply GIT. For $\cM = \Lambda \Coh^d(X)$, one compactifies $Y$ in some way and then shows that those new points are unstable. This complicates the story, especially because the new points do not lie over points of $\Lambda\Coh^d(X)$. For $\cM = \Pair^d(X)$, one takes the quotient of an open subscheme $Y^\circ \subset Y$ by a free $\bG_m$-action in order to get a scheme that is projective over $Q_n$.

In the infinite dimensional GIT framework, there are analogous stacks of rational objects $\cM_{rat}$ along with a forgetful morphism $\cM_{rat} \to \Coh^d(X)_{rat}$. The fibers of the morphism $\cM \to \cM_{rat}$ over a $T$-point $T \to \cM_{rat}$ are closed sub-ind-schemes of the affine Grassmannians $\Gr_{X,D,\cE}$, where the pair $(D,\cE)$ corresponds to the composition $T \to \cM_{rat} \to \Coh^d(X)_{rat}$ (see Proposition \ref{prop: representability of affine grassmannian for other moduli}). The argument for monotonicity of the numerical invariant on $\Coh^d(X)$ implies monotonicity of the restriction of the numerical invariant to $\cM$ with no further effort.

Although we have chosen to focus on $\Pair^d(X)$ and $\Lambda \Coh^d(X)$ here, there are many more examples of moduli stacks that are affine over $\Coh^d(X)$, such that the relevant affine Grassmannian is a closed sub-ind-scheme of $\Gr_{X,D,\cE}$. Examples include moduli of holomorphic triples \cite{bpp-holomorphic-triples} and more generally moduli spaces related to decorated sheaves \cite{gomezsols.tensors, schmitt-quivers, schmitt-decorated}. In light of our discussion above, the results of this paper reduce the analysis of all of these moduli problems to verifying two types of boundedness conditions: HN boundedness and boundedness of the semistable locus.

One drawback of the infinite dimensional GIT approach is that does not automatically produce a projective moduli space. It is tempting to relate positivity of the line bundle on the affine Grassmannian to ampleness of the line bundle after descending to the good moduli space of the semistable locus. This is possible in certain classical examples, like the moduli of $\SL_n$-bundles on a curve, but at the moment we do not know of any general ampleness results in this framework.

\begin{remark}
The construction of the affine Grassmannians $\Gr_{X,D,\cE}$ makes use of quot-schemes, just like the GIT construction of the moduli space for $\Coh^d(X)^{\rm ss}$. However, the GIT approach uses quot-schemes parameterizing $d$-dimensional quotients, whereas $\Gr_{X,D,\cE}$ is a colimit of quot-schemes parameterizing $(d-1)$-dimensional quotients. Hence, despite this superficial resemblance, the approaches are different. The original construction by Seshadri of the moduli of vector bundles on a curve also uses $(d-1)$-dimensional quot-schemes \cite{seshadri-unitary}, but the relationship with the approach in this paper is not clear.
\end{remark}

\textbf{Acknowledgements:} The first author acknowledges the support of a Simons Collaboration Grant. The second author acknowledges support by NSF grants DMS-1454893 and DMS-2001071. The first and third authors acknowledge the support of the NSF CAREER grant DMS-1945478. We would like to thank Harold Blum, Tom\'as L. G\'omez, Jochen Heinloth, Felix Janda, Yuchen Liu, Chenyang Xu, and Alfonso Zamora for helpful conversations.
\end{section}

\begin{section}{Preliminaries} \label{section: preliminaries}
\begin{subsection}{Notation} \label{subsection: notation}
We work over a fixed quasi-compact base scheme $S$. For any two $S$-schemes $Y$ and $T$, we will always denote $Y \times_S T$ by $Y \times T$. For any $S$-scheme $T$ we let $\text{Aff}_T$ denote the category of (absolutely) affine schemes over $T$. We will write $\text{Sch}_T$ to denote the category of all schemes over $T$. For any map of schemes $Y \rightarrow T$ and any $t \in T$ we write $Y_t$ to denote the fiber $Y \times_T \Spec(k(t))$ over $t$. All of the sheaves we consider will be quasi-coherent. Whenever we write ``$\mathcal{O}_Y$-module" we mean a quasi-coherent $\mathcal{O}_Y$-module. Let $\mathcal{F}$ be an $\mathcal{O}_Y$-module. For any point $t \in T$, we sometimes write $\mathcal{F}_t$ for the restriction of $\mathcal{F}$ to the fiber $Y_t$. This should not be confused with the identical notation employed for the stalk of a sheaf. We will emphasize which use we have in mind whenever it is not clear from context.  

In order to simplify notation for pullbacks, we use the following convention. Whenever we have a map of schemes $f: Y \rightarrow T$ and a sheaf $\mathcal{G}$ on $T$, we will write $\mathcal{G}|_{Y}$ to denote the pullback $f^{*} \mathcal{G}$ whenever the morphism $f$ is clear from context.

We fix once and for all a scheme $X$ that is projective and of finite presentation over $S$. $\pi: X \rightarrow S$ will denote the structure morphism. We also fix an $S$-ample line bundle $\mathcal{O}(1)$ on $X$.

If $\mathcal{M}$ is an algebraic stack, we denote by  $\left\lvert \mathcal{M} \right\rvert$ its underlying topological space \cite[\href{https://stacks.math.columbia.edu/tag/04Y8}{Tag 04Y8}]{stacks-project}. We will often work with rational line bundles in $\Pic(\cM)_{\bQ}: = \Pic(\cM) \otimes_{\bZ} {\bQ}$; we may omit the adjective ``rational" whenever it is clear from context. We also fix some notation for certain stacks defined in \cite{halpernleistner2018structure, alper2019existence}. We will denote by $\Theta$ the quotient stack $[\mathbb{A}^1_{\mathbb{Z}} /\mathbb{G}_m]$. We use the convention that $\mathbb{G}_m$ acts linearly on $\mathbb{A}^1_{\mathbb{Z}} = \Spec(\mathbb{Z}[t])$ so that $t$ is given weight $-1$.

\begin{notn} \label{notation: Y spaces}
Let $R$ be a discrete valuation ring with uniformizer $\varpi$. We define \[Y_{\Theta_{R}} \vcentcolon = \Spec(R[t]) \text{ and } Y_{\overline{ST}_R} := \Spec(R[s,t]/(st-\varpi)),\] equipped with the $\bG_m$-action that assigns $t$ weight $-1$ and $s$ weight $1$. The isomorphism class of the $\mathbb{G}_m$-scheme $Y_{\overline{ST}_{R}}$ is independent of the choice of uniformizer $\varpi$. We denote $\Theta_R = [Y_{\Theta_R} / \bG_m]$ and $\overline{ST}_R = [Y_{\overline{ST}_R}/\bG_m]$. 
\end{notn}
Note that $Y_{\Theta_R}$ and $Y_{\overline{ST}_R}$ each contain a unique $\bG_m$-invariant closed point cut out by the ideals $(t, \varpi)$ and $(s,t)$ respectively. We denote this closed point by $0$ in both cases.

\end{subsection}

\begin{subsection}{Pure sheaves}
We are interested in variations of the moduli stack of pure sheaves on $X$. Here we state some of the relevant definitions and set up some notation.
\begin{defn} \label{defn: absolute torsion-free} 
Let $Y$ be a scheme of finite type over a field. Let $p \in Y$ be a point. We say that $p$ has dimension $d$ if the closure $\overline{p} \subset Y$ is a variety of dimension $d$. 
Let $\mathcal{F}$ be a coherent $\mathcal{O}_Y$-module. We say that $\mathcal{F}$ is a pure sheaf of dimension $d$ if all the associated points of $\mathcal{F}$ have dimension $d$.
\end{defn}
The definition of pure sheaf above is equivalent to the property that $\mathcal{F}$ has support of dimension $d$ and does not contain any nontrivial subsheaves supported on a scheme of dimension smaller than or equal to $d-1$. The latter is the definition of pure sheaf given in \cite[Def. 1.1.2]{huybrechts.lehn}.
 \begin{defn} \label{defn: relative torsion-free} 
Let $\pi: Y \to T$ be a finite type morphism of schemes. We say that a sheaf $\cF$ on $Y$ is $T$-pure of dimension $d$ if it is $T$-flat, finitely presented, and for all $t \in T$ we have that $\mathcal{F}_t$ is a pure sheaf of dimension $d$ on $Y_{t}$.
\end{defn}

\begin{defn}[Moduli of pure sheaves]
The stack $\Coh^{d}(X)$ is the pseudofunctor from $(\text{Sch}/S)^{op}$ to groupoids defined as follows. For every $S$-scheme $T$, we set
\[ \Coh^{d}(X)\, (T) \; = \; \left[ \begin{matrix} \; \; \text{groupoid of $T$-pure  $\mathcal{O}_{X\times T}$-modules of dimension $d$} \; \; \; \end{matrix} \right]\]
\end{defn}
The following proposition is certainly well-known. We include a proof since we are not aware of an explicit reference.
\begin{prop} \label{prop: stack of pure sheaves is algebraic}
$\Coh^{d}(X)$ is an algebraic stack with affine diagonal and locally of finite presentation over $S$.
\end{prop}
\begin{proof}
$\Coh^{d}(X)$ is a subfunctor of $\mathcal{C}\! \mathit{oh}_{X/S}$ as defined in \cite[\href{https://stacks.math.columbia.edu/tag/09DS}{Tag 09DS}]{stacks-project}. By \cite[\href{https://stacks.math.columbia.edu/tag/08KA}{Tag 08WC}]{stacks-project}, $\mathcal{C}\! \mathit{oh}_{X/S}$ is an algebraic stack with affine diagonal. It is locally of finite presentation over $S$ by \cite[\href{https://stacks.math.columbia.edu/tag/08KD}{Tag 08KD}]{stacks-project} and \cite[\href{https://stacks.math.columbia.edu/tag/0CMY}{Tag 0CMY}]{stacks-project}. The proposition follows because the inclusion $\Coh^{d}(X) \hookrightarrow \mathcal{C}\! \mathit{oh}_{X/S}$ is an open immersion by Theorem 12.2.1 (iii) on page 179 of \cite{egaiv}. 
\end{proof} 

\begin{defn} \label{defn: m_n line bundle}
Let $n \in \mathbb{Z}$. Let $\cF_{univ}$ denote the universal sheaf on $\Coh^{d}(X) \times X$, and let $\pi_{\Coh^{d}(X)}: \Coh^{d}(X) \times X \to \Coh^{d}(X)$ denote the first projection. We set $ M_n \vcentcolon = \text{det} \, R\pi_{\Coh^{d}(X) \, *}\left(\mathcal{F}_{univ}(n) \right)$.
\end{defn}
\begin{remark}
The symbol $\text{det}$ denotes the determinant in the $K$-theoretic sense. Note that this makes sense because the derived pushforward is a perfect complex \cite[\href{https://stacks.math.columbia.edu/tag/0A1H}{Tag 0A1H}]{stacks-project}.
\end{remark}

\begin{remark}
Let $T$ be a scheme and let $f: T \to \Coh^{d}(X)$ be a morphism corresponding to a $T$-flat sheaf $\cF$ on $X_T$. Then, by \cite[\href{https://stacks.math.columbia.edu/tag/0A1D}{Tag 0A1D}]{stacks-project} the pullback $f^{*}(M_n)$ can be similarly described as $\text{det} \, R\pi_{T \, *}\left(\mathcal{F}(n) \right)$.
\end{remark}

\begin{defn} \label{defn: knudsen-mumford coefficients}
By \cite[Thm. 4]{knudsen-mumford}, the line bundle $M_n \in {\rm Pic}(\Coh^d(X))$ is a polynomial in the variable $n$ of degree $d+1$ with values in $\Pic(\Coh^d(X))$. More precisely,
\[
M_n = \bigotimes_{i = 0}^{d+1} b_i^{\binom{n}{i}} =  \left(b_{d+1}^{\frac{1}{(d+1)!}}\right)^{n^{d+1}} \otimes \left(b_d^{\frac{1}{d!}} \otimes b_{d+1}^{-\frac{1}{2 \cdot (d-1)!}}\right)^{n^d} \otimes (\text{lower order in }n)
\]
for certain line bundles $b_i$, which can be expressed as $b_i \vcentcolon = \bigotimes_{j = 0}^i M_{j}^{(-1)^{i-j} \binom{i}{j}}$ using the theory of discrete Taylor series.
\end{defn}
\end{subsection}

\begin{subsection}{Hilbert polynomials and Gieseker semistability} \label{subsection: hilbert polynomials background}
The purpose of this subsection is to recall some notions from \cite{huybrechts.lehn} and set some notation in place. Let $Y$ be a scheme that is projective over a field. Fix an ample line bundle $\mathcal{O}(1)$ of $Y$. Let $\mathcal{F}$ be a pure sheaf of dimension $d$ on $Y$. We denote by $P_{\mathcal{F}}(n)$ the Hilbert polynomial of $\mathcal{F}$ \cite[Lemma 1.2.1]{huybrechts.lehn}. This is a polynomial of degree $d$ in the variable $n$. It can be written in the form
\[ P_{\mathcal{F}}(n) = \sum_{k = 0}^{d} \frac{a_k(\mathcal{F})}{k!} n^k\]
for some sequence of rational numbers $a_k(\mathcal{F})$. We set $\rk_{\mathcal{F}} \vcentcolon = a_d(\mathcal{F})$ (notice that this is called the multiplicity in \cite{huybrechts.lehn}, their notion of rank differs from this by a factor of $a_d(\cO_Y)$). It turns out that $\rk_{\mathcal{F}}$ is always a positive integer. We define the reduced Hilbert polynomial to be $\overline{p}_{\mathcal{F}} \vcentcolon = \frac{1}{\rk_{\mathcal{F}}} P_{\mathcal{F}}$. For every $0 \leq i \leq d-1$, we set $\widehat{\mu}_i(\mathcal{F}) \vcentcolon = \frac{a_{i}(\mathcal{F})}{a_d(\mathcal{F})}$. We call this the $i^{th}$ slope of $\mathcal{F}$. The $(d-1)^{th}$ slope $\widehat{\mu}_{d-1}(\mathcal{F})$ is the Mumford slope of $\mathcal{F}$ as in \cite{huybrechts.lehn}; we denote it simply by $\widehat{\mu}(\mathcal{F})$.
\begin{defn}
A pure sheaf $\mathcal{F}$ of dimension $d$ on $Y$ is called Gieseker semistable if for all non-trivial subsheaves $\mathcal{E} \subset \mathcal{F}$ we have $\overline{p}_{\mathcal{E}}(n) \leq \overline{p}_{\mathcal{F}}(n)$ for $n\gg0$.
\end{defn}
For simplicity of notation, we will use the following convention from now on.
\begin{defn} \label{defn: inequality of polynomials}
For any two polynomials $p_1, p_2 \in \mathbb{R}[n]$, we write $p_1 \geq p_2$ if $p_1(n) \geq p_2(n)$ for $n\gg0$.
\end{defn}
In order to check that $\mathcal{F}$ is Gieseker semistable, it suffices to show that $\overline{p}_{\mathcal{E}} \leq \overline{p}_{\mathcal{F}}$ for all nontrivial subsheaves $\mathcal{E}$ such that $\mathcal{F}/ \mathcal{E}$ is also pure of dimension $d$. This is a consequence of \cite[Prop. 1.2.6]{huybrechts.lehn}.

\begin{defn}
For any rational polynomial $P \in \mathbb{Q}[n]$ of degree $d$, let $\Coh^d(X)_P$ denote the subfunctor of $\Coh^d(X)$ consisting of families whose fibers all have Hilbert polynomial $P$.
\end{defn}
Since the Hilbert polynomial is locally constant in flat families of sheaves, $\Coh^{d}(X)$ can be written as a disjoint union of open and closed substacks $\Coh^{d}(X) = \bigsqcup_{P \in \mathbb{Q}[n]} \Coh^{d}(X)_{P}$. %Using this, we can define a new ``normalized" family of line bundles on $\Coh^{d}(X)$.
\begin{defn} \label{defn: l_n line bundle}
Let $n \in \mathbb{Z}$. We define $L_n$ to be the line bundle on $\Coh^{d}(X)$ that is defined on each $\Coh^{d}(X)_{P}$ as follows. We set
\[ L_n|_{\Coh^{d}(X)_{P}} \vcentcolon = M_n|_{\Coh^{d}(X)_{P}} \otimes \left( b_d|_{\Coh^{d}(X)_{P}}\right)^{-\otimes \, \overline{p}(n)} \]
Here the reduced Hilbert polynomial $\overline{p}$ on $\Coh^{d}(X)_{P}$ is defined to be the unique scalar multiple of $P$ with leading coefficient $\frac{1}{d!}$.
\end{defn}
\end{subsection}

\begin{subsection}{Some related moduli problems} \label{section: related moduli problems}
\begin{subsubsection}{Moduli of pairs}
Fix a finitely presented sheaf $\mathcal{A}$ on $X$.
\begin{defn}
The moduli stack of pure pairs $\Pair_{\mathcal{A}}^{d}(X)$ is the pseudofunctor from $(\text{Sch}/S)^{op}$ into groupoids defined as follows. For any $S$-scheme $T$, we set
\[ \Pair_{\mathcal{A}}^{d}(X)\, (T) \; = \; \left[ \begin{matrix} \; \; \text{groupoid of tuples $(\mathcal{F}, \alpha)$, where $\mathcal{F} \in \Coh^{d}(X)(T)$} \; \;  \\ \; \; \; \text{and $\alpha$ is a morphism $\alpha: \mathcal{A}|_{X_T} \rightarrow \mathcal{F}$} \; \; \; \; \end{matrix} \right]\]
We require the isomorphisms in $\Pair_{\mathcal{A}}^{d}(X)$ to be compatible with the morphism $\alpha$.
\end{defn}
There is a natural morphism $\Pair_{\mathcal{A}}^{d}(X) \rightarrow \Coh^{d}(X)$ that forgets the morphism $\alpha$.
\begin{prop} \label{prop: stack pairs algebraic}
The forgetful morphim $\Pair_{\mathcal{A}}^{d}(X) \rightarrow \Coh^{d}(X)$ is schematic, affine and of finite presentation. In particular, $\Pair_{\mathcal{A}}^{d}(X)$ is an algebraic stack with affine diagonal and locally of finite presentation over $S$.
\end{prop}
The following lemma will be used to prove Proposition \ref{prop: stack pairs algebraic}.
\begin{lemma} \label{lemma: representability of hom sheaves}
Let $Y \rightarrow T$ be a morphism of schemes that is proper and of finite presentation. Let $\mathcal{F}$ and $\mathcal{G}$ be quasi-coherent $\mathcal{O}_Y$-modules. Suppose that $\mathcal{G}$ is finitely presented and $T$-flat. Let $\text{Hom}(\mathcal{F}, \mathcal{G})$ denote the functor from $(\text{Sch}/T)^{op}$ into sets that sends a $T$-scheme $L$ to
\[ \text{Hom}(\mathcal{F}, \mathcal{G}) \, (L) \vcentcolon = \text{Hom}_{\mathcal{O}_{Y_L}}(\mathcal{F}|_{Y_L}, \, \mathcal{G}|_{Y_L}) \]
Then, 
\begin{enumerate}[(a)]
    \item $\text{Hom}(\mathcal{F}, \mathcal{G})$ is represented by a scheme that is relatively affine over $T$.
    \item If $\mathcal{F}$ is finitely presented, then $\text{Hom}(\mathcal{F}, \mathcal{G})$ is of finite presentation over $T$.
    \item If $\mathcal{F}$ is of finite type, then the section $0: T \to \text{Hom}(\mathcal{F}, \mathcal{G})$ induced by the $0$ morphism is a closed immersion of finite presentation.
\end{enumerate}
\end{lemma}
\begin{proof}
Parts (a) and (b) are a special case of \cite[\href{https://stacks.math.columbia.edu/tag/08K6}{Tag 08K6}]{stacks-project}.

For (c), suppose first that $\mathcal{F}$ is finitely presented. Then (c) follows because the composition $T \xrightarrow{0} \text{Hom}(\mathcal{F}, \mathcal{G}) \to T$ is clearly of finite presentation and the diagonal $\text{Hom}(\mathcal{F}, \mathcal{G}) \hookrightarrow \text{Hom}(\mathcal{F}, \mathcal{G}) \times_{T} \text{Hom}(\mathcal{F}, \mathcal{G})$ is of finite presentation by \cite[\href{https://stacks.math.columbia.edu/tag/0818}{Tag 0818}]{stacks-project}.

For the general case, since the statement is local on $T$ we can assume without loss of generality that $T$ is affine. Therefore $Y$ is quasi-compact and separated. By \cite[\href{https://stacks.math.columbia.edu/tag/086M}{Tag 086M}]{stacks-project}, $\mathcal{F}$ admits a surjection $\mathcal{E} \twoheadrightarrow \mathcal{F}$ from a finitely presented sheaf $\mathcal{E}$. Observe that the naturally induced diagram
\begin{figure}[H]
\centering
\begin{tikzcd}
  T \ar[d, "id", labels=left] \ar[r, "0"] &  \text{Hom}(\mathcal{F}, \mathcal{G}) \ar[d, symbol = \hookrightarrow] \\   T \ar[r, "0"] &  \text{Hom}(\mathcal{E}, \mathcal{G})
\end{tikzcd}
\end{figure}
\noindent
is Cartesian. Since $\mathcal{E}$ is finitely presented, we already know that $0: T \hookrightarrow \text{Hom}(\mathcal{E}, \mathcal{G})$ is a closed immersion of finite presentation. Therefore, the same holds for $0 : T \hookrightarrow \text{Hom}(\mathcal{F}, \mathcal{G})$.
\end{proof}
\begin{proof}[Proof of Proposition \ref{prop: stack pairs algebraic}]
Let $T$ be an $S$-scheme. Let $T \rightarrow \Coh^{d}(X)$ be a morphism represented by a $T$-pure sheaf $\mathcal{F}$ of dimension $d$ on $X_T$. By definition, the fiber product $\Pair_{\mathcal{A}}^{d}(X) \times_{\Coh^{d}(X)} T$ is the functor $\text{Hom}(\mathcal{A}|_{X_T}, \, \mathcal{F})$ over the scheme $T$. By Lemma \ref{lemma: representability of hom sheaves}, this is represented by a scheme that is relatively affine and of finite presentation over $T$.
\end{proof}
We will abuse notation and denote by $M_n, L_n, b_d$ the pullbacks to $\Pair_{\mathcal{A}}^{d}(X)$ of the corresponding line bundles on $\Coh^{d}(X)$ under the forgetful morphism defined above.
\end{subsubsection}

\begin{subsubsection}{Moduli of $\Lambda$-modules}
In this section we will work with a sheaf $\Lambda$ of rings of differential operators for the morphism $\pi: X \to S$, as in \cite[\S2]{simpson-repnfundamental-I}. We add an extra finite presentation condition in order to work in the non-Noetherian setting.
\begin{defn} \label{defn: rings diff operators}
A sheaf of finitely presented rings of differential operators $\Lambda$ on $X$ relative to $S$ is an associative unital $\mathcal{O}_{X}$-algebra with a filtration $\Lambda_{0} \subset \Lambda_1 \subset ... \subset \Lambda$ satisfying the following conditions.
\begin{enumerate}[(1)]
    \item $\Lambda = \bigcup_{j \geq 0} \Lambda_j$. Furthermore, $\Lambda_i \cdot \Lambda_j \subset \Lambda_{i+j}$ for all $i,j$.
    \item The image of the unit $\mathcal{O}_{X} \rightarrow \Lambda$ is equal to $\Lambda_0$.
    \item The image of $\pi^{-1}(\mathcal{O}_{S})$ in $\mathcal{O}_{X}$ is contained in the center of $\Lambda$.
    \item The left and right $\mathcal{O}_{X}$-module structures on $\Lambda_i /\Lambda_{i-1}$ are equal.
    \item The sheaves of $\mathcal{O}_{X}$-modules $\Lambda_i / \Lambda_{i-1}$ are finitely presented.
    \item The sheaf of graded $\mathcal{O}_{X}$-algebras $\text{Gr}(\Lambda) \vcentcolon = \bigoplus_{j\geq 0} \Lambda_i / \Lambda_{i-1}$ is generated by $\Lambda_1/\Lambda_0$.
    \item Let $T_{\mathcal{O}_{X}}(\Lambda_1/\Lambda_0)$ denote the free associative tensor $\mathcal{O}_{X}$-algebra. The kernel of the natural surjective morphism $T_{\mathcal{O}_{X}}(\Lambda_1/\Lambda_0) \twoheadrightarrow \text{Gr}(\Lambda)$ is locally finitely generated as a two-sided ideal.
\end{enumerate}
\end{defn}
\begin{remark}
The definition of sheaf of rings of differential operators in \cite[\S 2]{simpson-repnfundamental-I} only includes conditions (1) through (6) in the definition above. We include the additional condition (7) to show that the stack $\Lambda \Coh^{d}(X)$ of $\Lambda$-modules is locally of finite presentation over the base $S$ in Proposition \ref{prop: stack of lambda modules is artin}. Without assuming (7), the proof of Proposition \ref{prop: stack of lambda modules is artin} shows that $\Lambda \Coh^d(X)$ is locally of finite type over $S$, and hence locally of finite presentation if $S$ is Noetherian. Therefore when $S$ is Noetherian all of the results that follow hold with condition (7) omitted from Definition \ref{defn: rings diff operators}, which is the context in \cite{simpson-repnfundamental-I}.

If the associated graded algebra $\text{Gr}(\Lambda)$ is commutative and the base scheme $S$ is Noetherian, then condition (7) is automatically satisfied by Hilbert's basis theorem. We do not know if (7) is automatic more generally when the base $S$ is Noetherian.
\end{remark}
The proof of \cite[Lemma 2.2]{simpson-repnfundamental-I} implies that each subsheaf $\Lambda_{i}$ is finitely presented for both the left and right $\mathcal{O}_{X}$-module structures. It follows that $\Lambda$ is quasi-coherent for both the left and right $\mathcal{O}_{X}$-module structures.

For the remainder of this article, we fix a sheaf $\Lambda$ of finitely presented rings of differential operators in Definition \ref{defn: rings diff operators}. 

Let $f: T \rightarrow S$ be a morphism of schemes. By \cite[Lemma 2.5, 2.6]{simpson-repnfundamental-I} there is a canonical isomorphism between the pullbacks $(f_{X}^*)_{l} \Lambda_i$ and $(f_{X}^*)_{r} \Lambda_i$ with respect to the left and right $\mathcal{O}_{X}$-module structures. There is also a canonical isomorphism between the left and right pullbacks of $\Lambda$. We will implicitly use this isomorphism and write $\Lambda_{i}|_{X_{T}}$ and $\Lambda|_{X_{T}}$ to denote either of these pullbacks. It is shown in \cite[Lemma 2.5, 2.6]{simpson-repnfundamental-I} that $\Lambda|_{X_{T}}$ is a sheaf of rings of differential operators on $X_{T}$ relative to $T$. In particular, we have a filtration $\left(\Lambda|_{X_{T}}\right)_0 \subset \left(\Lambda|_{X_{T}}\right)_1 \subset ... \subset \Lambda|_{X_T}$, where $\left(\Lambda|_{X_{T}}\right)_j$ is defined to be the image of the morphism $\Lambda_{j}|_{X_{T}} \rightarrow \Lambda|_{X_{T}}$.
\begin{defn} \label{defn: lambda module}
Let $\mathcal{E}$ be an $S$-flat finitely presented $\mathcal{O}_{X}$-module. A $\Lambda$-module structure on $\mathcal{E}$ is a module structure on $\mathcal{E}$ for the sheaf of rings $\Lambda$ that is compatible with the $\mathcal{O}_{X}$-module structure. If $\mathcal{E}$ is equipped with such structure, then we say that $\mathcal{E}$ is a $\Lambda$-module.
\end{defn}
Let $\mathcal{E}$ be a $\Lambda$-module on $X$. Suppose that we are given a morphism of schemes $f: T \rightarrow S$. \cite[Lemma 2.7]{simpson-repnfundamental-I} implies that the pullback $\mathcal{E}|_{X_{T}}$ can be naturally equipped with a $\Lambda|_{X_{T}}$-module structure. We use this fact to define a stack of $\Lambda$-modules on $X$.
\begin{defn} \label{defn: stack of lambda modules}
The moduli stack $\Lambda\Coh^{d}(X)$ of $\Lambda$-modules is the pseudofunctor from $(\text{Sch}/S)^{op}$ into groupoids defined as follows. For any $S$-scheme $T$, we set
\[ \Lambda\Coh^{d}(X)\, (T) \; = \; \left[ \begin{matrix} \; \; \text{groupoid of $T$-pure sheaves $\mathcal{E}$ of dimension $d$ on $X_{T}$}
\; \; \; \\ \; \; \;  \text{equipped with the structure of a $\Lambda|_{X_{T}}$-module} \; \; \; \end{matrix} \right]\]
The isomorphisms in the groupoid are required to preserve the $\Lambda|_{X_{T}}$-module structure.
\end{defn}

For each $j >1$, let $\mathcal{K}_j$ denote the left and right $\mathcal{O}_{X}$-module of finite type fitting into the following short exact sequence:
\[ 0 \rightarrow \mathcal{K}_{j} \rightarrow (\Lambda_1)^{\otimes j} \rightarrow \Lambda_j \rightarrow 0.\]
Here the tensor product is in the sense of $(\cO_X,\cO_X)$-bimodules, and the morphism $(\Lambda_1)^{\otimes j} \rightarrow \Lambda_j$ is given by multiplication.

\begin{lemma} \label{defn: finitely presented ring diff operators}
There exists some $j_{gen}\gg 0$ such that for all $j \geq j_{gen}$ we have $\cK_{j} = \sum_{a+b=j-j_{gen}} \Lambda_1^{\otimes a} \otimes \cK_{j_{gen}} \otimes \Lambda_1^{\otimes b}$.
\end{lemma}
\begin{proof}
By quasi-compactness, it suffices to show that the kernel of the surjection $\bigcup_{j} \Lambda_1^{\otimes j} \twoheadrightarrow \Lambda$ induced by multiplication is a locally finitely generated two-sided ideal, where the morphisms for the union are given by $\Lambda_1^{\otimes j} \xrightarrow{\text{id}_{\Lambda^{\otimes j}} \otimes 1_{\Lambda}} \Lambda_1^{\otimes j+1 }$. The algebra $\bigcup_{j}\Lambda_1^{\otimes j}$ is naturally filtered; a standard argument reduces to showing that the corresponding initial ideal in $\text{Gr}(\bigcup_{j}\Lambda_1^{\otimes j})$ is locally finitely generated. Note that $\text{Gr}(\bigcup_{j}\Lambda_1^{\otimes j})$ is the free tensor associative $\mathcal{O}_{X}$-algebra $T_{\mathcal{O}_{X}} (\Lambda_1/\Lambda_0)$. 

The initial ideal is the kernel of the associated graded morphism $T_{\mathcal{O}_{X}} (\Lambda_1/\Lambda_0) \to \text{Gr}(\Lambda)$ obtained from our original morphism $\bigcup_{j} \Lambda_1^{\otimes j} \twoheadrightarrow \Lambda$ of filtered algebras. By condition (7) in Definition \ref{defn: rings diff operators}, this kernel is a locally finitely generated two-sided ideal, as desired.
\end{proof}

Let $f: T \rightarrow S$ be an $S$-scheme. Let $\mathcal{E}$ be a finitely-presented $T$-flat $\mathcal{O}_{X_{T}}$-module. Suppose that $\mathcal{E}$ is given the structure of a $\Lambda|_{X_{T}}$-module. Then we get a morphism $a: \Lambda_1|_{X_{T}} \otimes \mathcal{E} \rightarrow \mathcal{E}$ of left $\mathcal{O}_{X_{T}}$-modules given by the composition
\[ \Lambda_1|_{X_{T}} \otimes \mathcal{E} \rightarrow \Lambda|_{X_{T}} \otimes \mathcal{E} \rightarrow \mathcal{E}. \]
By \cite[Lemma 2.8]{simpson-repnfundamental-I}, the image of $\Lambda_{1}|_{X_{T}}$ in $\Lambda|_{X_{T}}$ locally generates $\Lambda|_{X_{T}}$ as an algebra. This implies that the morphism $a$ completely determines the $\Lambda|_{X_{T}}$-module structure on $\mathcal{E}$. The next proposition characterizes which morphisms $a: \Lambda_1|_{X_{T}} \otimes \mathcal{E} \rightarrow \mathcal{E}$ arise from such a $\Lambda|_{X_{T}}$-structure. This gives a useful alternative description of $\Lambda$-modules.
\begin{prop} \label{prop: characterization of lambda modules}
Let $f: T \rightarrow S$ be an $S$-scheme. Let $\mathcal{E}$ be a finitely-presented $T$-flat $\mathcal{O}_{X_{T}}$-module. Let $a: \Lambda_1|_{X_{T}} \otimes \mathcal{E} \rightarrow \mathcal{E}$ be a morphism of left $\mathcal{O}_{X_{T}}$-modules. Then $a$ arises from a $\Lambda|_{X_{T}}$-module structure on $\mathcal{E}$ (as described above) if and only if the following conditions are satisfied
\begin{enumerate}[({$\Lambda$}1)]
    \item The composition $\mathcal{E} \xrightarrow{\sim} \mathcal{O}_{X_{T}} \otimes \mathcal{E} \xrightarrow{unit|_{X_{T}} \otimes \text{id}} \Lambda_1|_{X_{T}} \otimes \mathcal{E} \xrightarrow{a} \mathcal{E}$ is the identity.\\
    
    \item For all $j >1$, the composition $\mathcal{K}_{j}|_{X_{T}} \otimes \mathcal{E} \rightarrow (\Lambda_1|_{X_{T}})^{\otimes j} \otimes \mathcal{E} \rightarrow \mathcal{E}$ is $0$. Here the last morphism $(\Lambda_1|_{X_{T}})^{\otimes j} \otimes \mathcal{E} \rightarrow \mathcal{E}$ is the natural one obtained by applying  $j$-times the morphism $a$, starting with the right-most copy of $\Lambda_1|_{X_{T}}$ and ending at the left-most copy.
\end{enumerate}
In particular, the set of $\Lambda|_{X_{T}}$-module structures on $\mathcal{E}$ is in natural bijection with morphisms $a: \Lambda_1|_{X_{T}} \otimes \mathcal{E} \rightarrow \mathcal{E}$ satisfying conditions $(\Lambda 1)$ and $(\Lambda 2)$ above. If $j_{gen}$ is as Lemma \ref{defn: finitely presented ring diff operators}, then $(\Lambda2)$ is equivalent to the composition $\mathcal{K}_{j_{gen}}|_{X_{T}} \otimes \mathcal{E} \rightarrow (\Lambda_1|_{X_{T}})^{\otimes j_{gen}} \otimes \mathcal{E} \rightarrow \mathcal{E}$ being $0$.
\end{prop}
\begin{proof}
If $a$ comes from a $\Lambda|_{X_{T}}$-module structure, then it follows from construction that it must satisfy conditions $(\Lambda 1)$ and $(\Lambda 2)$ above. $(\Lambda 1)$ follows from the fact that the unit of $\Lambda|_{X_{T}}$ acts as the identity, and $(\Lambda 2)$ follows from the compatibility of the $\Lambda|_{X_{T}}$-module structure with multiplication.

Conversely, suppose that $a: \Lambda_1|_{X_{T}} \otimes \mathcal{E} \rightarrow \mathcal{E}$ is a morphism of left $\mathcal{O}_{X_{T}}$-modules satisfying conditions $(\Lambda 1)$ and $(\Lambda 2)$. By $(\Lambda 2)$, the morphism $(\Lambda_1|_{X_{T}})^{\otimes j} \otimes \mathcal{E} \rightarrow \mathcal{E}$ factors through the quotient $\left((\Lambda_1|_{X_{T}})^{\otimes j}/ \mathcal{K}_j\right) \otimes \mathcal{E} \cong \Lambda_j|_{X_{T}} \otimes \mathcal{E}$. Hence we get a sequence of morphisms $a_j: \Lambda_{j}|_{X_{T}} \otimes \mathcal{E} \rightarrow \mathcal{E}$ compatible with the natural maps $\Lambda_j|_{X_{T}} \rightarrow \Lambda_{j+1}|_{X_{T}}$. This yields a well-defined morphism $\colim_{j\geq 0} \; \Lambda_j|_{X_{T}} \otimes \mathcal{E} \rightarrow \mathcal{E}$. Since tensoring and pulling back commutes with taking colimits, this is equivalent to an action morphism $\Lambda|_{X_{T}} \otimes \mathcal{E} \rightarrow \mathcal{E}$. By construction, this action map is compatible with multiplication on $\Lambda|_{X_{T}}$. Moreover, $(\Lambda 1)$ implies that the unit in $\Lambda|_{X_{T}}$ acts as the identity. We conclude that this is a $\Lambda|_{X_{T}}$-module structure on $\mathcal{E}$. It follows from construction that $a$ is the composition $\Lambda_1|_{X_{T}} \otimes \mathcal{E} \rightarrow \Lambda|_{X_{T}} \otimes \mathcal{E} \rightarrow \mathcal{E}$, as desired. The last statement in the proposition follows because by assumption $\mathcal{K}_{j_{gen}}$ generates all $\mathcal{K}_j$ for $j \geq j_{gen}$.
\end{proof}

There is a natural forgetful morphism of pseudofunctors $\Lambda\Coh^{d}(X) \rightarrow \Coh^{d}(X)$ that forgets the $\Lambda$-module structure.
\begin{prop} \label{prop: stack of lambda modules is artin}
The forgetful morphism $\Lambda\Coh^{d}(X) \rightarrow \Coh^{d}(X)$ is schematic, affine and of finite presentation. In particular, $\Lambda\Coh^{d}(X)$ is an algebraic stack with affine diagonal and locally of finite presentation over $S$.
\end{prop}
\begin{proof}
 Let $T$ be an $S$-scheme. Let $T \rightarrow \Coh^{d}(X)$ be represented by a $T$-pure sheaf $\mathcal{F}$ of dimension $d$ on $X_{T}$. We need to show that $\Lambda\Coh^{d}(X) \times_{\Coh^{d}(X)} T$ is relatively affine and of finite type over $T$.

By Proposition \ref{prop: characterization of lambda modules}, the fiber product $\Lambda\Coh^{d}(X) \times_{\Coh^{d}(X)} T$ is the subfunctor of $\text{Hom}(\Lambda_1|_{X_{T}}\otimes \mathcal{F}, \, \mathcal{F})$ that classifies morphisms $a: \Lambda_1|_{X_{T}} \otimes \mathcal{F} \rightarrow \mathcal{F}$ of left $\mathcal{O}_{X_{T}}$-modules satisfying conditions $(\Lambda 1)$ and $(\Lambda 2)$. By Lemma \ref{lemma: representability of hom sheaves}, the functor $\text{Hom}(\Lambda_1|_{X_{T}}\otimes \mathcal{F}, \, \mathcal{F})$ is represented by a scheme that is relatively affine and of finite presentation over $T$.

Set $H \vcentcolon = \text{Hom}(\Lambda_1|_{X_{T}}\otimes \mathcal{F}, \, \mathcal{F})$. Let $a_{univ}: \Lambda_1|_{X_{H}} \otimes \mathcal{F}|_{X_{H}} \rightarrow \mathcal{F}|_{X_{H}}$ denote the universal morphism on $X_{H}$. The composition 
\[ \mathcal{F}|_{X_{H}} \xrightarrow{\sim} \mathcal{O}_{X_{H}} \otimes \mathcal{F}|_{X_{H}} \xrightarrow{unit|_{X_{H}} \otimes \text{id}} \Lambda_1|_{X_{H}} \otimes \mathcal{F}|_{X_{H}} \xrightarrow{a_{univ}} \mathcal{F}|_{X_{H}}\]
defines a section $f: H \rightarrow \text{Hom}(\mathcal{F},\, \mathcal{F})\times_T H$ of the projection $\text{Hom}(\mathcal{O}_{X_{T}},\, \mathcal{F})\times_T H \rightarrow H$. The constant identity morphism $\text{id}_{\mathcal{F}|_{X_{H}}}$ induces another section. The subfunctor of $H$ where $(\Lambda 1)$ is satisfied is the locus $Z_{0}$ where $f - \text{id}_{\mathcal{F}|_{X_{H}}}$ agrees with the $0$ section. Since the $0$ section is a closed immersion of finite presentation by Lemma \ref{lemma: representability of hom sheaves}, it follows that the locus $Z_0$ is represented by a closed subscheme of finite presentation over $H$.

Let $j_{gen}$ be as in \Cref{defn: finitely presented ring diff operators}. A similar reasoning shows that the locus $Z_{j_{gen}} \hookrightarrow H$ where the following composition vanishes
    \[ \mathcal{K}_{j_{gen}}|_{X_{T}} \otimes \mathcal{E} \rightarrow (\Lambda_1|_{X_{T}})^{\otimes j_{gen}} \otimes \mathcal{E} \rightarrow \mathcal{E}\]
is represented by a closed subscheme of finite presentation. We conclude that $\Lambda\Coh^{d}(X) \times_{\Coh^{d}(X)} T$ is represented by the closed subscheme $Z_0 \cap Z_{j_{gen}}$ of finite presentation in $H$. Since $H$ is relatively affine and of finite presentation over $T$, then so is $\Lambda\Coh^{d}(X) \times_{\Coh^{d}(X)} T$.
\end{proof}
\end{subsubsection}
\end{subsection}
\begin{subsection}{$\Theta$-stratifications and numerical invariants} \label{subsection: background theta stratifications}
We recall some of the theory of $\Theta$-stratifications and numerical invariants introduced in \cite{halpernleistner2018structure}. Let $\mathcal{M}$ be an algebraic stack with quasi-affine diagonal and locally of finite presentation over $S$.

Let $k$ be a field. We denote by $0 \in \mathbb{A}^1_{k}$ the fixed point given by the vanishing of the coordinate $t$ in $\mathbb{A}^1_{k}$. We will abuse notation and also denote by $0$ the corresponding $k$-point of $\Theta_k$. Similarly, we let $1$ be the $k$-point of $\mathbb{A}^1_{k}$ given by the vanishing of $t-1$, and also denote by $1$ its image in the quotient $\Theta_k$. 

For any $k$-valued point $p \in \mathcal{M}(k)$, a filtration of $p$ is a morphism $f: \Theta_k \rightarrow \mathcal{M}$ along with an isomorphism $f(1) \simeq p$ \cite[\S 1]{halpernleistner2018structure}. The stack of filtrations $\text{Filt}(\mathcal{M}) = \Map(\Theta, \, \mathcal{M})$ is represented by an algebraic stack locally of finite presentation over $S$ \cite[Prop. 1.1.2]{halpernleistner2018structure}. There is a morphism $\text{ev}_1: \text{Filt}(\mathcal{M}) \to \mathcal{M}$ given by evaluating at $1 \in \Theta$.

A graded point of the stack $\mathcal{M}$ is a morphism $(B\mathbb{G}_m)_{k} \rightarrow \mathcal{M}$, where $k$ is a field. The mapping stack of graded points $\text{Grad}(\mathcal{M}) \vcentcolon = \text{Map}(B\mathbb{G}_m, \, \mathcal{M})$ is also an algebraic stack locally of finite presentation over $S$ \cite[Prop. 1.1.2]{halpernleistner2018structure}.

Let $\mathcal{X}$ be an open substack of $\mathcal{M}$. A weak $\Theta$-stratum of $\mathcal{X}$ is a union of connected components of $\text{Filt}(\mathcal{X})$ such that the restriction of $\text{ev}_1$ is finite and radicial. A $\Theta$-stratum is a weak $\Theta$-stratum such that $\text{ev}_1$ is a closed immersion. We can think of a $\Theta$-stratum as a closed substack of $\mathcal{X}$ that is identified with some connected components of $\text{Filt}(\mathcal{X})$.
\begin{defn} \label{defn: theta stratification}
A (weak) $\Theta$-stratification of $\mathcal{M}$ consists of a collection of open substacks $(\mathcal{M}_{\leq c})_{c \in \Gamma}$ indexed by a totally ordered set $\Gamma$. We require the following conditions to be satisfied
\begin{enumerate}[(1)]
\item $\mathcal{M}_{\leq c} \subset \mathcal{M}_{\leq c'}$ for all $c< c'$.
\item $\mathcal{M} = \bigcup_{c \in \Gamma} \mathcal{M}_{\leq c}$.
\item For all $c$, there exists a (weak) $\Theta$-stratum $\mathfrak{S}_c \subset \text{Filt}(\mathcal{M}_{\leq c})$ of $\mathcal{M}_{\leq c}$ such that
\[ \mathcal{M}_{\leq c} \setminus \text{ev}_1(\mathfrak{S}_c) = \bigcup_{c' < c} \mathcal{M}_{\leq c'}\]
\item For every point $p \in \mathcal{M}$, the set $\left\{ c \in \Gamma \, \mid \, p \in \mathcal{M}_{\leq c}\right\}$ has a minimal element.
\end{enumerate}
\end{defn}
In order to define $\Theta$-stratifications for our stacks of interest, we will need to introduce some numerical invariants in the sense of \cite[Defn. 0.0.3]{halpernleistner2018structure}. Our numerical invariants will be valued in the ring $\mathbb{R}[n]$ of polynomials in the variable $n$ with coefficients in $\mathbb{R}$. The order described in Definition \ref{defn: inequality of polynomials} equips $\mathbb{R}[n]$ with the structure of a totally ordered $\mathbb{R}$-vector space.

For $q\geq 1$, let $B\mathbb{G}^q_m \vcentcolon = \left[\Spec(\mathbb{Z}) / \mathbb{G}^q_m \right]$. Let $k$ be a field. We say that a morphism $g: B(\mathbb{G}^q_m)_{k} \rightarrow \mathcal{M}$ is nondegenerate if the induced homomorphism $\gamma: (\mathbb{G}_{m}^q)_{k} \rightarrow \text{Aut}(g|_{\Spec(k)})$ has finite kernel. A polynomial numerical invariant is an assignment of a function $\nu_{\gamma}: \mathbb{R}^{q} \setminus \{0\} \to \mathbb{R}[n]$ for each such nondegenerate morphism, satisfying some compatibility conditions.
\begin{defn} \label{defn: numerical poly invariant}
A polynomial numerical invariant $\nu$ on a stack $\mathcal{M}$ is an assignment defined as follows. Let $k$ be a field and let $p \in \mathcal{M}$. Let $\gamma: (\mathbb{G}_m^q)_{k} \rightarrow \text{Aut}(p)$ be a homomorphism of $k$-groups with finite kernel. $\nu$ assigns to this data a scale-invariant function $\nu_{\gamma}: \mathbb{R}^q\setminus \{0\} \rightarrow \mathbb{R}[n]$ such that
\begin{enumerate}[(1)]
    \item $\nu_{\gamma}$ is unchanged under field extensions $k \subset k'$.
    \item $\nu$ is locally constant in algebraic families. In other words, let $T$ be a scheme. Let $\xi: T \rightarrow \mathcal{M}$ be a morphism and let $\gamma: (\mathbb{G}_{m}^{q})_{T} \rightarrow \text{Aut}(\xi)$ be a homomorphism of $T$-group schemes with finite kernel. Then as we vary $t \in T$, the function $\nu_{\gamma_{t}}$ is locally constant in $T$.
    \item Given a homomorphism $\phi: (\mathbb{G}^w_{m})_k \rightarrow (\mathbb{G}^q_{m})_{k}$ with finite kernel, the function $\nu_{\gamma \circ \phi}$ is the restriction of $\nu_{\gamma}$ along the inclusion $\mathbb{R}^w \hookrightarrow \mathbb{R}^q$ induced by $\phi$.
\end{enumerate}
\end{defn}

We say that a filtration $f: \Theta_{k} \rightarrow \mathcal{M}$ is nondegenerate if the restriction $f|_0: [0/\, (\mathbb{G}_m)_{k}] \to \mathcal{M}$ is nondegenerate. We regard $\nu$ as a function on the set of nondegenerate filtrations by defining $\nu(f) \vcentcolon = \nu_{f|_{0}}(1) \in \bR[n]$. Given a point $p \in |\mathcal{M}|$, we say that $p$ is \emph{semistable} if all nondegenerate filtrations $f$ with $f(1)=p$ satisfy $\nu(f) \leq 0$. Otherwise we say that $p$ is \emph{unstable}. Note that although \Cref{defn: numerical poly invariant} involves data for all $q \geq 1$, only the $q=1$ data is used to define semistability.

\begin{remark}
We will only consider nondegenerate filtrations, because these are the ones relevant for stability. For the rest of the article, we will sometimes omit the adjective ``nondegenerate".
\end{remark}

Next we explain a useful way to construct polynomial numerical invariants on a stack. Fix a sequence of rational line bundles $(\mathcal{L}_n)_{n \in \mathbb{Z}}$, where each $\mathcal{L}_n$ is in the rational Picard group $Pic(\mathcal{M})\otimes_{\mathbb{Z}} \mathbb{Q}$ of the stack. For each morphism $g: (B\mathbb{G}^q_m)_{k} \rightarrow \mathcal{M}$, the pullback line bundle $g^*(\mathcal{L}_{n})$ amounts to a rational character in $X^*(\mathbb{G}^q_m)\otimes_{\mathbb{Z}} \mathbb{Q}$. Under the natural identification $X^*(\mathbb{G}^q_m)\otimes_{\mathbb{Z}} \mathbb{Q} \cong \mathbb{Q}^q$, we can interpret this as a $q$-tuple of rational numbers $(w_n^{(i)})_{i=1}^q$, which we call the weight of $g^*(\mathcal{L}_n)$. It is often the case that one can choose the line bundles $(\mathcal{L}_n)_{n \in \mathbb{Z}}$ in such a way that for each fixed $i$ the weight $w_n^{(i)}$ is a polynomial in $\mathbb{Q}[n]$. If this is the case, then for every $g$ we can define an $\mathbb{R}$-linear function $L_{g} : \mathbb{R}^q \to \mathbb{R}[n]$ given by
\[ L_{g}\left((r_i)_{i=1}^q\right) = \sum_{i=1}^q r_i \cdot w_n^{(i)} \]
In order to obtain a scale invariant $\nu$, we use a rational quadratic norm on graded points as in \cite[Defn. 4.1.12]{halpernleistner2018structure}. This consists of an assignment of a positive definite quadratic norm $b_{\gamma}(-)$ with rational coefficients defined on $\mathbb{R}^q$ for each choice of $p \in \mathcal{M}(k)$ and homomorphism $\gamma: (\mathbb{G}_m^q)_{k} \rightarrow \text{Aut}(p)$ with finite kernel. Just as in the definition of numerical invariants, we require that this assignment is stable under field extension and locally constant on algebraic families. Also, we require that for any homomorphism $\phi: (\mathbb{G}^w_{m})_k \rightarrow (\mathbb{G}^q_{m})_{k}$ with finite kernel the quadratic norms $b_{\gamma}$ and $b_{\gamma \circ \phi}$ are compatible with the corresponding induced inclusion $\mathbb{R}^m \hookrightarrow \mathbb{R}^q$.

Given a sequence of rational line bundles $(\mathcal{L}_n)_{n \in \mathbb{Z}}$ as above and a rational quadratic norm on graded points $b$, we can define a numerical invariant $\nu$ as follows. For all nondegenerate $g: (B\mathbb{G}^q_m)_{k} \rightarrow \mathcal{M}$ with corresponding morphism $\gamma: (\mathbb{G}^q_m)_{k} \to \text{Aut}(g|_{\Spec(k)})$, we set
\[ \nu_{\gamma}(\vec{r}) = \frac{L_{g}(\vec{r})}{\sqrt{b_{\gamma}(\vec{r})}}  \]

We next explain how a polynomial numerical invariant $\nu$ can define a (weak) $\Theta$-stratification on $\mathcal{M}$ (see \cite[\S 4.1]{halpernleistner2018structure} for more details). For any unstable point, we set $M^{\nu}(p)$ to be the supremum of $\nu(f)$ over all filtrations $f$ with $f(1)=p$ (if such supremum exists). If $p$ is semistable, then by convention we set $M^{\nu}(p) = 0$. For any $c \in \mathbb{R}[n]_{\geq 0}$, we set $\mathcal{M}_{\leq c}$ to be the set of all points $p$ satisfying $M^{\nu}(p) \leq c$, and we let $\cM^{\nu\dash\rm{ss}} := \cM_{\leq 0}$ denote the set of semistable points. For any unstable point $p$, a filtration $f$ of $p$ is called a Harder-Narasimhan filtration if $\nu(f) = M^{\nu}(p)$.

We say that $\nu$ defines a (weak) $\Theta$-stratification if: 1) every unstable point has a Harder-Narasimhan filtration that is unique up to pre-composing with a ramified covering $\Theta \to \Theta$; and 2) $\cM_{\leq c}$ are open substacks of $\cM$ coming from a (weak) $\Theta$-stratification such that the (weak) stratum $\mathfrak{S}_c \subset \Filt(\cM_{\leq c})$ is an open and closed substack of Harder-Narasimhan filtrations $f$ with $\nu(f)=c$. If $\nu$ defines a $\Theta$-stratification, then the Harder-Narasimhan filtration of any point is defined over the field of definition of that point, but if it is only a weak $\Theta$-stratification, the Harder-Narsimhan filtration might only be defined over a finite purely inseparable field extension.

A natural question to ask is: When does a numerical invariant $\nu$ define a (weak) $\Theta$-stratification as described above, and when does $\mathcal{M}^{\nu \dash \mathrm{ss}}$ admit a good moduli space in the sense of \cite{alper-good-moduli}? The following theorem provides sufficient criteria.
\begin{thm}[{\cite[Theorem B]{halpernleistner2018structure}}] \label{thm: theta stability paper theorem}
Let $\nu$ be a polynomial numerical invariant on $\mathcal{M}$ defined by a sequence of rational line bundles and a norm on graded points, as explained above.
\begin{enumerate}
    \item If $\nu$ is strictly $\Theta$-monotone, then it defines a weak $\Theta$-stratification of $\mathcal{M}$ if and only if it satisfies the HN boundedness condition. If moreover $S$ is defined over $\mathbb{Q}$, then $\nu$ defines a $\Theta$-stratification.
    \item Suppose that all of the conditions in (1) above are satisfied and $S$ is defined over $\mathbb{Q}$. Furthermore, assume $\nu$ is strictly $S$-monotone and that the semistable locus $\mathcal{M}^{\nu \dash \mathrm{ss}}$ can be written as a disjoint union of bounded open substacks. Then $\mathcal{M}^{\nu \dash \mathrm{ss}}$ has a separated good moduli space.
\end{enumerate}
If all of the above are satisfied and if $\mathcal{M}$ satisfies the existence part of the valuative criterion for properness for complete discrete valuation rings relative to $S$ \cite[\href{https://stacks.math.columbia.edu/tag/0CLK}{Tag 0CLK}]{stacks-project}, then each quasi-compact open and closed substack of $\mathcal{M}^{\nu \dash \mathrm{ss}}$ admits a good moduli space that is proper over $S$.
\end{thm}

We end this subsection by explaining each of the hypotheses that one needs to check in Theorem \ref{thm: theta stability paper theorem}. Let $\kappa$ be a field and let $a \geq 1$ be an integer. We denote by $\mathbb{P}^1_{\kappa}[a]$ the $\mathbb{G}_m$-scheme $\mathbb{P}^{1}_{\kappa}$ equipped with the $\mathbb{G}_m$-action determined by the equation $t \cdot [x:y] = [t^{-a}x : y]$. We set $0 = [0: 1]$ and $\infty = [1:0]$. Recall the schemes $Y_{\Theta_R}$ and $Y_{\overline{ST}_R}$ of Notation \ref{notation: Y spaces}.

We will use simplified versions of the monotonicity conditions in \cite[\S 5]{halpernleistner2018structure}. We refer the reader there for the general definitions, which allow $\Sigma$ to be an orbifold.
\begin{defn} \label{defn: strictly theta monotone and STR monotone} A polynomial numerical invariant $\nu$ on $\mathcal{M}$ is strictly $\Theta$-monotone (resp. strictly $S$-monotone) if the following condition holds. 

Let $R$ be any complete discrete valuation ring and set $\mathfrak{X}$ to be $\Theta_{R}$ (resp. $\overline{ST}_{R}$). Choose a map $\varphi: \mathfrak{X} \setminus 0 \rightarrow \mathcal{M}$. Then, after maybe replacing $R$ with a finite DVR extension, there exists a reduced and irreducible $\mathbb{G}_m$-equivariant scheme $\Sigma$ with maps $f: \Sigma \rightarrow Y_{\mathfrak{X}}$ and $\widetilde{\varphi}: \left[\Sigma/ \, \mathbb{G}_m\right] \rightarrow  \mathcal{M}$ such that
\begin{enumerate}[({M}1)]
    \item The map $f$ is proper, $\mathbb{G}_m$-equivariant, and its restriction induces an isomorphism $f : \, \Sigma_{Y_{\mathfrak{X}} \setminus 0} \xrightarrow{\sim} Y_{\mathfrak{X}} \setminus 0$.
    \item The following diagram commutes
\begin{figure}[H]
\centering
\begin{tikzcd}
  \left[\left(\Sigma_{Y_{\mathfrak{X}} \setminus 0}\right)/ \, \mathbb{G}_m \right] \ar[rd, "\widetilde{\varphi}"] \ar[d, "f"'] & \\   \mathfrak{X} \setminus 0 \ar[r, "\varphi"'] &  \mathcal{M}
\end{tikzcd}
\end{figure}
    \item Let $\kappa$ denote a finite extension of the residue field of $R$. For any $a \geq 1$ and any finite $\mathbb{G}_m$-equivariant morphism $\mathbb{P}^1_{\kappa}[a] \to \Sigma_{0}$, we have $\nu\left( \;\widetilde{\varphi}|_{\left[\infty / \mathbb{G}_m\right]} \;\right) >  \nu\left(\; \widetilde{\varphi}|_{\left[0 / \mathbb{G}_m\right]} \;\right)$.
\end{enumerate}
\end{defn}

\begin{defn}[HN Boundedness] \label{defn: HN boundedness}
We say that a polynomial numerical invariant $\nu$ satisfies the HN boundedness condition if the following is always satisfied:

Let $T$ be an affine Noetherian scheme. Choose a morphism $g: T \rightarrow \mathcal{M}$. Then there exists a quasi-compact open substack $\mathcal{U}_{T} \subset \mathcal{M}$ such that the following holds.
For all geometric points $t \in T$ with residue field $k$ and all nondegenerate filtrations $f: \Theta_{k} \rightarrow \mathcal{M}$ of the point $g(t)$ with $\nu(f)>0$, there exists another filtration $f'$ of $g(t)$ satisfying $\nu(f') \geq \nu(f)$ and $f'|_{0} \in \mathcal{U}_{T}$.
\end{defn}
In plain words, this says that for the purposes of maximizing $\nu(f)$ among all filtrations of points in a bounded family, it suffices to consider only $f$ such that the associated graded $f|_{0}$ lies in some other (possibly larger) bounded family.
\end{subsection}
\end{section}
\begin{section}{Rational maps and affine Grassmannians}  \label{section: affine grassmannian pure sheaves}
In this section we define a stack of rational maps $\mathcal{M}_{rat}$ associated to each of the moduli stacks $\mathcal{M}$ we are considering. There is a natural morphism $\mathcal{M} \to \mathcal{M}_{rat}$ whose fibers are higher dimensional analogues of affine Grassmannians. Just as in the classical case of vector bundles on a curve, our affine Grassmannians are ind-projective strict-ind-schemes over $S$. 

We will also identify a family of line bundles on the affine Grassmannian that are eventually relatively ample on each projective stratum. The material from this section will be one of the ingredients in the infinite dimensional GIT argument we will use later on to prove monotonicity of the numerical invariant on $\cM$.
\begin{subsection}{Regular principal subschemes and twists}
\begin{defn}
Let $Y$ be a scheme. We say that a closed subscheme $D\hookrightarrow Y$ is principal if the corresponding ideal sheaf $\mathcal{I}_D$ is locally a principal ideal.
\end{defn}

\begin{example}
We give some examples of principal subschemes to illustrate that not all principal subschemes are Cartier divisors. First, the identity subscheme $Y \subset Y$ is always principal, cut out by the $0$ function. For another example, if we let $Y:= \Spec(k[x,y]/(xy))$, then the subscheme $D\subset Y$ cut out by the function $x$ is a principal subscheme. 
\end{example}

\begin{remark} \label{remark: being principal is fpqc local}
A subscheme $D$ is principal if and only if the first Fitting ideal $\text{Fit}_1(\mathcal{I}_{D})$ is the structure sheaf $\mathcal{O}_{Y}$. Indeed, if $\mathcal{I}_{D}$ is locally generated by one element, then $\text{Fit}_1(\mathcal{I}_{D}) = \mathcal{O}_{Y}$ by \cite[\href{https://stacks.math.columbia.edu/tag/07ZA}{Tag 07ZA}(1)]{stacks-project}. The converse follows by \cite[\href{https://stacks.math.columbia.edu/tag/07ZC}{Tag 07ZC}]{stacks-project}. Since the formation of Fitting ideals of sheaves commutes with base-change \cite[\href{https://stacks.math.columbia.edu/tag/07ZA}{Tag 07ZA}(3)]{stacks-project}, and $f^*(\mathcal{I}_D) = \mathcal{I}_{f^{-1}(D)}$ for a flat morphism $f : Y' \to Y$, it follows that the property of being a principal subscheme can be checked flat locally.
\end{remark}
Let $D_1$ and $D_2$ be two principal subschemes cut out by the $\mathcal{O}_Y$-ideals $\mathcal{I}_{D_1}$ and $\mathcal{I}_{D_2}$. We define the sum $D_1 + D_2$ to be the locally principal scheme cut out by the ideal $\mathcal{I}_{D_1} \cdot \mathcal{I}_{D_2}$.

\begin{defn}
Let $T$ be an $S$-scheme. Let $D \hookrightarrow X_T$ be a principal subscheme. Let $\mathcal{F}$ be a $T$-flat finitely presented $\mathcal{O}_{X_T}$-module. We say that $D$ is $\mathcal{F}$-regular if for all $t \in T$ the fiber $D_t$ does not contain any associated point of $\mathcal{F}|_{X_t}$.
\end{defn}

The following lemma collects some useful properties of the notion of $\mathcal{F}$-regular subscheme.

\begin{lemma} \label{lemma: properties of regularity}
Let $C, D \hookrightarrow X$ be two principal subschemes of the finitely presented scheme $X \to S$. Let $\mathcal{F}$ be an $S$-flat finitely presented $\mathcal{O}_{X}$-module. Suppose that both $C,D$ are $\mathcal{F}$-regular. Then,
\begin{enumerate}[(a)]
    \item $C+D$ is $\mathcal{F}$-regular.
    \item For any morphism of schemes $f:T \rightarrow S$, the preimage $(f_X)^{-1}(D)$ under the base-change morphism $f_X : X_T \to X_S$ is $\mathcal{F}|_{X_T}$-regular.
    \item Let $\sigma: \mathcal{I}_{D} \rightarrow \mathcal{O}_{X}$ denote the inclusion of the ideal sheaf. Then, the morphism $\sigma \otimes \text{id}_{\mathcal{F}}: \mathcal{I}_{D} \otimes \mathcal{F} \rightarrow \mathcal{F}$ is injective. Moreover, the quotient $\mathcal{F} / \, (\mathcal{I}_{D} \otimes \mathcal{F})$ is $S$-flat.
    \item $\mathcal{I}_{D} \otimes \mathcal{F}$ is locally isomorphic to $\mathcal{F}$. 
\end{enumerate}
\end{lemma}
\begin{proof}
\quad
\begin{enumerate}[(a)]
    \item For each $s \in S$, the support of $(C+D)_s$ is the union of the supports of $C_s$ and $D_s$. If $C_s$ and $D_s$ do not contain any associated point of $\mathcal{F}|_{X_s}$, then neither does the union.
    \item Let $t \in T$, and set $s = f(t)$. By \cite[\href{https://stacks.math.columbia.edu/tag/05DC}{Tag 05DC}]{stacks-project}, the set of associated points of $\mathcal{F}|_{X_{t}}$ are contained in the preimage of the associated points of $\mathcal{F}|_{X_s}$ under the morphism $f_{t}: X_{t} \to X_{s}$. This implies the claim, because the fiber $\left((f_X)^{-1}D\right)_{t}$ is the preimage under $f_{t}$ of the fiber $D_{s}$.
    \item Since the question is local, we can assume that $X$ is affine and the ideal $\mathcal{I}_{D}$ is generated by a global section $\mathcal{O}_{X} \xrightarrow{u} \mathcal{O}_{X}$. Let $\text{Ann}(u) \subset \mathcal{O}_{X}$ denote the ideal annihilating $u$. By definition, the sequence $ \text{Ann}(u) \hookrightarrow \mathcal{O}_{X} \xrightarrow{u} \mathcal{O}_{X}$ is exact in the middle, and the cokernel of the first map is identified with $\mathcal{I}_{D}$. Consider the complex $ \text{Ann}(u) \otimes \mathcal{F} \to \mathcal{F} \xrightarrow{u} \mathcal{F}$ obtained by tensoring with $\mathcal{F}$. It suffices to show that the second morphism $u: \mathcal{F} \to \mathcal{F}$ is injective and that $\mathcal{F}/ \, (u\cdot \mathcal{F})$ is $S$-flat.
    
    Since $D$ is $\mathcal{F}$-regular, it follows that for all $s \in S$ the restriction $\mathcal{F}|_{X_s} \xrightarrow{u|_{X_s}} \mathcal{F}|_{X_s}$ is injective. Because $\mathcal{F}$ is $S$-flat, the slicing criterion for flatness \cite[\href{https://stacks.math.columbia.edu/tag/00ME}{Tag 00ME}]{stacks-project} implies that $\mathcal{F} \xrightarrow{u} \mathcal{F}$ is injective and $\mathcal{F}/ \, (u\cdot \mathcal{F})$ is $S$-flat.
    
    \item Locally, we can choose $u$ as in part (c). We have $\mathcal{F} \otimes \mathcal{I}_{D} \cong \text{coker}\left[ \text{Ann}(u) \otimes \mathcal{F} \to \mathcal{F} \right]$. It follows from part (c) that $\text{Ann}(u) \otimes \mathcal{F} \to \mathcal{F}$ is the $0$ morphism, since the image is contained in the kernel of the injection $\mathcal{F} \xrightarrow{u} \mathcal{F}$. Therefore $\mathcal{F} \otimes \mathcal{I}_{D} \cong \mathcal{F}$.
\end{enumerate}
\end{proof}
Let $D \hookrightarrow X$ be an $\mathcal{F}$-regular principal subscheme. Let $U = X \setminus D$ be the open complement. The open immersion $j : U \hookrightarrow X$ is locally given by the localization of a generator of the ideal sheaf $\mathcal{I}_{D}$ that cuts out $D$. In particular $j$ is affine. We can use this local interpretation of $j$ plus Lemma \ref{lemma: properties of regularity}(c) to see that the unit $\mathcal{F} \rightarrow j_* j^* \mathcal{F}$ is injective. Consider the isomorphism $\mathcal{I}_{D} \otimes j_* j^* \mathcal{F} \rightarrow j_*j^*\mathcal{F}$ induced by multiplication.
\begin{defn}[Twists]
With notation as above, we define $\mathcal{F}(-D)$ to be the image of $\mathcal{I}_{D} \otimes \mathcal{F}$ under the multiplication morphism $\mathcal{I}_{D} \otimes j_* j^* \mathcal{F} \rightarrow j_*j^*\mathcal{F}$.

Similarly, we define $\mathcal{F}(D)$ to be the maximal subsheaf $\cE \subset j_* j^* \mathcal{F}$ such that the image of $\mathcal{I}_{D} \otimes \cE \to j_* j^* \mathcal{F}$ lands in $\mathcal{F} \subset j_* j^* \mathcal{F}$.
\end{defn}
By definition, there is an infinite chain of injections
\[ ... \hookrightarrow \mathcal{F}(-nD) \hookrightarrow ... \hookrightarrow \mathcal{F}(-2D) \hookrightarrow \mathcal{F}(-D) \hookrightarrow \mathcal{F} \hookrightarrow \mathcal{F}(D) \hookrightarrow ... \hookrightarrow \mathcal{F}(nD) \hookrightarrow ...\]
By using the affine local interpretation of $j$ as a localization, it can be seen that there is a canonical identification
$j_{*} \, j^* \, \mathcal{F} \cong \underset{m \in \mathbb{Z}}{\colim} \; \mathcal{F}(mD)$.

Lemma \ref{lemma: properties of regularity} (d) implies that $\mathcal{F}(nD)$ is locally isomorphic to $\mathcal{F}$ for all $n \in \mathbb{Z}$. In particular, if $\mathcal{F}$ is $S$-pure of dimension $d$, then all the twists $\mathcal{F}(nD)$ are $S$-pure of dimension $d$.
\end{subsection}

\begin{subsection}{Stacks of rational maps}
\begin{defn} \label{defn: torsion-free rational maps}
We denote by $\Coh^{d}(X)_{\text{rat}}$ the pseudofunctor from $(\text{Aff}_S)^{op}$ to categories defined as follows. For each affine scheme $T$ in $\text{Aff}_S$, the objects of $\Coh^{d}(X)_{\text{rat}}\,(T)$ are pairs $(D, \mathcal{E})$, where 
\begin{enumerate}[(1)]
    \item $\mathcal{E}$ is a $T$-pure sheaf of dimension $d$ on $X_T$.
    \item $D \hookrightarrow X_{T}$ is an $\mathcal{E}$-regular principal subscheme of $X_{T}$.
\end{enumerate}
Let $A= (D_1, \mathcal{E}_1)$ and $B = (D_2, \mathcal{E}_2)$ be two objects in $\Coh^{d}(X)_{\text{rat}}\,(T)$. The set of morphisms $\text{Mor}_{\Coh^{d}(X)_{\text{rat}}\,(T)}(A, B)$ consists of pairs $(i, \psi)$, where
\begin{enumerate}[(1)]
    \item $i: D_1 \hookrightarrow D_2$ is an inclusion.
    \item $\psi: \mathcal{E}_2 \rightarrow \mathcal{E}_{1}$ is a monomorphism such that the restriction $\psi|_{X \setminus D_2}$ is an isomorphism.
\end{enumerate}
\end{defn}

Note that all morphisms in $\Coh^{d}(X)_{\text{rat}}\,(T)$ are monic.
 
\begin{remark}
In order to see that this yields a well-defined functor, it is necessary to check that $\psi$ remains a monomorphism after base-changing $T$. This follows because the cokernel of $\psi$ is $T$-flat, by an argument similar to the one for (i) inside the proof of Lemma \ref{lemma: representability truncated version grassmannian} below.
\end{remark}

The pseudofunctor $\Coh^{d}(X)_{\text{rat}}$ is a stack in the fpqc topology. Indeed, descent data amounts to a cocycle with isomorphisms in $\Coh^{d}(X)_{\text{rat}}$. By definition these are the usual naive isomorphisms of pairs $(D, \mathcal{F})$. Therefore, descent follows from the theory of fpqc descent for quasi-coherent sheaves and the fact that the property of being a principal subscheme is fpqc local (Remark \ref{remark: being principal is fpqc local}).

We conclude that if $Y$ is an affine scheme in $\text{Aff}_{S}$ and $H$ is an $S$-flat algebraic group scheme acting on $Y$, then the data of a morphism of pseudofunctors $\left[Y/H \right] \rightarrow \Coh^{d}(X)_{\text{rat}}$ amounts to a pair $(D, \mathcal{F})$, where
\begin{enumerate}[(1)]
    \item $\mathcal{F}$ is an $H$-equivariant $Y$-pure sheaf of dimension $d$ on $X_Y$;
    \item $D \hookrightarrow X_{Y}$ is an $H$-stable $\mathcal{F}$-regular principal subscheme.
\end{enumerate}

\begin{remark} \label{rem: variants of coh rat}
It is possible to develop the theory of infinite dimensional GIT using variants of this definition of $\Coh^d(X)_{\rm{rat}}$. For instance, one could define objects to be those of $\Coh^d(X)$, and morphisms to be equivalence classes consisting of pairs $(D \subset X, \psi : \cE_1|_{X \setminus D} \to \cE_{2}|_{X \setminus D})$, where $D$ is a principal subscheme that is $\cE_1$-regular and $\cE_2$-regular, and $\psi$ is an isomorphism. Two such morphisms are equivalent if they agree on the intersection of the open subsets where they are defined. This is a pseudofunctor valued in groupoids, but it does not satisfy smooth descent. The rational filling condition in Definition \ref{defn: rational filling} would be more difficult to state without smooth descent.
\end{remark}

We can similarly define stacks of rational maps for each of the other stacks introduced in \ref{section: related moduli problems}.
\begin{defn} \label{defn: rational maps pairs}
We denote by $\Pair_{\mathcal{A}}^{d}(X)_{\text{rat}}$ the pseudofunctor from $(\text{Aff}_S)^{op}$ to categories defined as follows. For each affine scheme $T$ in $\text{Aff}_S$, the objects of $\Pair_{\mathcal{A}}^{d}(X)_{\text{rat}}\,(T)$ consist of triples $(D,\mathcal{E}, \alpha)$, where 
\begin{enumerate}[(1)]
    \item $(D, \mathcal{E})$ is an element of $\Coh^{d}(X)_{rat}$.
    \item $\alpha$ is a morphism $\alpha: \mathcal{A}|_{X_T \setminus D} \rightarrow \mathcal{F}|_{X_T \setminus D}$.
\end{enumerate}
A morphism from $(D_1, \mathcal{E}_1, \alpha_1)$ to $(D_2, \mathcal{E}_2, \alpha_2)$ is a morphism $(i,\psi) : (D_1, \mathcal{E}_1) \to (D_2, \mathcal{E}_2)$ such that $\psi|_{X_{T} \setminus D_2}$ is compatible with the sections $\alpha_1|_{X_{T} \setminus D_2}$ and $\alpha_2$.
\end{defn}

\begin{defn} \label{defn: rational maps lambda modules}
We denote by $\Lambda\Coh^{d}(X)_{\text{rat}}$ the pseudofunctor from $(\text{Aff}_S)^{op}$ to categories defined as follows. For each affine scheme $T$ in $\text{Aff}_S$, the objects of $\Lambda\Coh^{d}(X)_{\text{rat}}\,(T)$ consist of triples $(D,\mathcal{E}, a)$, where 
\begin{enumerate}[(1)]
    \item $(D, \mathcal{E})$ is an element of $\Coh^{d}(X)_{rat}$.
    \item $a$ is a morphism $a: \Lambda_1|_{X_{T} \setminus D}\otimes \mathcal{E}|_{X_{T} \setminus D} \rightarrow \mathcal{E}|_{X_{T} \setminus D}$ of left $\mathcal{O}_{X_{T} \setminus D}$-modules that endows $\mathcal{E}|_{X_{T} \setminus D}$ with the structure of a $\Lambda|_{X_{T} \setminus D}$-module (as in Proposition \ref{prop: characterization of lambda modules}).
\end{enumerate}
A morphism from $(D_1, \mathcal{E}_1, a_1)$ to $(D_2, \mathcal{E}_2, a_2)$ is a morphism $(i,\psi) : (D_1, \mathcal{E}_1) \to (D_2, \mathcal{E}_2)$ such that $\psi|_{X_{T} \setminus D_2}$ is compatible with the action morphisms $a_1|_{X_{T} \setminus D_2}$ and $a_2$.
\end{defn}
Let $\mathcal{M}= \Pair_{\mathcal{A}}^{d}(X)$ or $\Lambda\Coh^{d}(X)$. The theory of fpqc descent for quasi-coherent sheaves implies that $\mathcal{M}_{\text{rat}}$ satisfies fpqc descent. This shows that the data of a pseudofunctor $\left[ Y / H \right] \rightarrow \mathcal{M}_{\text{rat}}$ from a quotient stack $\left[ Y / H \right]$ can be concretely described as an $H$-equivariant pair $(D, \mathcal{F}) \in \Coh^{d}(X)_{rat}(Y)$, and the corresponding $H$-equivariant structure defined on $X_{T} \setminus D$ in each case.
\end{subsection}

\begin{subsection}{The affine Grassmannian for pure sheaves}
There is a morphism $\Coh^{d}(X) \rightarrow \Coh^{d}(X)_{\text{rat}}$ that takes a pure sheaf $\mathcal{E}$ on $X_T$ and maps it to the pair $(\emptyset,  \mathcal{E})$. In this subsection we will describe the ``fibers" of this morphism.

Since $\Coh^{d}(X)_{\text{rat}}$ is not a category fibered in groupoids, it is useful to work with the comma category instead of the usual fiber products of categories.
\begin{defn}
Suppose that we are given a diagram of pseudofunctors from $(\text{Aff}_{S})^{op}$ into categories as follows.
\begin{figure}[H]
\centering
\begin{tikzcd}
  &  \mathfrak{X} \ar[d, "f_{\mathfrak{X}}"] \\   \mathfrak{Y} \ar[r, "f_{\mathfrak{Y}}"] &  \mathfrak{Z}
\end{tikzcd}
\end{figure}
We define the right comma fiber product $\mathfrak{X} \vec{\times}_{\mathfrak{Z}} \mathfrak{Y}$ to be a pseudofunctor from $(\text{Aff}_{S})^{op}$ into categories. For every $T \in \text{Aff}_{S}$, the objects of $\mathfrak{X} \, \vec{\times}_{\mathfrak{Z}} \mathfrak{Y}(T)$ are triples $(x,y,g)$, where $x \in \mathfrak{X}(T)$, $y \in \mathfrak{Y}(T)$ and $g$ is a morphism $g: f_{\mathfrak{X}}(x) \to f_{\mathfrak{Y}} (y)$.

A morphism $(x_1, y_1, g_1)\to (x_2, y_2, g_2)$ consists of a pair $(\psi_{\mathfrak{X}}, \psi_{\mathfrak{Y}})$ of morphisms $\psi_{\mathfrak{X}}: x_1 \to x_2$ and $\psi_{\mathfrak{Y}}: y_1 \to y_2$ such that the following diagram is commutative
\begin{figure}[H]
\centering
\begin{tikzcd}
  f_{\mathfrak{X}}(x_1) \ar[r, "f_{\mathfrak{X}}(\psi_{\mathfrak{X}})"] \ar[d, "g_1", labels=left] &  f_{\mathfrak{X}}(x_2) \ar[d, "g_2"] \\   f_{\mathfrak{Y}}(y_1) \ar[r, "f_{\mathfrak{Y}}(\psi_{\mathfrak{Y}})"] &  f_{\mathfrak{Y}}(y_2)
\end{tikzcd}
\end{figure}
The composition of two composable pairs $(\psi_{\mathfrak{X}}, \psi_{\mathfrak{Y}})$, $(\xi_{\mathfrak{X}}, \xi_{\mathfrak{Y}})$ is given by the pointwise composition $(\xi_{\mathfrak{X}} \circ \psi_{\mathfrak{X}}, \xi_{\mathfrak{Y}} \circ\psi_{\mathfrak{Y}})$.
\end{defn}
By construction, $\mathfrak{X} \vec{\times}_{\mathfrak{Z}} \mathfrak{Y}$ is equipped with two projection morphisms $\pi_{\mathfrak{X}}: \mathfrak{X} \vec{\times}_{\mathfrak{Z}} \mathfrak{Y} \to \mathfrak{X}$ and $\pi_{\mathfrak{Y}}: \mathfrak{X} \vec{\times}_{\mathfrak{Z}} \mathfrak{Y} \to \mathfrak{Y}$ plus a $2$-morphism $H: f_{\mathfrak{X}} \circ \pi_{\mathfrak{X}} \Rightarrow f_{\mathfrak{Y}} \circ \pi_{\mathfrak{Y}}$ fitting into the diagram
\begin{figure}[H]
\centering
\begin{tikzcd}
   \mathfrak{X} \vec{\times}_{\mathfrak{Z}} \mathfrak{Y} \ar[d, "\pi_{\mathfrak{Y}}", labels=left] \ar[r, "\pi_{\mathfrak{X}}"] & \mathfrak{X} \ar[d, "f_{\mathfrak{X}}"] \ar[dl, symbol = \Longrightarrow] \\
  \mathfrak{Y} \ar[r, "f_{\mathfrak{Y}}", labels=below]  & \mathfrak{Z}
\end{tikzcd}
\end{figure}
This data is final for pseudofunctors $\mathfrak{C}$ equipped with morphisms $g_{\mathfrak{X}} : \mathfrak{C} \to \mathfrak{X}$ and  $g_{\mathfrak{Y}}: \mathfrak{C} \to \mathfrak{Y}$ plus a $2$-morphism $G: f_{\mathfrak{X}} \circ g_{\mathfrak{X}} \Rightarrow f_{\mathfrak{Y}} \circ g_{\mathfrak{Y}}$.

Fix the choice of an $S$-point in $\Coh^{d}(X)_{\text{rat}}(S)$ represented by a pair $(D, \mathcal{E})$.
\begin{defn}
The affine Grassmannian $\text{Gr}_{X, D, \mathcal{E}}$ is defined to be the right comma fiber product $\Coh^{d}(X) \, \vec{\times}_{\Coh^{d}(X)_{\text{rat}}} S$.
\end{defn}
For any $T \in \rm{Aff}_S$, denote $\cE_T := \cE|_{X_T}$. Let $U := X \setminus D$ and denote the inclusion $j : U \hookrightarrow X$. Then $\text{Gr}_{X, D, \mathcal{E}} (T)$ is the groupoid of pairs $(\mathcal{F}, \psi)$, where
\begin{enumerate}[(1)]
    \item $\mathcal{F}$ is a $T$-pure sheaf of dimension $d$ on $X_T$ such that $D_T$ is $\mathcal{F}$-regular.
    \item $\psi$ is monomorphism $\psi: \mathcal{E}_T \, \rightarrow  \, \mathcal{F}$ such that the restriction $\psi|_{U_{T}}$ is an isomorphism.
\end{enumerate}
An isomorphism of pairs between pairs $(\mathcal{F}_1, \psi_1)$ and $(\mathcal{F}_2, \psi_2)$ is an isomorphism $\mathcal{F}_1 \xrightarrow{\sim} \mathcal{F}_2$ that identifies the morphisms $\psi_1$ and $\psi_2$. Letting $j_T : U_T \hookrightarrow X_T$ denote the base change of $j$, Lemma \ref{lemma: properties of regularity} implies that $\mathcal{F}_i \rightarrow j_{T \, *} \, j_T^* \, \mathcal{F}_i$ is a monomorphism for $i=1,2$. It follows that an isomorphism between $(\cF_1,\psi_1)$ and $(\cF_2,\psi_2)$ is unique if it exists, so $\Gr_{X,D,\cE}$ can be regarded as a sheaf of sets.

The main result of this subsection is the following.
\begin{prop}\label{prop: ind-representability of step 1 grassmannian}
$\text{Gr}_{X, D, \mathcal{E}}$ is represented by a strict-ind scheme that is ind-projective over $S$.
\end{prop}
In order to prove Proposition \ref{prop: ind-representability of step 1 grassmannian}, we will proceed as in \cite[Thm. 1.1.3]{zhu_affine_grassmannians}. For any $(\mathcal{F}, \psi) \in \text{Gr}_{X, D, \mathcal{E}}(T)$, we use the isomorphism $j_{T *}j_{T}^*\psi : j_{T *}j_{T}^*\cE_T \to j_{T *}j_{T}^* \cF$ in order to view $\mathcal{F}$ as a subsheaf of $j_{T\, *} \, j_T^* \, \mathcal{E}_T$. We have $\mathcal{E}_{T} \subset \mathcal{F} \subset j_{T\, *} \, j_T^* \, \mathcal{E}_T$.
\begin{defn}
Let $N$ be a positive integer. We define $\text{Gr}^{\leq N}_{X, D, \mathcal{E}}$ to be the subfunctor of $\text{Gr}_{X, D, \mathcal{E}}$ that sends $T \in \text{Aff}_S$ to the set
\[ \text{Gr}^{\leq N}_{X, D, \mathcal{E}} (T) \; = \; \left\{ \begin{matrix} \; \; \; \text{pairs $(\mathcal{F}, \psi)$ in  $\text{Gr}_{X, D, \mathcal{E}}(T)$ such that} \; \; \; \\ \text{$\mathcal{E}_{T} \subset \mathcal{F} \subset \mathcal{E}_T(ND_T)$} \end{matrix} \right\}\]
\end{defn}
Note that there is a natural inclusion $\text{Gr}^{\leq N}_{\mathcal{E}, X/S} \subset \text{Gr}^{\leq M}_{\mathcal{E}, X/S}$ whenever $N \leq M$. Proposition \ref{prop: ind-representability of step 1 grassmannian} is a direct consequence of the following two lemmas.
\begin{lemma} \label{lemma: truncated grassmannian exhaust}
We have $\text{Gr}_{X, D, \mathcal{E}} = \underset{N>0}{\colim} \; \text{Gr}^{\leq{N}}_{X, D, \mathcal{E}}$ (as presheaves on $\text{Aff}_S$).
\end{lemma}
\begin{proof}
Fix $T$ in $\text{Aff}_S$. Let $(\mathcal{F}, \psi)$ be in $\text{Gr}_{X, D, \mathcal{E}}(T)$. We would like to show that there exists some $N >0$ such that $(\mathcal{F}, \psi)$ is in $\text{Gr}^{\leq N}_{X, D, \mathcal{E}}$. The scheme $X_T$ is quasi-compact, because it is of finite type over the affine scheme $T$. After passing to a finite affine cover, we can assume that $X_T$ is affine and $D_T$ is cut out by an element $x \in \mathcal{O}_{X_T}$. Then $\mathcal{F}$ and $\mathcal{E}_T$ are finitely presented modules without $x$-torsion. Let $(e_i)_{i \in I}$ be a finite set of generators for $\mathcal{F}$.

By assumption, we have $\mathcal{F}\left[\frac{1}{x}\right] = \mathcal{E}_T\left[\frac{1}{x}\right]$. This implies that for all $i \in I$, there exists a positive integer $n_i$ such that $x^{n_i}e_i \in \mathcal{E}_{T}$. Set $N \vcentcolon = \text{max}_{i \in I} \, (n_i)$. Then $\mathcal{F} \subset x^{-N}\mathcal{E}_T = \mathcal{E}_T(N D_T)$.
\end{proof}

\begin{lemma} \label{lemma: representability truncated version grassmannian}
For each $N>0$, the functor $\text{Gr}^{\leq{N}}_{X, D, \mathcal{E}}$ is represented by a disjoint union of schemes that are projective and of finite presentation over $S$. Moreover, for all $N \leq M$ the inclusion $\text{Gr}^{\leq N}_{X, D, \mathcal{E}} \subset \text{Gr}^{\leq M}_{X, D, \mathcal{E}}$ is a closed immersion.
\end{lemma}
\begin{proof}
We will use an auxiliary functor $\text{Q}^N_{X, D, \mathcal{E}}$. For $N>0$, we set $\text{Q}^N_{X, D, \mathcal{E}}$ to be the Quot-scheme $\text{Quot}_{X/S}\left(\mathcal{E}(ND)\, / \, \mathcal{E}\right)$ parametrizing $S$-flat quotients of $\mathcal{E}(ND)\, / \, \mathcal{E}$. A well-known theorem of Grothendieck states that $\text{Quot}_{X/S}\left(\mathcal{E}(ND)\, / \, \mathcal{E}\right)$ is represented by a disjoint union of projective schemes over $S$ \cite[Thm. 1.5.4]{olsson-stacks-book}.

For all $N>0$, we have a natural inclusion of functors $t_N: \text{Q}^N_{X, D, \mathcal{E}} \hookrightarrow \text{Q}^{N+1}_{X, D, \mathcal{E}}$ given as follows. Let $T$ in $\text{Aff}_S$ and $\mathcal{G}$ a quotient sheaf in $\text{Q}^N_{X, D, \mathcal{E}}(T)$. Let $q$ be the composition
\[ \mathcal{E}_T(ND_{T}) \; \twoheadrightarrow \; \frac{\mathcal{E}_T(ND_{T})}{\mathcal{E}_T} \; \twoheadrightarrow \; \mathcal{G}\]
Let $\text{Ker}(q)$ denote the kernel. Then, we set $t_N(T)(\mathcal{G})$ to be the quotient $\mathcal{E}_{T}((N+1)D_{T}) / \, \text{Ker}(q)$. We claim that $t_N(T)(\mathcal{G})$ is a well-defined element of $\text{Q}^{N+1}_{X, D, \mathcal{E}}$. The only thing to check is that the quotient $t_N(T)(\mathcal{G})$ is flat over $T$. Note that $t_N(T)(\mathcal{G})$ fits into a short exact sequence
\[ 0 \, \longrightarrow \, \mathcal{G} \, \longrightarrow \, t_N(T)(\mathcal{G}) \, \longrightarrow \, \frac{\mathcal{E}_T((N+1)D_{T})}{\mathcal{E}_T(ND_{T})} \, \longrightarrow \, 0 \]
Since $\mathcal{G}$ is $T$-flat, it suffices to show that $\mathcal{E}_T((N+1)D_{T})\, / \, \mathcal{E}_T(ND_{T})$ is $T$-flat. This is a consequence of Lemma \ref{lemma: properties of regularity}(c). Therefore, we have described a well-defined morphism $t_N: \text{Q}^N_{X, D, \mathcal{E}} \hookrightarrow \text{Q}^{N+1}_{X, D, \mathcal{E}}$. Since each Quot-scheme $\text{Q}^N_{X, D, \mathcal{E}}$ is a disjoint union of proper schemes over $S$, it follows from \cite[\href{https://stacks.math.columbia.edu/tag/01W6}{Tag 01W6}]{stacks-project} that $t_N$ is proper. By \cite[\href{https://stacks.math.columbia.edu/tag/04XV}{Tag 04XV}]{stacks-project} the proper monomorphism $t_{N}$ is a closed immersion.

To complete the proof of the lemma, we shall show that there is a collection of isomorphisms of functors $\phi_N: \text{Gr}^{\leq N}_{X, D, \mathcal{E}} \, \xrightarrow{\sim} \, \text{Q}^N_{X, D, \mathcal{E}}$ such that the following diagram commutes for any $N>0$.
\begin{equation} \label{eqn: grassmannian embedding}
\xymatrix{
  \Gr_{X, D, \mathcal{E}}^{\leq N+1} \ar[r]^{\phi_{N+1}} &  \text{Q}^{N+1}_{X, D, \mathcal{E}}\\   \Gr_{X, D, \mathcal{E}}^{\leq N} \ar@{^{(}->}[u] \ar[r]^{\phi_{N}} &  \text{Q}^N_{X, D, \mathcal{E}} \ar@{^{(}->}[u]
}
\end{equation}

For any $T \in \rm{Aff}_S$ and $T$-point $\mathcal{E}_T\subset \mathcal{F} \subset \mathcal{E}_T(ND_{T})$ of $\text{Gr}^{\leq N}_{X, D, \mathcal{E}}$, we set $\phi_N(T)(\mathcal{F})$ to be the quotient $ \frac{\mathcal{E}_T(ND_{T})}{\mathcal{E}_T} \; \twoheadrightarrow \; \frac{\mathcal{E}_T(ND_T)}{\mathcal{F}}$. On the other hand, suppose that $\mathcal{E}_T(ND_{T}) \,/\, \mathcal{E}_T \; \twoheadrightarrow \; \mathcal{G}$ is a $T$-point of $\text{Q}^N_{X, D, \mathcal{E}}$. We set $\tau_N(T)(\mathcal{G})$ to be the kernel of the composition $\mathcal{E}_T(ND_{T}) \; \twoheadrightarrow \; \frac{\mathcal{E}_T(ND_{T})}{\mathcal{E}_T} \; \twoheadrightarrow \; \mathcal{G}$. We need to check the following:
\begin{enumerate}[(i)]
    \item $\phi_N$ gives a well-defined natural transformation $\phi_N: \text{Gr}^{\leq N}_{X, D, \mathcal{E}} \, \longrightarrow \, \text{Q}^N_{X, D, \mathcal{E}}$.
    \item $\tau_N$ gives a well-defined natural transformation $\tau_N: \text{Q}^N_{X, D, \mathcal{E}} \, \longrightarrow \,\text{Gr}^{\leq N}_{X, D, \mathcal{E}}$.
    \item $\phi_N$ and $\tau_N$ are inverse of each other.
    \item The diagram \eqref{eqn: grassmannian embedding} commutes.
% \begin{figure}[H]
% \centering
% \begin{tikzcd}
%   \text{Gr}_{X, D, \mathcal{E}}^{\leq N+1} \ar[r, "\phi_{N+1}"] &  \text{Q}^{N+1}_{X, D, \mathcal{E}}\\   \text{Gr}_{X, D, \mathcal{E}}^{\leq N} \ar[u, symbol= \xhookrightarrow{\;\;}] \ar[r, "\phi_{N}"] &  \text{Q}^N_{X, D, \mathcal{E}} \ar[u, symbol= \xhookrightarrow{\;\;}]
% \end{tikzcd}
% \end{figure}
\end{enumerate}
The claims $(iii)$ and $(iv)$ are immediate from construction. %We are left to prove $(i)$ and $(ii)$.

\medskip
\noindent{\textit{Proof of (i):}}
\medskip

%\begin{enumerate}[(i)]
%    \item
    
Let $T$ in $\text{Aff}_S$ and $(\mathcal{F}, \psi)$ in $\text{Gr}^{\leq N}_{X, D, \mathcal{E}}(T)$. We have to show that $\mathcal{E}_T(ND_{T}) \, / \, \mathcal{F}$ is $T$-flat. We have a short exact sequence
\[0 \, \longrightarrow \frac{\mathcal{E}_T(ND_{T})}{\mathcal{F}} \, \longrightarrow \, \frac{j_{T \, *} \, j_T^* \, \mathcal{E}_T}{\mathcal{F}} \, \longrightarrow \, \frac{j_{T \, *} \, j_T^{*} \mathcal{E}_T}{\mathcal{E}_T(ND_{T})} \, \longrightarrow \, 0\]

So it suffices to show that both $\frac{j_{T \, *} \, j_T^* \, \mathcal{E}_T}{\mathcal{F}}$ and $\frac{j_{T \, *} \, j_T^{*} \mathcal{E}_T}{\mathcal{E}_T(ND_{T})}$ are $T$-flat. 

Let's show that $\frac{j_{T \, *} \, j_T^* \, \mathcal{E}_T}{\mathcal{F}}$ is $T$-flat. Note that $\psi$ induces an isomorphism $j_{T \, *} \, j_T^* \, \mathcal{E} \cong j_{T \, *} \, j_T^* \, \mathcal{F}$. Hence $\frac{j_{T \, *} \, j_T^* \, \mathcal{E}}{\mathcal{F}}$ is isomorphic to $\frac{j_{T \, *} \, j_T^* \, \mathcal{F}}{\mathcal{F}}$. Since $j_{T \, *} \, j_T^* \, \mathcal{F} = \underset{m}{\colim} \; \mathcal{F}(mD_T)$, we have $\frac{j_{T \, *} \, j_T^* \, \mathcal{F}}{\mathcal{F}} = \underset{m >0}{\colim} \; \frac{\mathcal{F}(mD_T)}{\mathcal{F}}$.

Since filtered colimits of $T$-flat sheaves are $T$-flat, it suffices to show that $\frac{\mathcal{F}(mD_T)}{\mathcal{F}}$ is $T$-flat for all $m>0$. This follows because $\frac{\mathcal{F}(mD_T)}{\mathcal{F}}$ is an iterated extension of sheaves of the form $\frac{\mathcal{F}(iD_T)}{\mathcal{F}((i+1)D_{T})}$, which are $T$-flat by Lemma \ref{lemma: properties of regularity}.
The proof that $\frac{j_{T \, *} \, j_T^{*} \mathcal{E}_T}{\mathcal{E}_T(ND_{T})}$ is $T$-flat is the same.
    
%    \item
\medskip
\noindent{\textit{Proof of (ii):}}
\medskip

Let $T$ in $\text{Aff}_S$ and $(\varphi, \mathcal{G})$ in $\text{Q}^{ N}_{X, D, \mathcal{E}}(T)$. We need to check that $\tau_N(T)(\mathcal{G})$ is $T$-pure of dimension $d$. There is a short exact sequence 
\[ 0 \, \longrightarrow \, \tau_N(T)(\mathcal{G}) \, \longrightarrow \, \mathcal{E}_T(ND_{T}) \, \longrightarrow \, \mathcal{G} \, \longrightarrow 0 \]
Since both $\mathcal{E}_T(ND_{T})$ and $\mathcal{G}$ are $T$-flat, we conclude that  $\tau_N(T)(\mathcal{G})$ is $T$-flat.
Let $t \in T$. Since $\mathcal{G}$ is $T$-flat, the short exact sequence above remains exact when we restrict to the fiber $X_t$
\[ 0 \, \longrightarrow \, \tau_N(T)(\mathcal{G})|_{X_t} \, \longrightarrow \, \mathcal{E}_T(ND_{T})|_{X_t} \, \longrightarrow \, \mathcal{G}|_{X_t} \, \longrightarrow 0 \]
Since $\mathcal{E}_T(ND_{T})|_{X_t}$ is pure of dimension $d$, we conclude that its subsheaf $\tau_N(T)(\mathcal{G})|_{X_t}$ is pure of dimension $d$ as well. This shows that $\tau_N(T)(\mathcal{G})$ $T$-pure of dimension $d$.
%\end{enumerate}
\end{proof}

\begin{remark}
As discussed in Remark \ref{rem: variants of coh rat}, there are several variants of the definition of $\Gr_{X,D,\cE}$ that one could use to formulate infinite dimensional GIT. For instance, instead of a $T$-point being defined by a morphism $\psi : \cE_T \to \cF$ whose restriction to $U_T$ is an isomorphism, one could simply ask for the data of an isomorphism $\cE_T|_{U_T} \to \cF|_{U_T}$. The resulting affine Grassmannian is a colimit of the sequence of closed immersions
\[
\cdots \hookrightarrow \Gr_{X,D,\cE(nD)} \hookrightarrow \Gr_{X,D,\cE((n-1)D)} \hookrightarrow \cdots
\]
induced by the maps $\cE((n-1)D) \to \cE(nD)$.

Another variant is to take the fibers of the canonical morphism $\Coh^d(X) \to \Coh^d(X)_{rat}$ for the alternative definition of $\Coh^d(X)_{rat}$ discussed in Remark \ref{rem: variants of coh rat}. The resulting affine Grassmannian is still an ind-projective ind-scheme, but it is much larger. This construction is closer in spirit to the ``rational affine Grassmannian" $\text{Gr}_{\text{Ran} \, X}$ studied in \cite[\S 5]{gaitsgory-contractibility-rational}.
\end{remark}

\end{subsection}
\begin{subsection}{$L_n^{\vee}$ is asymptotically ample on $\text{Gr}_{X, D, \mathcal{E}}$} \label{subsection: asymptotic ampleness line bundle pure sheaves}
There is a natural forgetful morphism $\text{Gr}_{X, D, \mathcal{E}} \rightarrow \Coh^{d}(X)$  that sends a pair $(\mathcal{F}, \psi)$ to the sheaf $\mathcal{F}$. We will abuse notation and denote by $M_n, L_n$ and $b_d$ the line bundles on $\text{Gr}_{X, D, \mathcal{E}}$ obtained by pulling back the corresponding line bundles on $\Coh^{d}(X)$ under the forgetful morphism (see Definitions \ref{defn: m_n line bundle}, \ref{defn: knudsen-mumford coefficients} and \ref{defn: l_n line bundle}).

\begin{defn} \label{defn: locus of grassmannian with given hilbert polynomial}
Let $P \in \mathbb{Q}[x]$ be an integer valued polynomial. We will denote by $\text{Gr}^{P}_{X, D, \mathcal{E}}$ the subfunctor of $\text{Gr}_{X, D, \mathcal{E}}$ that sends any $T$ in $\text{Aff}_S$ to the set
\[ \text{Gr}^{P}_{X, D, \mathcal{E}} (T) \; = \; \left\{ \begin{matrix} \; \; \; \text{pairs $(\mathcal{F}, \psi)$ in  $\text{Gr}_{X, D, \mathcal{E}}(T)$ such that} \; \; \; \\ \text{$P_{\mathcal{F}|_{X_t}} = P$ for all points $t \in T$} \end{matrix} \right\}\]
Similarly, we define $\text{Gr}^{\leq N, P}_{X, D, \mathcal{E}} \vcentcolon = \, \text{Gr}^{P}_{X, D, \mathcal{E}} \, \cap \, \text{Gr}^{\leq N}_{X, D, \mathcal{E}}$.
\end{defn}

\begin{prop} \label{prop: line bundle affine grassmannian torsion-free is asymptotically ample}
Fix $N>0$ and $P \in \mathbb{Q}[x]$. The functor $\text{Gr}^{\leq N, P}_{X, D, \mathcal{E}}$ is represented by a projective scheme over $S$. Furthermore, there exists $m \in \mathbb{N}$ (depending on $N$ and $P$) such that for all $n \geq m$ the line bundle $L_{n}^{\vee}|_{\text{Gr}_{X, D, \mathcal{E}}^{\leq N, P}}$ is $S$-ample.
\end{prop}

In order to prove this proposition, we will define auxiliary families of line bundles on $\text{Gr}_{X, D, \mathcal{E}}^{\leq N, P}$. For the rest of this subsection we fix $N \in \mathbb{N}$ and $P \in \mathbb{Q}[x]$. We denote by $\overline{p}$ the reduced Hilbert polynomial corresponding to $P$ ($\overline{p}$ is the unique constant multiple of $P$ with leading coefficient $\frac{1}{d!}$).

\begin{defn} \label{defn: auxiliary line bundles affine grassmannian torsion-free}
Let $n \in \mathbb{N}$. Let $p: \text{Gr}^{\leq N, P}_{X, D, \mathcal{E}} \rightarrow S$ denote the structure morphism  of $\text{Gr}^{\leq N, P}_{X, D, \mathcal{E}}$. We define the line bundle $\widetilde{M}_n$ on $\text{Gr}_{X, D, \mathcal{E}}^{\leq N, P}$ by
\[\widetilde{M}_n \vcentcolon = M_n^{\vee} \; \otimes \; p^*\,  \text{det} \, R\pi_*\left(\mathcal{E}(ND)(n)  \right)\]

Furthermore, we define $\widetilde{b}_{d}\vcentcolon = \bigotimes_{i = 0}^d (\widetilde{M}_{i})^{(-1)^{d-i} \binom{d}{i}}$ and set $\widetilde{L}_n \vcentcolon = \widetilde{M}_n \otimes (\widetilde{b}_{d})^{- \otimes \overline{p}(n)}$.
\end{defn}

 Proposition \ref{prop: line bundle affine grassmannian torsion-free is asymptotically ample} now follows from Lemmas \ref{lemma: shifting indexes mu line bundle torsion-free grassmannian} and \ref{lemma: ampleness for line bundle with big enough index torsion-free grassmannian} below.

\begin{lemma} \label{lemma: shifting indexes mu line bundle torsion-free grassmannian}
Fix $N>0$ and $P \in \mathbb{Q}[x]$. The line bundle $L_{n}^{\vee}$ is $S$-ample on $\text{Gr}_{X, D, \mathcal{E}}^{\leq N, P}$ if and only if $\widetilde{L}_n$ is $S$-ample.
\end{lemma}
\begin{proof}
It suffices to show that the line bundles $L_n^{\vee}$ and $\widetilde{L}_{n}$ differ by a line bundle coming from the base $S$. This follows from the fact that $M_n^{\vee}$ and $\widetilde{M}_n$ differ by a line bundle coming from $S$, which implies the same for $\widetilde{b}_d$ and $\widetilde{L}_n$.
\end{proof}
\begin{lemma} \label{lemma: ampleness for line bundle with big enough index torsion-free grassmannian}
Fix $N>0$ and $P \in \mathbb{Q}[x]$. Then $\text{Gr}^{\leq N, P}_{X, D, \mathcal{E}}$ is a projective scheme over $S$. Moreover, there exists $m \in \mathbb{N}$ (depending on $N$) such that for all $n \geq m$ the line bundle $\widetilde{L}_n$ is $S$-ample on $\text{Gr}^{\leq N, P}_{X, D, \mathcal{E}}$.
\end{lemma}
\begin{proof}
The isomorphism $\phi_N$ described in \Cref{lemma: representability truncated version grassmannian} identifies $\text{Gr}_{X, D, \mathcal{E}}^{\leq N, P}$ with the component $\text{Quot}^{H}_{X/S}\left(\mathcal{E}(ND) / \mathcal{E}\right)$ of the $\text{quot}$-scheme determined by the Hilbert polynomial $H = P_{\mathcal{E}(ND)} - P$. In particular, $\text{Gr}^{\leq N, P}_{X, D, \mathcal{E}}$ is projective over $S$.

It suffices to prove the claim when the base is a field, because $S$ is quasi-compact, ampleness is an open condition on the base \cite[Corollaire 9.6.4]{egaiv}, and ampleness can be checked on fibers \cite[Corollaire 9.6.5]{egaiv}. After base-change we can assume that $k$ is infinite. We will keep these assumptions for the rest of the proof.

We start by showing that $\widetilde{M}_n$ is ample for $n$ big enough. We can apply cohomology and base-change \cite[\href{https://stacks.math.columbia.edu/tag/0A1D}{Tag 0A1D}]{stacks-project} to the fiber diagram 
\begin{figure}[H]
\centering
\begin{tikzcd}
  X \times_{S} \text{Gr}_{X, D, \mathcal{E}}^{\leq N, P} \ar[r] \ar[d, "\pi'", labels=left] &  X \ar[d, "\pi"] \\   \text{Gr}_{X, D, \mathcal{E}}^{\leq N,P} \ar[r, "p"] &  S
\end{tikzcd}
\end{figure}
to conclude that $ p^*\,  \text{det} \, R\pi_*\left(\mathcal{E}(ND)(n)  \right)  \; \cong \; \text{det} \, R\pi'_{*}\left((\mathcal{E}(ND))_{\text{Gr}^{\leq N, P}_{X, D, \mathcal{E}}}(n)  \right)$, where we are using the notation $(\mathcal{E}(ND))_{\text{Gr}^{\leq N, P}_{X, D, \mathcal{E}}}$ to denote the base-change $\mathcal{E}(ND)|_{X \times_S \text{Gr}^{\leq N, P}_{X, D, \mathcal{E}}}$. Therefore,
\[ \widetilde{M}_n = M_n^{\vee} \; \otimes \; \text{det} \, R\pi'_{*}\left((\mathcal{E}(ND))_{\text{Gr}^{\leq N, P}_{X, D, \mathcal{E}}}(n) \right)\]
Let $(\mathcal{F}_{univ}, \psi_{univ})$ denote the universal bundle on $X \times \text{Gr}_{X, D, \mathcal{E}}^{\leq N, P}$. The additivity of the determinant implies that for all $n \in \mathbb{Z}$ we have
\[ \widetilde{M}_n = \text{det} \, R \pi'_{*}\left( \left((\mathcal{E}(ND))_{\text{Gr}_{X, D, \mathcal{E}}^{\leq N, P}} \, / \,\mathcal{F}_{univ} \right) (n)\, \right)\]

Under the isomorphism $\phi_{N} : \text{Gr}_{X, D, \mathcal{E}}^{\leq N, P} \xrightarrow{\sim} \text{Quot}^{H}_{X}\left(\mathcal{E}(ND) / \mathcal{E}\right)$, the line bundle $\widetilde{M}_n$ corresponds to $\text{det} \, R\pi'_{*}\left(\mathcal{Q}_{univ}(n)\right)$ on the quot-scheme. Here 
\[\mathcal{Q}_{univ} \vcentcolon = \frac{\left(\mathcal{E}(ND) / \mathcal{E}\right)_{\text{Gr}_{X, D, \mathcal{E}}^{\leq N, P}}}{\mathcal{F}_{univ}/\mathcal{E}_{\text{Gr}_{X, D, \mathcal{E}}^{\leq N, P}}}\]
denotes the universal quotient. It is known \cite[Prop. 2.2.5]{huybrechts.lehn} that $\text{det} \, R\pi'_{*}\left(\mathcal{Q}_{univ}(n)\right)$ is ample for $n$ big enough. It follows that $\widetilde{M}_n$ is ample for $n$ big enough. To conclude the proof, we show that there exists some positive integer $r$ such that $\widetilde{b}_{d}^{\otimes r}$ is trivial, which implies that $\widetilde{L}_n^{\otimes r} \cong \widetilde{M}_n^{\otimes r}$ is ample.

Notice that the support of $\mathcal{E}(ND) / \mathcal{E}$ has dimension $\leq d-1$. After replacing $X$ with this support, we can assume that $X$ has dimension $\leq d-1$. Now \cite[Theorem 4]{knudsen-mumford} shows that the family of line bundles $\widetilde{M}_n$ is a polynomial of degree $d$ in the variable $n$ with coefficients in the Picard group. By definition we see that $\widetilde{b}_{d}$ is the leading coefficient of this polynomial. This leading coefficient $\widetilde{b}_{d}$ gets replaced by a multiple $\widetilde{b}_{d}^{\otimes r^{d}}$ whenever we replace $\mathcal{O}(1)$ by a multiple $\mathcal{O}(r)$. Hence, we can assume without loss of generality that $\mathcal{O}(1)$ is very ample.

We can use $\mathcal{O}(1)$ to embed $X$ into a projective space $\mathbb{P}^m_{k}$. Choose a linear subspace $L$ of dimension $m-d$ in $\mathbb{P}^{m}_{k}$ that is disjoint from $X$. Then the corresponding projection $\mathbb{P}_k^{m} \setminus L \rightarrow \mathbb{P}_k^{d-1}$ restricts to a finite map $f: X \rightarrow \mathbb{P}_k^{d-1}$. The Quot-scheme $\text{Quot}^{H}_{X}\left(\mathcal{E}(ND) / \mathcal{E}\right)$ can be identified with the functor parameterizing flat families of quotients of $f_* \, (\mathcal{E}(ND)/ \mathcal{E})$ as a module over the coherent $\cO_{\bP^n_k}$-algebra $f_\ast(\cO_{X})$. Such a quotient is determined by the quotient of underlying $\cO_{\bP^n_k}$-modules, so the morphism $\text{Quot}^{H}_{X}\left(\mathcal{E}(ND) / \mathcal{E}\right) \to \text{Quot}^{H}_{\mathbb{P}_k^{d-1}}\left(f_* \, (\mathcal{E}(ND)/ \mathcal{E}) \, \right)$ that forgets the $f_\ast(\cO_X)$-module structure is a proper monomorphism, and hence a closed immersion. Note that the formation of $\widetilde{M}_n$ is compatible with this closed immersion, because derived pushforward commutes with Tor-independent base change \cite[\href{https://stacks.math.columbia.edu/tag/0A1H}{Tag 0A1H}]{stacks-project}. We have therefore reduced to the case when $X = \mathbb{P}^{d-1}_k$.

Set $\text{Quot} \vcentcolon = \text{Quot}^{H}_{\mathbb{P}_k^{d-1}}\left(f_* \, (\mathcal{E}(ND) / \mathcal{E}) \, \right)$. Let $\mathcal{Q}_{univ}$ denote the universal quotient on $\mathbb{P}_k^{d-1} \times \text{Quot}$. Let $p$ denote the structure morphism $p: \mathbb{P}_k^{d-1} \times \text{Quot} \rightarrow \text{Quot}$. By definition we have $\widetilde{M}_n = \text{det} \, Rp_*\left(\mathcal{Q}_{univ}(n)\right)$. By \Cref{lemma: linear resolution in projective space} below applied to $T = \text{Quot}$, there is a finite resolution
\[ 0 \rightarrow p^*\mathcal{V}_l(a_l) \rightarrow p^*\mathcal{V}_{l-1}(a_{l-1}) \rightarrow \cdots \rightarrow p^*\mathcal{V}_{1}(a_1) \rightarrow \mathcal{Q}_{univ} \rightarrow 0\]
where each $\mathcal{V}_{i}$ is a vector bundle on $\text{Quot}$. By using additivity of the determinant and the projection formula, we conclude that
\[\widetilde{M}_n =   \bigotimes_{i=1}^{l} \text{det}(\mathcal{V}_{i})^{\otimes (-1)^{i+1} \chi(\mathcal{O}(n+a_i))}\]
Here $\chi(\mathcal{O}(n+a_i))$ denotes the Euler characteristic of $\mathcal{O}(n+a_i)$ on $\mathbb{P}_k^{d-1}$. This is a polynomial of degree $d-1$ in the variable $n$. So we see that $\widetilde{M}_n$ is a polynomial of degree $d-1$ in the variable $n$ with coefficients in the Picard group of $\text{Quot}$. By definition $\widetilde{b}_{d}$ is the coefficient of $n^d$, which is trivial as desired.
\end{proof}

Here we include the resolution result needed for the proof of the last lemma.
\begin{lemma} \label{lemma: linear resolution in projective space}
Let $T$ be a Noetherian scheme and $d$ a positive integer. Consider the projective space $p: \mathbb{P}^{d}_{T} \to T$. Let $\mathcal{Q}$ be a $T$-flat coherent sheaf on $\mathbb{P}^d_{T}$. Then, there exists an integer $n \in \mathbb{Z}$ and a tuple of vector bundles $(\mathcal{V}_i)_{i=0}^{d}$ on $T$ such that $\mathcal{Q}$ admits a resolution of the form
\[ 0 \rightarrow p^*\mathcal{V}_d(-n-d) \rightarrow p^*\mathcal{V}_{d-1}(-n-d+1) \rightarrow \cdots \rightarrow p^{*}\mathcal{V}_1(-n-1) \rightarrow p^*\mathcal{V}_{0}(-n) \rightarrow \mathcal{Q} \rightarrow 0 \]
\end{lemma}
\begin{proof}
Consider the diagram given by the first and second projections.
\begin{figure}[H]
\centering
\begin{tikzcd}
  & \mathbb{P}^d_{T} \times_{T} \mathbb{P}^d_{T} \ar[dl, "q_1", labels = above left] \ar[dr, "q_2"] &  \\
  \mathbb{P}^d_{T} & & \mathbb{P}^d_{T}
\end{tikzcd}
\end{figure}
Let $\mathcal{O}_{\Delta}$ denote the structure sheaf of the diagonal $\mathbb{P}^d_{T} \hookrightarrow \mathbb{P}^d_{T} \times_{T} \mathbb{P}^d_{T}$. By Beilinson's resolution of the diagonal \cite{beilinson-resolution-diagonal} \cite[Lemma 8.27]{huybrechts-fourier-mukai-book}, $\mathcal{O}_{\Delta}$ admits a locally free resolution
\begin{gather*}
 0 \to \Omega^d_{\mathbb{P}^d_{T}/T}(d) \, \boxtimes \, \mathcal{O}(-d) \to  \Omega^{d-1}_{\mathbb{P}^d_{T}/T}(d-1) \, \boxtimes \, \mathcal{O}(-d+1) \to \cdots \to \Omega^1_{\mathbb{P}^d_{T}/T}(1) \, \boxtimes \, \mathcal{O}(-1) \to \mathcal{O} \, \boxtimes \, \mathcal{O} \to \mathcal{O}_{\Delta} \to 0 
\end{gather*}
For any integer $n$, we can tensor with the pullback $q_1^{*} \mathcal{Q}(n)$to obtain the following acyclic complex in $\mathbb{P}^d_{T} \times_{T} \mathbb{P}^d_{T}$
\begin{gather*}
    \mathcal{C}^{\bullet}_n = \; \left[ \, 0 \to \left(\Omega^d_{\mathbb{P}^d_{T}/T} \otimes \mathcal{Q}(n+d)\right)\, \boxtimes \, \mathcal{O}(-d) \to \cdots \to \left(\Omega^1_{\mathbb{P}^d_{T}/T} \otimes \mathcal{Q}(n+1) \right) \, \boxtimes \, \mathcal{O}(-1) \to \mathcal{Q}(n) \, \boxtimes \, \mathcal{O} \to \mathcal{O}_{\Delta} \otimes (q_1)^* \mathcal{Q}(n) \to 0 \, \right]
\end{gather*}
Choose $n \gg 0$ so that $R^jp_{*}\left(\Omega^i_{\mathbb{P}^d_{T}/T} \otimes \mathcal{Q}(n+i)\right) = 0$ for all $j \geq 1$ and $i \geq 0$. Then by applying \cite[\href{https://stacks.math.columbia.edu/tag/0A1H}{Tag 0A1H}]{stacks-project} and truncation, it follows that 
\[\mathcal{V}_i \vcentcolon = p_*\left(\Omega^i_{\mathbb{P}^d_{T}/T} \otimes \mathcal{Q}(n+i)\right)\]
is a locally-free $\mathcal{O}_T$-module for all $i \geq 0$. Also, the projection formula implies that $R^j (q_2)_*\left( \Omega^i_{\mathbb{P}^d_{T}/T} \otimes \mathcal{Q}(n+i) \boxtimes \mathcal{O}(-i) \right)= 0$ for all $j \geq 1$ and $i \geq 0$. The first page of the spectral sequence for the derived pushforward $R(q_2)_* \, \mathcal{C}^{\bullet}_n$ of the complex $\mathcal{C}^{\bullet}_n$ is concentrated in a single column and shows that the pushforward complex
\begin{gather*}
    (q_2)_* \, \mathcal{C}^{\bullet}_n = \; \left[ \, 0 \to p^*\mathcal{V}_d(-d) \to p^*\mathcal{V}_{d-1}(-d+1) \to \cdots \to p^{*}\mathcal{V}_1(-1) \to p^*\mathcal{V}_{0} \to \mathcal{Q}(n) \to 0 \,\right]
\end{gather*} 
is acyclic. We obtain our desired resolution by tensoring with $\mathcal{O}(-n)$.
\end{proof}

\end{subsection}

\begin{subsection}{Affine Grassmannians for related moduli} \label{section: affine grassmannians}
Set $\mathcal{M}$ to be either $\Pair_{\mathcal{A}}^{d}(X)$ or $\Lambda\Coh^{d}(X)$. Let $(D , \mathcal{E}, \beta)$ be an $S$-point of the stack of rational maps $\mathcal{M}_{rat}$.
\begin{defn} \label{defn: Affine Grassmannians for related moduli}
We define the affine Grassmannian $\text{Gr}_{\mathcal{M}, X, D, \mathcal{E}, \beta}$ to be the right comma fiber product $\mathcal{M} \, \vec{\times}_{\mathcal{M}_{rat}} S$.
\end{defn}
\begin{notn}
In latter sections we will occasionally simplify notation and write $\text{Gr}_{\mathcal{M}}$ instead of $\text{Gr}_{\mathcal{M}, X, D, \mathcal{E}, \beta}$ (omitting the data $X, D, \cE, \beta$) whenever this extra data is clear from context.
\end{notn}
By definition, $\text{Gr}_{\mathcal{M}, X, D, \mathcal{E}, \beta}(T)$ is the set of triples $(\mathcal{F}, \psi, \tilde{\beta})$ where $(\mathcal{F}, \psi) \in \text{Gr}_{X, D, \mathcal{E}}(T)$ and $\widetilde{\beta}$ is a morphism $\mathcal{A}|_{X_{T}} \to \mathcal{F}$ (resp. a $\Lambda|_{X_{T}}$-module structure on $\mathcal{F}$) such that $\widetilde{\beta}|_{X_{T} \setminus D_{T}}$ is identified with $\beta|_{X_{T} \setminus D_{T}}$ under the isomorphism $\psi|_{X_{T} \setminus D_{T}}$. There is a forgetful morphism 
\[\text{Gr}_{\mathcal{M}, X, D, \mathcal{E}, \beta} \rightarrow \text{Gr}_{X, D, \mathcal{E}}, \; \; \; \; (\mathcal{F}, \psi, \widetilde{\beta}) \mapsto (\mathcal{F}, \psi)\]
\begin{prop} \label{prop: representability of affine grassmannian for other moduli}
The forgetful morphism exhibits $\text{Gr}_{\mathcal{M}, X, D, \mathcal{E}, \beta}$ as a closed ind-subscheme of $\text{Gr}_{X, D, \mathcal{E}}$.
\end{prop}
\begin{proof}
We set $U$ to be the open complement $X \setminus D$, $T$ an $S$-scheme, and $T \rightarrow \text{Gr}_{X, D, \mathcal{E}}$ a morphism represented by a pair $(\mathcal{F}, \psi)$.

\smallskip
\noindent{\textit{Proof for $\mathcal{M} = \Pair_{\mathcal{A}}^{d}(X)$:}}
\smallskip

Consider the composition $\gamma = \left(\psi|_{U_{T}}\right)^{-1} \circ \beta|_{U_{T}} : \mathcal{A}|_{U_{T}} \rightarrow \mathcal{F}|_{U_{T}}$.
The fiber product $\text{Gr}_{\mathcal{M}, X, D, \mathcal{E}, \beta} \times_{\text{Gr}_{X, D, \mathcal{E}}} T$ is the functor that takes a $T$-scheme $Y$ to the set of morphisms $\widetilde{\beta}: \mathcal{A}|_{X_{Y}} \rightarrow \mathcal{F}|_{X_{Y}}$ such that $\widetilde{\beta}|_{U_{Y}} = \beta|_{U_{Y}}$. Since $D_{T}$ is $\mathcal{F}$-regular, if such $\widetilde{\beta}$ exists then it is unique. 

There exists some $m\gg0$ and a morphism $ \mathcal{A}|_{X_{T}} \rightarrow \mathcal{F}(mD_{T})$ extending $\gamma$. Consider the composition
\[ \xi: \mathcal{A}|_{X_{T}} \rightarrow \mathcal{F}(mD_{T}) \twoheadrightarrow \mathcal{F}(mD_{T}) / \mathcal{F} \]

Then $\text{Gr}_{\mathcal{M}, X, D, \mathcal{E}, \beta} \times_{\text{Gr}_{X, D, \mathcal{E}}} T$ is the subfunctor of $T$ consisting of morphisms $Y \rightarrow T$ such that $\xi|_{X_{Y}} = 0$. We know that $\mathcal{F}(mD_{T}) / \mathcal{F}$ is $T$-flat by Lemma \ref{lemma: properties of regularity} (c). Therefore, Lemma \ref{lemma: representability of hom sheaves} can be applied to show that the functor $\text{Hom}(\mathcal{A}|_{X_{T}}, \, \mathcal{F}(mD_{T}) / \mathcal{F})$ is representable by a separated scheme over $T$. The morphisms $\xi$ and $0$ determine two sections of the structure morphism $ \text{Hom}(\mathcal{A}|_{X_{T}}, \, \mathcal{F}(mD_{T}) / \mathcal{F}) \rightarrow T$. We conclude that $\text{Gr}_{\Pair, X, D, \mathcal{E}, \beta} \times_{\text{Gr}_{X, D, \mathcal{E}}} T$ is represented by the closed subscheme of $T$ where these two sections agree.

\smallskip
\noindent{\textit{Proof for $\mathcal{M} = \Lambda\Coh^{d}(X)$:}}
\smallskip

We can use the $\Lambda|_{U_{T}}$-module structure on $\mathcal{E}|_{U_{T}}$ and the isomorphism $\psi|_{U_{T}}: \mathcal{E}|_{U_{T}} \xrightarrow{\sim} \mathcal{F}|_{U_{T}}$ in order to define a $\Lambda|_{U_{T}}$-module structure on $\mathcal{F}|_{U_{T}}$. By Proposition \ref{prop: characterization of lambda modules}, this amounts to the data of a morphism $b: \Lambda_1|_{U_{T}} \otimes \mathcal{F}|_{U_{T}} \rightarrow \mathcal{F}|_{U_{T}}$ satisfying conditions $(\Lambda 1)$ and $(\Lambda 2)$ in that proposition. 

The fiber product $\text{Gr}_{\Lambda\Coh, X, D, \mathcal{E}, a} \times_{\text{Gr}_{X, D, \mathcal{E}}} T$ is the functor that takes a $T$-scheme $Y$ to the set of morphisms $\widetilde{b}:  \Lambda_{X_{Y}} \otimes \mathcal{F} \rightarrow \mathcal{F}$ that extend $b|_{U_{Y}}$ and satisfy $(\Lambda 1)$ and $(\Lambda 2)$. 

Since $D_{T}$ is $\mathcal{F}$-regular, it follows that $\widetilde{b}$ is uniquely determined whenever it exits. By the same uniqueness argument, if $\widetilde{b}$ exist then $(\Lambda 1)$ and $(\Lambda 2)$ are automatically satisfied, because they are already satisfied for the restriction $b|_{U_{Y}}$. This shows that $\text{Gr}_{\Lambda\Coh, X, D, \mathcal{E}, \beta} \times_{\text{Gr}_{X, D, \mathcal{E}}} T$ is the subfunctor of $T$ consisting of morphisms $Y \rightarrow T$ such that the morphism $b$ admits an extension to $X_{Y}$. The same argument as in the case of pairs shows that this is represented by a closed subscheme of $T$, as desired.
\end{proof}
These affine Grassmannians admit natural morphisms to $\mathcal{M}$ given by forgetting the isomorphism $\psi$ defined away from $D$. We have just seen that there is a closed immersion into $\text{Gr}_{X, D, \mathcal{E}}$ such that the following diagram commutes
\begin{figure}[H]
\centering
\begin{tikzcd}
  \text{Gr}_{\mathcal{M}, X, D, \mathcal{E}, \beta} \ar[r] \ar[d, symbol = \xhookrightarrow{\;\;\;\;\;\;}] & \mathcal{M} \ar[d, "Forget"] \\  \text{Gr}_{X, D, \mathcal{E}} \ar[r] & \Coh^{d}(X)
\end{tikzcd}
\end{figure}

We use the forgetful morphism $\text{Gr}_{\mathcal{M}, X, D, \mathcal{E}, \beta} \rightarrow \mathcal{M}$ to restrict the line bundle $L_n$ to $\text{Gr}_{\mathcal{M}, X, D, \mathcal{E}, \beta}$. This is the same as the pullback of the line bundle $L_n$ on $\text{Gr}_{X, D, \mathcal{E}}$. Therefore, we obtain the following corollary as an immediate consequence of Proposition \ref{prop: line bundle affine grassmannian torsion-free is asymptotically ample}.
\begin{coroll} \label{coroll: ampleness line bundle affine grassmannian}
Let $N \in \mathbb{N}$ and $P \in \mathbb{Q}[x]$. Then, there exists some $m \gg0$ such that for all $n \geq m$, the line bundle $L_{n}^{\vee}$ on $\text{Gr}_{\mathcal{M}, X, D, \mathcal{E}, \beta}^{\leq N, P} \vcentcolon = \text{Gr}_{\mathcal{M}, X, D, \mathcal{E}, \beta} \cap \text{Gr}_{X, D, \mathcal{E}}^{\leq N, P}$ is $S$-ample.
\qed
\end{coroll}
\end{subsection}
\end{section}

\begin{section}{Monotonicity via ``infinite dimensional GIT"} \label{section: infinite dimensional GIT}
In this section, we define a polynomial numerical invariant $\nu$ on the moduli stacks we are considering. Furthermore, we prove that $\nu$ is strictly $\Theta$-monotone and strictly $S$-monotone.
\begin{subsection}{The numerical invariant} \label{section: numerical invariant}
We refer the reader to Subsection \ref{subsection: background theta stratifications} for a discussion of filtrations and graded points. In the case of the stack $\Coh^{d}(X)$, a $k$-valued point is the same thing as morphism $\Spec(k) \rightarrow S$ and a pure sheaf $\mathcal{F}$ of dimension $d$ on $X_k$, and a filtration of $\mathcal{F}$ in the stack $\Coh^{d}(X)$ is a $\mathbb{G}_m$-equivariant $\mathbb{A}^1_k$-flat relative torsion-free sheaf $\widetilde{\mathcal{F}}$ on $X \times \mathbb{A}^1_k$. Using the Rees construction \cite[Prop. 1.0.1]{halpernleistner2018structure} we can view this as a sequence $(\mathcal{F}_m)_{m \in \mathbb{Z}}$ of subsheaves of $\mathcal{F}$ satisfying
\begin{enumerate}[(a)]
    \item $\mathcal{F}_{m+1} \subset \mathcal{F}_{m}$.
    \item $\mathcal{F}_{m}/\mathcal{F}_{m+1}$ is pure of dimension $d$.
    \item $\mathcal{F}_m = 0$ for $m \gg0$ and $\mathcal{F}_m = \mathcal{F}$ for $m\ll0$.
\end{enumerate}
The line bundle $L_n|_{\Theta_k}$ for such a filtration will be a $\mathbb{G}_m$-equivariant line bundle on $\mathbb{A}^1_k$. All such line bundles come from the base $[\Spec(k)/\mathbb{G}_m]$, and so they are classified by the isomorphism class of the fiber $L_n|_{0}$ as a one dimensional $\mathbb{G}_m$-representation. These isomorphism classes are just characters of $\mathbb{G}_m$, which are classified by an integer called the weight. We will denote this integer by $\wt\left( L_n|_{0} \right)$.

By \cite[\href{https://stacks.math.columbia.edu/tag/0E9I}{Tag 0A1D}]{stacks-project}, the restriction $M_n|_{0}$ of the line bundle $M_n$ (Definition \ref{defn: m_n line bundle}) is isomorphic to $\text{det} \, \left(R\pi_{k \, *} \, (\widetilde{\mathcal{F}}|_{0})(n)\right)$. We have an equality of graded sheaves $\widetilde{\mathcal{F}}|_{0} = \bigoplus_{m \in \mathbb{Z}} \left(\mathcal{F}_m/\mathcal{F}_{m+1}\right)$. Here $\mathbb{G}_m$ acts on $\mathcal{F}_m/\mathcal{F}_{m+1}$ with weight $m$. We can take derived pushforward and determinant in order to conclude that
\[ \wt \left( M_n|_{0}\right) = \sum_{m \in \mathbb{Z}} \left(P_{\mathcal{F}_m}(n) - P_{\mathcal{F}_{m+1}}(n) \right) \cdot m\]
In particular, this yields $\wt(b_{d+1}|_0) = 0$, hence the expansion in \Cref{defn: knudsen-mumford coefficients} implies that
\[ \wt \left( b_d|_{0}\right) = d! \cdot \left( \text{coefficient of } n^d \text{ in } \wt(M_n|_0) \right) = \sum_{m \in \mathbb{Z}} \left(\rk_{\mathcal{F}_m} - \rk_{\mathcal{F}_{m+1}}  \right) \cdot m\]
We conclude that $\wt (L_n|_{0})$ is given by:
\[ \wt \left( L_n|_{0}\right) = \wt(M_n|_0) - \overline{p}_{\mathcal{F}}(n) \cdot \wt(b_d|_0) = \sum_{m \in \mathbb{Z}}m \cdot (\overline{p}_{\mathcal{F}_{m}/ \mathcal{F}_{m+1}}(n) - \overline{p}_{\mathcal{F}}(n)) \cdot \rk_{\mathcal{F}_m / \mathcal{F}_{m+1}}\]
One can also give the following alternative formula.

\begin{prop} \label{prop: computation of weight for torsion-free sheaves}
Let $k$ be a field over $S$. Let $\mathcal{F}$ be a pure sheaf of dimension $d$ on $X_k$. Let $\left(\mathcal{F}_m\right)_{m \in \mathbb{Z}}$ be a filtration of $\mathcal{F}$ in the stack $\Coh^{d}(X)$. Then, we have
\[\wt \left( L_n|_{0}\right) = \sum_{m \in \mathbb{Z}} \left(\overline{p}_{\mathcal{F}_m}(n)- \overline{p}_{\mathcal{F}}(n) \right) \cdot \rk_{\mathcal{F}_m}. \]
\end{prop}
\begin{proof}
We have seen that $\wt \left( L_n|_{0}\right) = \sum_{m \in \mathbb{Z}}m \cdot \left(P_{\mathcal{F}_{m}/ \mathcal{F}_{m+1}}(n) - \frac{\rk_{\mathcal{F}_{m}/\mathcal{F}_{m+1}}}{\rk_{\mathcal{F}}} P_{\mathcal{F}}(n) \right)$. By condition (c) above, we can express this as a finite sum
\[ \wt \left( L_n|_{0}\right) = \sum_{m =-N}^N m \cdot \left(P_{\mathcal{F}_{m}/ \mathcal{F}_{m+1}}(n) - \frac{\rk_{\mathcal{F}_{m}/\mathcal{F}_{m+1}}}{\rk_{\mathcal{F}}} P_{\mathcal{F}}(n) \right)\]
Here $N$ is a big positive integer, and we have $\mathcal{F}_N = 0$ and $\mathcal{F}_{-N} = \mathcal{F}$. By summation by parts and additivity of the Hilbert polynomials, we can rewrite the sum above as
\[ \wt \left( L_n|_{0}\right) = \sum_{m= -N+1}^N \left(P_{\mathcal{F}_m}(n) - \frac{\rk_{\mathcal{F}_m}}{\rk_{\mathcal{F}}}P_{\mathcal{F}}(n)\right) \]
We can express this in terms of reduced Hilbert polynomials to conclude that
\[ \wt \left( L_n|_{0}\right) =  \sum_{m \in \mathbb{Z}} \left(\overline{p}_{\mathcal{F}_m}(n)- \overline{p}_{\mathcal{F}}(n) \right) \cdot \rk_{\mathcal{F}_m} \]
\end{proof}
 
 Before defining our numerical invariant, we need to specify a rational quadratic norm on graded points of $\Coh^{d}(X)$.
\begin{defn} \label{defn: norm filtrations torsion-free sheaves}
Let $g: (B\mathbb{G}_m)^q_k \rightarrow \Coh^{d}(X)$ be a $\mathbb{Z}^q$-graded pure sheaf $\overline{\mathcal{F}} = \bigoplus_{\vec{m} \in \mathbb{Z}^q} \overline{\mathcal{F}}_{\vec{m}}$ of dimension $d$ on $X_k$. We define $b(g)$ to be the positive definite rational quadratic form on $\mathbb{R}^q$ given by
\[b(g)(v) \vcentcolon = \sum_{\vec{m} \in \mathbb{Z}^q} \rk_{\overline{\mathcal{F}}_{\vec{m}}} \cdot (\vec{m} \cdot_{\text{std}} v)^2\]
Here $\cdot_{\text{std}}$ denotes the standard inner product on $\mathbb{R}^q$.
\end{defn}
If $q=1$ above, then the rational quadratic form $b(g)$ on $\mathbb{R}$ is uniquely determined by its value at $1$. From now on, it will be convenient to abuse notation and write $b(g)$ to denote the value $b(g)(1)$ for any graded point $g: (B\mathbb{G}_m)_k \rightarrow \Coh^{d}(X)$.

Given the sequence of line bundles $L_n$ on $\Coh^{d}(X)$ and the norm on graded points $b$, we can define a numerical invariant $\nu$ as explained in Subsection \ref{subsection: background theta stratifications}. For our discussion we only need to understand the corresponding polynomial $\nu(f)$ assigned to a given nondegenerate filtration $f$. 
\begin{defn} \label{defn: poly numerical invariant on coh}
Let $f: \Theta_k \rightarrow \Coh^{d}(X)$ be a nondegenerate filtration given by $(\mathcal{F}_m)_{m\in \mathbb{Z}}$. We define the numerical invariant $\nu(f)$ to be the polynomial in $\mathbb{R}[n]$ given by
\[ \nu(f) \vcentcolon = \frac{ \wt \left(L_n|_{0}\right)}{\sqrt{b(f|_{0})}} = \frac{ \sum_{m \in \mathbb{Z}} m \cdot \left(\overline{p}_{\mathcal{F}_m/ \mathcal{F}_{m+1}}-\overline{p}_{\mathcal{F}} \right) \cdot \rk_{\mathcal{F}_m/\mathcal{F}_{m+1}}}{\sqrt{\left(\sum_{m\in \mathbb{Z}} \rk_{\mathcal{F}_m/\mathcal{F}_{m+1}} \cdot m^2\right)}}\]
\end{defn}
We use the same formula to define a numerical invariant $\nu$ for each of the moduli stacks described in Subsection \ref{section: related moduli problems}. In each case, the corresponding sequence of line bundles $L_n$ and the rational quadratic norm $b$ are pulled back using the forgetful morphism $\Coh^{d}(X)$.

We end this subsection by recording one simple observation that will be useful in our proof of monotonicity. Namely, we observe that the rational quadratic norm $b$ is well-defined on the stack of rational maps.
\begin{lemma} \label{lemma: quadratic norm descends rational maps}
Let $\mathcal{M} = \Coh^{d}(X), \, \Pair_{\mathcal{A}}^{d}(X)$ or $\Lambda \Coh^{d}(X)$. Let $\varphi_1, \varphi_2: (B \mathbb{G}_m)_{k} \to \mathcal{M}_{rat}$ be two graded points. Let $g_1$ and $g_2$ denote the underlying graded pure sheaves of dimension $d$ on $X_{k}$ corresponding to $\varphi_1$ and $\varphi_2$. If there exists a 2-morphism $\varphi_1 \Rightarrow \varphi_2$, then $b(g_1) = b(g_2)$.
\end{lemma}
\begin{proof}
Suppose that $g_1$ (resp. $g_2$) is represented by the graded sheaf $\bigoplus_{m \in \mathbb{Z}} \overline{\mathcal{F}}_m$ (resp.  $\bigoplus_{m \in \mathbb{Z}} \overline{\mathcal{F}}'_m$) defined on $X_{k}$. By the definition of $b$, it suffices to check that $\rk_{\overline{\mathcal{F}}_m} = \rk_{\overline{\mathcal{F}}'_m}$ for all $m \in \mathbb{Z}$. The 2-morphism $\varphi_1 \Rightarrow \varphi_2$ yields a morphism of graded sheaves $\psi: g_1 \to g_2$ that restricts to an isomorphism away from a principal subscheme $D \hookrightarrow X_{k}$, which is regular for both $g_1$ and $g_2$. This implies that the ranks of corresponding graded terms agree, as desired.
\end{proof}
\end{subsection}

\begin{subsection}{Rational filling conditions for torsion-free sheaves} \label{section: rational filling conditions}
Let $R$ be a complete discrete valuation ring over $S$, with residue field $\kappa$ and fraction field $K$. We refer the reader to Subsection \ref{subsection: notation} for the definitions of $\Theta_{R}$, $\overline{ST}_{R}$. The following is the main condition we will study in this subsection.
\begin{defn} \label{defn: rational filling}
Let $\mathcal{M}$ be one of the stacks defined in \ref{section: related moduli problems}. Set $\mathfrak{X}$ to be either $\Theta_{R}$ of $\overline{ST}_{R}$. We say that $\mathcal{M}$ admits $\mathfrak{X}$ rational filling if the following is satisfied: 
for all morphisms $f: \mathfrak{X}\setminus 0 \rightarrow \mathcal{M}$, there exists a morphism $g: \mathfrak{X} \rightarrow \mathcal{M}_{rat}$ and a 2-commutative diagram of pseudofunctors
\begin{figure}[H]
\centering
\begin{tikzcd}
  \mathfrak{X} \setminus 0 \ar[d, "j", labels=left] \ar[r, "f"] & \mathcal{M} \ar[d] \ar[dl, symbol = \Longrightarrow] \\
  \mathfrak{X} \ar[r, dashrightarrow, "g", labels=below]  & \mathcal{M}_{rat}
\end{tikzcd}
\end{figure}
\noindent where $j$ is the open immersion of stacks $\mathfrak{X} \setminus 0 \hookrightarrow \mathfrak{X}$. Note that the direction of the $2$-morphism matters, because $\cM_{rat}$ is not valued in groupoids.
\end{defn}

\begin{remark}
If we set $Y$ to be $Y_{\Theta_{R}}$ (resp. $Y_{\overline{ST}_{R}}$), as defined in \Cref{notation: Y spaces}, then a $\mathbb{G}_m$-equivariant morphism $Y \setminus 0 \to \cM$ is the same thing as a morphism $\Theta_R \setminus 0 \to \cM$ (resp. $\overline{ST}_{R} \setminus 0 \to \cM$). We will use this alternative description in proofs without further elaboration, and work with $\mathbb{G}_m$-equivariant morphisms from $Y$ and $Y \setminus 0$.
\end{remark}

We start by reducing rational filling conditions to the case of the stack of pure sheaves $\Coh^{d}(X)$.
\begin{lemma} \label{lemma: rational filling conditions for other moduli follows from coh}
Suppose that the stack $\Coh^{d}(X)$ admits $\Theta_{R}$ rational filling (resp. $\overline{ST}_{R}$ rational filling). Then both of the stacks $\Pair_{\mathcal{A}}^{d}(X)$ and $\Lambda\Coh^{d}(X)$ admit $\Theta_{R}$ rational filling (resp. $\overline{ST}_{R}$-rational filling).
\end{lemma}
\begin{proof}
Set $Y$ to be $Y_{\Theta_{R}}$ (resp. $Y_{\overline{ST}_{R}}$), as defined in \Cref{notation: Y spaces}. We denote by $W$ the open complement of $0$ in $Y$. 

\smallskip
\noindent{\textit{Proof for $\Pair^d_{\mathcal{A}}(X)$:}}
\smallskip

Suppose that we are given a $\mathbb{G}_m$-equivariant morphism $W = Y \setminus 0 \rightarrow \Pair_{\mathcal{A}}^{d}(X)$. This consists of a $\mathbb{G}_m$-equivariant $W$-pure sheaf $\mathcal{F}$ of dimension $d$ on $X_{W}$ and a $\mathbb{G}_m$-equivariant morphism $\beta: \mathcal{A}|_{X_{W}} \rightarrow \mathcal{F}$.
Assume that $\Coh^{d}(X)$ satisfies $\Theta_{R}$ rational filling (resp. $\overline{ST}_{R}$ rational filling). Then we can obtain the following:
\begin{enumerate}[(1)]
    \item A $\mathbb{G}_m$-equivariant $Y$-pure sheaf $\mathcal{E}$ of dimension $d$ on $X_{Y}$.
    \item A $\mathbb{G}_m$-stable principal subscheme $D \hookrightarrow X_{Y}$ that is $\mathcal{E}$-regular.
    \item A monomorphism $\psi: \mathcal{E}|_{X_{W}} \hookrightarrow \mathcal{F} $ such that $\psi|_{X_{W} \setminus D_{W}}$ is an isomorphism.
\end{enumerate}
It suffices to show that the morphism $\alpha$ obtained by the composition
\[ \alpha: \mathcal{A}|_{X_{W} \setminus D_{W}} \xrightarrow{\beta|_{X_{W} \setminus D_{W}}} \mathcal{F}|_{X_{W} \setminus D_{W}} \xrightarrow{(\psi|_{X_{W} \setminus D_{W}})^{-1}} \mathcal{E}|_{X_{W} \setminus D_{W}} \]
extends to a $\mathbb{G}_m$-equivariant morphism $\widetilde{\alpha}: \mathcal{A}|_{X_{Y} \setminus D} \rightarrow \mathcal{E}|_{X_{Y} \setminus D}$.

Let $U$ denote the open complement $X_{Y} \setminus D$. Note that the closed fiber $U_{0}$ is cut out by a two term regular sequence coming from $Y$ (either $(\varpi, t)$ or $(s,t)$). Since $\mathcal{E}|_{U}$ is $Y$-flat, the sequence is also a regular sequence for the sheaf $\mathcal{E}|_{U}$. This shows that $\mathcal{E}|_{U}$ has depth at least $2$ at all points of the closed fiber $U_{0}$. \cite[\href{https://stacks.math.columbia.edu/tag/0AV5}{Tag 0AV5}]{stacks-project} implies that $\mathop{\mathcal{H}\! \mathit{om}}\left(\mathcal{A}|_{U}, \, \mathcal{E}|_{U}\right)$ also has depth at least $2$ at all points of $U_{0}$. By \cite[\href{https://stacks.math.columbia.edu/tag/0E9I}{Tag 0E9I}]{stacks-project} applied to the open subscheme $U_{W}$ inside $U$, it follows that the morphism $\alpha$ extends uniquely to a morphism $\widetilde{\alpha}: \mathcal{A}|_{U} \rightarrow \mathcal{E}|_{U}$. Moreover, this section is automatically $\mathbb{G}_m$-equivariant. This can be seen by applying \cite[\href{https://stacks.math.columbia.edu/tag/0E9I}{Tag 0E9I}]{stacks-project} again to the pullbacks of $\mathop{\mathcal{H}\! \mathit{om}}\left(\mathcal{A}|_{U}, \, \mathcal{E}|_{U}\right)$ to the product $\mathbb{G}_m \times U$ under the action and projection morphisms.

\smallskip
\noindent{\textit{Proof for $\Lambda\Coh^d(X)$:}}
\smallskip

Suppose that we are given a $\mathbb{G}_m$-equivariant morphism $W=Y \setminus 0 \rightarrow \Lambda\Coh^{d}(X)$. This consists of a $\mathbb{G}_m$-equivariant $\Lambda|_{X_{W}}$-module $\mathcal{F}$ that is $W$-pure of dimension $d$ on $X_{W}$. If $\Coh^{d}(X)$ satisfies the appropriate rational filling conditions, then we can find a triple $(\mathcal{E}, D, \psi)$ as in the previous case of pairs. Set $U = X_{Y} \setminus D$. The isomorphism $\psi|_{U_{W}} : \mathcal{F}|_{U_{W}} \rightarrow \mathcal{E}|_{U_{W}}$ can be used to equip $\mathcal{E}|_{U_{W}}$ with a $\mathbb{G}_m$-equivariant $\Lambda|_{U_{W}}$-module structure. We are left to show that this can be extended to a $\mathbb{G}_{m}$-equivariant $\Lambda|_{U}$-module structure on $\mathcal{E}|_{U}$. 

By Proposition \ref{prop: characterization of lambda modules}, a $\mathbb{G}_m$-equivariant $\Lambda|_{U_{W}}$-module structure on $\mathcal{E}|_{U_{W}}$ is equivalent to the data of a $\mathbb{G}_m$-equivariant morphism $a: \Lambda_1|_{U_{W}} \otimes \mathcal{E}|_{U_{W}} \rightarrow \mathcal{E}|_{U_{W}}$ satisfying conditions $(\Lambda 1)$ and $(\Lambda 2)$. We just need to check that $a$ extends to a $\mathbb{G}_m$-equivariant morphism $\widetilde{a}: \Lambda_1|_{U} \otimes \mathcal{E}|_{U} \rightarrow \mathcal{E}|_{U}$ satisfying $(\Lambda 1)$ and $(\Lambda 2)$.

Since $\mathcal{E}|_{U}$ has depth at least $2$ at all points of $U_{0}$, we can apply \cite[\href{https://stacks.math.columbia.edu/tag/0AV5}{Tag 0AV5}]{stacks-project} to conclude that $\mathop{\mathcal{H}\! \mathit{om}}\left(\Lambda_1|_{U}\otimes \mathcal{E}|_{U}, \, \mathcal{E}|_{U}\right)$ also has depth at least $2$ at all points of $U_{0}$. Now we can use \cite[\href{https://stacks.math.columbia.edu/tag/0E9I}{Tag 0E9I}]{stacks-project} in the same way as in the case of pairs to conclude that there is a unique $\mathbb{G}_m$-equivariant extension $\widetilde{a}: \Lambda_1|_{U} \otimes \mathcal{E}|_{U} \rightarrow \mathcal{E}|_{U}$ of the action morphism $a$. We are left to check conditions $(\Lambda 1)$ and $(\Lambda 2)$. We showed in the course of the proof of Proposition \ref{prop: stack of lambda modules is artin} that these are closed conditions on the base $Y$. Since they hold over the schematically dense open subset $W \subset Y$, they automatically hold over the whole of $Y$.
\end{proof}

We now proceed to prove the rational filling conditions in the case when the fibers of $X \to S$ are geometrically integral of dimension $d$. In this case pure sheaves of dimension $d$ are the same as torsion-free sheaves.
\begin{lemma} \label{lemma: rational filling for projective space}
Suppose that the morphism $X \to S$ is flat with geometrically integral fibers of dimension $d$. The following stacks admit both $\Theta_{R}$ and $\overline{ST}_{R}$ rational filling.
\begin{enumerate}[(i)]
    \item $\Coh^{d}(X)$.
    \item $\Pair_{\mathcal{A}}^{d}(X)$.
    \item $\Lambda\Coh^{d}(X)$.
\end{enumerate}
\end{lemma}
\begin{proof}
By Lemma \ref{lemma: rational filling conditions for other moduli follows from coh}, it suffices to prove (i). Let $\mathfrak{X}$ denote either $\Theta_{R}$ or $\overline{ST}_{R}$. Set $Y \vcentcolon = Y_{\mathfrak{X}}$. Note that the closed fixed point $0$ is cut out by a regular sequence $(y_1, y_2)$ in $\mathcal{O}_{Y}$. We write $W$ for the open complement of $0$ in $Y$. Let $X_{R}$ denote the pullback under the morphism $\Spec(R) \rightarrow S$. We can further pullback using the structure morphism $Y \rightarrow \Spec(R)$ to obtain $X_{Y}$.

Suppose that we are given a $\mathbb{G}_m$-equivariant morphism $W= Y \setminus 0 \rightarrow \Coh^{d}(X)$, which amounts to a $\mathbb{G}_m$-equivariant $W$-flat family $\mathcal{F}$ of torsion-free sheaves on $X_{W}$. Observe that a principal subscheme of $X_{Y}$ is regular with respect to a nontrivial torsion-free sheaf if and only if all of its $Y$-fibers are Cartier divisors if and only if it is a $Y$-relative Cartier divisor \cite[\href{https://stacks.math.columbia.edu/tag/062Y}{Tag 062Y}]{stacks-project}.
Our goal is to find the following.
\begin{enumerate}[(1)]
    \item A $\mathbb{G}_m$-equivariant $Y$-pure sheaf $\mathcal{E}$ of dimension $d$ on $X_{Y}$.
    \item A $\mathbb{G}_m$-equivariant relative Cartier divisor $D \hookrightarrow X_{Y}$.
    \item A $\mathbb{G}_m$-equivariant morphism $\psi: \mathcal{E}|_{X_{W}} \to \mathcal{F}$ such that $\psi|_{X_{W} \setminus D_{W}}$ is an isomorphism.
\end{enumerate}
This automatically implies that $\psi|_{X_{W}}$ is a monomorphism, because $D_{W}$ is both $\mathcal{E}|_{X_{W}}$-regular and $\mathcal{F}$-regular. 

 Let $j$ denote the open immersion $j : W \hookrightarrow Y$. Let $j_{X}: X_{W} \hookrightarrow X_{Y}$ be the base-change. We observe that the pushforward $\widetilde{\mathcal{F}} \vcentcolon = (j_X)_{*} \, \mathcal{F}$ is $\mathbb{G}_m$-equivariant. $\widetilde{\mathcal{F}}$ is coherent and $Y$-flat by \cite[Lemmas 7.16, 7.17]{alper2019existence}, where we take the base ring $k= R$ and  we set the abelian category $\mathcal{A}$ to be the category $\text{QCoh}(X_{R})$ of quasi-coherent sheaves on $X_{R}$. Note that the hypothesis that the discrete valuation ring $R$ is essentially of finite type in \cite{alper-good-moduli} is not necessary in this case, since the abelian category $\cA = \text{QCoh}(X_{R})$ is Noetherian.

Since $X_{0}$ is integral, there exists some dense open subset $V_{0} \subset X_{0}$ such that $\widetilde{\cF}|_{V_{0}}$ is a free graded $\mathcal{O}_{V_{0}}$-module. Say $\widetilde{\mathcal{F}}|_{V_{0}} \cong \bigoplus_{i \in I} \mathcal{O}_{V_{0}}\langle m_i\rangle$, where $I$ is a finite indexing set and $\langle m_i \rangle$ denotes a shift in the $\mathbb{G}_m$-grading. Set $\mathcal{E} = \bigoplus_{i \in I} \mathcal{O}_{X_{Y}}\langle m_i\rangle$. This is a $\mathbb{G}_m$-equivariant $Y$-flat torsion-free sheaf on $X_{Y}$. 

There exists $n \gg 0$ and a non-zero section $s_0 \in H^{0}\left(X_{0}, \mathcal{O}(n)\right)$ such that $s_0$ vanishes at all points of the closed complement $X_{0} \setminus V_{0}$. By increasing $n$ if needed so that $H^1(X_Y, I_{X_0}(n)) = 0$, we can lift $s_0$ to a $\mathbb{G}_m$-invariant global section $s \in H^{0}(X_{Y}, \mathcal{O}(n))$. The section $s$ is non-zero when restricted to the fiber over every point of $Y$, because the non-vanishing locus of $s$ is open, $\bG_m$-stable, and contains a point 
 in the fiber over the fiber over $0$. Thus the $\bG_m$-equivariant principal subscheme $D \hookrightarrow X_Y$ cut out by $s$ is a relative Cartier divisor.

We denote by $U$ the affine open complement $X_{Y} \setminus D$. There is a $\mathbb{G}_m$-equivariant isomorphism $\varphi: \mathcal{E}|_{U_{0}} \xrightarrow{\sim} \widetilde{\mathcal{F}}|_{U_{0}}$, because $U_{0} \subset V_{0}$ by construction. Consider the exact sequence of graded sheaves
\[ 0 \to (y_1, y_2) \cdot\widetilde{\mathcal{F}}|_{U} \to \widetilde{\mathcal{F}}|_{U} \to \widetilde{\mathcal{F}}|_{U_{0}}\to 0\]
Applying $\text{Hom}(\mathcal{E}|_{U},-)$ we get an exact sequence
\[\text{Hom}\left(\mathcal{E}|_{U}, \widetilde{\mathcal{F}}|_{U}\right) \to \text{Hom}\left(\mathcal{E}|_{U_{0}}, \, \widetilde{\mathcal{F}}|_{U_{0}}\right) \to \text{Ext}^1\left(\mathcal{E}|_{U}, \,(y_1, y_2) \cdot\widetilde{\mathcal{F}}|_{U}\right)\]
The third term vanishes, because $\cE|_U$ is locally free and $U$ is affine, so we can lift $\varphi$ to a $\mathbb{G}_m$-equivariant map $\widetilde{\psi} : \mathcal{E}|_{U} \rightarrow \widetilde{\mathcal{F}}|_{U}$. After possibly replacing $\mathcal{E}$ with a subsheaf $\mathcal{E}(-nD)$, we can extend this further to a morphism $\widetilde{\psi} : \mathcal{E} \rightarrow \widetilde{\mathcal{F}}$. We shall prove the following claim.
\begin{claim}
There exists a $\mathbb{G}_m$-equivariant relative Cartier divisor $D' \hookrightarrow X_{Y}$ such that the restriction $\widetilde{\psi}|_{U \setminus D'}: \mathcal{E}|_{U \setminus D'} \rightarrow \widetilde{\mathcal{F}}|_{U \setminus D'}$ is an isomorphism.
\end{claim}
This claim will conclude the proof of the proposition, by replacing $D$ with the sum $D + D'$ and setting $\psi = \widetilde{\psi}|_{X_{W}}$. In order to show the claim, it suffices to find a $D'$ such that $\widetilde{\psi}|_{U \setminus D'}$ is surjective. This is because the kernel of $\widetilde{\psi}|_{U \setminus D'}$ will then be a torsion-free sheaf of rank $0$, which is therefore $0$.  

Consider the $\mathbb{G}_m$-equivariant cokernel $\mathcal{H}$ of $\widetilde{\psi}|_{U}$. By definition $\mathcal{H}$ is a sheaf on $U$ that is supported at the points where $\widetilde{\psi}|_{U}$ is not surjective. Let $Z_{\mathcal{H}}$ be the $\mathbb{G}_m$-equivariant closed subset of $U$ cut out by the $0$-Fitting ideal $\text{Fit}_0(\mathcal{H})$. It suffices to find a relative Cartier divisor $D' \hookrightarrow X_{Y}$ that contains $Z_{\mathcal{H}}$. By construction, $Z_{\mathcal{H}}$ does not contain any point in the fiber $U_{0}$. Let $\mathcal{O}_{U}$ denote the graded coordinate ring of the affine scheme $U$. The closed subscheme $Z_{\mathcal{H}}$ is cut out by a graded ideal $J$. Since $Z_{\mathcal{H}}$ does not meet the fiber $U_{0}$, we know that $y_1\mathcal{O}_{U} + y_2 \mathcal{O}_{U} + J  =\mathcal{O}_{U}$. This means that we can find elements $u,v \in \mathcal{O}_{U}$ and an element $i \in J$ such that $y_1u + y_2v +i = 1$. After passing to graded components, we can assume that $i \in J$ is homogeneous of degree $0$. This gives a $\mathbb{G}_m$-equivariant section $i: \mathcal{O}_{U} \rightarrow \mathcal{O}_{U}$. Note that $i$ cuts out a principal subscheme of $U$ that does not meet $U_{0}$ and contains $Z_{\mathcal{H}}$. There is some $m \gg 0$ such that we can extend $i$ to a homogeneous section $i_{Y}: \mathcal{O}_{X_{Y}} \rightarrow \mathcal{O}_{X_{Y}}(mD)$. Let $D'$ be the principal subscheme of $X_{Y}$ cut out by $i_{Y}$. The support of $D'_{0}$ is contained in $D_{0}$, so $D'_{0}$ is a Cartier divisor on $\mathbb{P}^d_{0}$. The same argument as the one for $D$ above shows that $D'$ is a relative Cartier divisor. By construction $D'$ contains $Z_{\mathcal{H}}$, so $\psi|_{U \setminus D'}$ is surjective.
\end{proof}
\end{subsection}

\begin{subsection}{Proof of monotonicity via ``infinite dimensional GIT"}
We refer the reader to \Cref{defn: strictly theta monotone and STR monotone} for the relevant definitions of strictly monotone. We prove the main theorems of this section using the affine Grassmannians we have defined. This type of argument is what we refer to as ``infinite dimensional GIT."
\begin{thm} \label{thm: theta monotonicity of stacks torsion-free sheaves}
Suppose that the morphism $X \to S$ is flat with geometrically integral fibers of dimension $d$. Then the invariant $\nu$ (\Cref{defn: poly numerical invariant on coh}) is strictly $\Theta$-monotone and strictly $S$-monotone on the stacks $\Coh^{d}(X)$, $\Pair_{\mathcal{A}}^{d}(X)$, and $\Lambda\Coh^{d}(X)$.
\end{thm}
\begin{proof}
We set $\mathcal{M}$ to be one of the stacks $\Coh^{d}(X)$, $\Pair_{\mathcal{A}}^{d}(X)$ or $\Lambda \Coh^{d}(X)$. Let $R$ be a complete discrete valuation ring with residue field $\kappa$. We set $Y$ to denote either $Y_{\Theta_{R}}$ of $Y_{\overline{ST}_{R}}$. Suppose that we are given a morphism $f: \left[\left( \, Y \setminus 0 \, \right)/ \,\mathbb{G}_m\right] \rightarrow \mathcal{M}$. The rational filling properties proved in Lemma \ref{lemma: rational filling for projective space} yield a $\mathbb{G}_m$-equivariant commutative diagram
\begin{figure}[H]
\centering
\begin{tikzcd}
  Y \setminus 0 \ar[ddr, bend right, "j"] \ar[drr, bend left, "f"]  \ar[dr, "\tau"]& & \\ & \mathcal{M}\, \vec{\times}_{\mathcal{M}_{rat}} Y \ar[d] \ar[r] & \mathcal{M} \ar[d] \\   & Y \ar[r, "g"] & \mathcal{M}_{rat}
\end{tikzcd}
\end{figure}
Here $j$ denotes the open immersion $j: Y\setminus 0 \hookrightarrow Y$. By definition, the right comma fiber product $\mathcal{M}\, \vec{\times}_{\mathcal{M}_{rat}} Y$ is an affine Grassmannian $\text{Gr}_{\mathcal{M}}$ (as in Definition \ref{defn: Affine Grassmannians for related moduli}, replacing $X \to S$ with $X_Y \to Y$). The diagram above can be rewritten as follows.
\begin{figure}[H]
\centering
\begin{tikzcd}
  Y \setminus 0 \ar[ddr, bend right, "j"] \ar[drr, bend left, "f"]  \ar[dr, "\tau"]& & \\ & \text{Gr}_{\mathcal{M}} \ar[d] \ar[r, "Forget"] & \mathcal{M} \ar[d] \\   & Y \ar[r, "g"] & \mathcal{M}_{rat}
\end{tikzcd}
\end{figure}
Since the data used to define $\text{Gr}_{\mathcal{M}}$ is $\mathbb{G}_m$-equivariant, the affine Grassmannian acquires a natural $\mathbb{G}_m$-action such that the structure morphism $\text{Gr}_{X_{Y}, D, \xi} \rightarrow Y$ is $\mathbb{G}_m$-equivariant. This action can be defined explicitly in each moduli problem by pulling back sheaves and their associated structures under the morphism induced by multiplication by an element of $\mathbb{G}_m$. This description shows that each $Y$-projective stratum $\text{Gr}_{\mathcal{M}}^{\leq N, P}$ is $\mathbb{G}_m$-stable. By construction, all of the morphisms in the commutative diagram above are $\mathbb{G}_m$-equivariant.

Since $Y \setminus 0$ is quasi-compact, there is a stratum $\text{Gr}_{\mathcal{M}}^{\leq N, P}$ through which $\tau$ factors. Hence we obtain the following $\mathbb{G}_m$-equivariant commutative diagram.
\begin{figure}[H]
\centering
\begin{tikzcd}
   & \text{Gr}_{\mathcal{M}}^{\leq N, P} \ar[r, "Forget"] \ar[d] & \mathcal{M} \\  \, Y \setminus 0 \,  \ar[ur, "\tau"] \ar[r, symbol = \hookrightarrow ] &  Y  & 
\end{tikzcd}
\end{figure}
Let $\Sigma \subset \text{Gr}_{\mathcal{M}}^{\leq N, \, P}$ denote the schematic closure of $Y \setminus 0$ in $\text{Gr}_{\mathcal{M}}^{\leq N, \, P}$. Note that $\Sigma$ is a reduced $\mathbb{G}_m$-scheme with a natural structure morphism to $Y$. The map $\Sigma \rightarrow Y$ is projective, because  $\text{Gr}_{\mathcal{M}}^{\leq N, \, P}$ is projective over $Y$ and $\Sigma$ is a closed subscheme of $\text{Gr}_{\mathcal{M}}^{\leq N, \, P}$. By construction the morphism $\Sigma \rightarrow Y$ is $\mathbb{G}_m$-equivariant and restricts to an isomorphism over $Y \setminus 0$. The composition $\Sigma \rightarrow \text{Gr}_{\mathcal{M}}^{\leq N, \, P} \rightarrow \mathcal{M}$ restricts to $f: Y \setminus 0 \rightarrow \mathcal{M}$ over the open subset $Y\setminus 0 \, \subset \, \Sigma$. Since everything is $\mathbb{G}_m$-equivariant, we obtain a morphism $\widetilde{\varphi}: \left[ \, \Sigma / \, \mathbb{G}_m \, \right] \rightarrow \mathcal{M}$ satisfying condition $(M 2)$ in Definition \ref{defn: strictly theta monotone and STR monotone}.

We are left to check condition $(M3)$ in \Cref{defn: strictly theta monotone and STR monotone}. By \Cref{coroll: ampleness line bundle affine grassmannian}, there exists some $m \gg 0$ such that $L_n^{\vee}|_{\text{Gr}_{\mathcal{M}}^{\leq N, P}}$ is $Y$-ample for all $n \geq m$. For any $a>1$, any field $\kappa$, and any equivariant line bundle on $\cL$ pm $\bP^1_\kappa[a]$, one has $\wt(\cL|_0)-\wt(\cL|_\infty) = a \deg(\cL)$, where $0 := \lim_{t \to 0} t\cdot x$ for a general point $x$ and $\infty$ is the other $\bG_m$-fixed point. In particular, for any finite $\bG_m$-equivariant morphism $\bP^1_\kappa[a] \to \Sigma_0$, $\wt(L_n|_\infty)>\wt(L_n|_0)$ for all $m \geq n$.
%Choose $a \geq 1$. Consider $\mathbb{P}^1_{\kappa}$ equipped with the $\mathbb{G}_m$-action determined by the equation $t \cdot [x:y] = [t^{-a}x: y]$. Suppose that we are given a finite $\mathbb{G}_m$-equivariant morphism $\mathbb{P}^1_{\kappa} \to \Sigma_{0}$. For all $n \geq m$, the restriction $L_n^{\vee}|_{\mathbb{P}^{1}_{\kappa}}$ is ample. This means that $L_n|_{\mathbb{P}^{1}_{\kappa}} \cong \mathcal{O}_{\mathbb{P}^1_{\kappa}}(-N_n)$ for some $N_n >0$. Note that $\mathcal{O}_{\mathbb{P}_{\kappa}^1}(-N_{n})$ admits a unique $\mathbb{G}_m$-equivariant structure up to twisting by a character, and we have $\wt \,\mathcal{O}_{\mathbb{P}^1_{\kappa}}(-N_{n})|_{\infty} > \wt \,\mathcal{O}_{\mathbb{P}^1_{\kappa}}(-N_{n})|_{0}$. 
%We have just shown that $\wt \, (L_n|_{\infty}) > \wt \, (L_n|_{0})$ for $n\gg0$.
%Hence, we will be done if we can prove $b\left(\widetilde{\varphi}|_{ \left[ \infty  / \mathbb{G}_m\right]} \right) =b\left(\widetilde{\varphi}|_{ \left[ 0 / \mathbb{G}_m\right]} \right)$.
The commutative diagram
\[
\begin{tikzcd}
  \left[\mathbb{P}^1_{\kappa}/ \, \mathbb{G}_m \right] \ar[d] \ar[r] & \mathcal{M} \ar[d] \ar[dl, symbol = \Longrightarrow] \\
  (B\mathbb{G}_m)_{\kappa}  \ar[r, "g_{0}"] & \mathcal{M}_{rat}
\end{tikzcd}
\]
shows that there exists a 2-morphism between the compositions $\left[ \infty  / \mathbb{G}_m\right] \xrightarrow{\widetilde{\varphi}} \mathcal{M} \to \mathcal{M}_{rat}$ and $\left[ 0  / \mathbb{G}_m\right] \xrightarrow{\widetilde{\varphi}} \mathcal{M} \to \mathcal{M}_{rat}$. By \Cref{lemma: quadratic norm descends rational maps}, this implies that  $b\left(\widetilde{\varphi}|_{ \left[ \infty  / \mathbb{G}_m\right]} \right) =b\left(\widetilde{\varphi}|_{ \left[ 0 / \mathbb{G}_m\right]} \right)$, and hence for $n\gg 0$,
\[ \nu\left(\widetilde{ \;\varphi}|_{ \left[ \infty  / \mathbb{G}_m\right]}\right) = \frac{\wt \, (L_n|_{\infty})}{\sqrt{b\left(\widetilde{\varphi}|_{ \left[ \infty  / \mathbb{G}_m\right]} \right)}} > \frac{\wt \,( L_n|_{0})}{\sqrt{b\left(\widetilde{\varphi}|_{ \left[ 0  / \mathbb{G}_m\right]} \right)}} = \nu\left(\widetilde{ \;\varphi}|_{ \left[ 0  / \mathbb{G}_m\right]}\right).\]

\end{proof}

We will use \Cref{thm: theta monotonicity of stacks torsion-free sheaves} to bootstrap to the more general case. We shall need the following lemma.
\begin{lemma} \label{lemma: pureness under pushforward finite map}
Let $Y$ be a Noetherian scheme, and let $p: Z \to P$ be a finite morphism of $Y$-schemes. Suppose that $Z$ and $P$ are schemes of finite type over $Y$ with $Y$-fibers of dimension $d$. Then, the pushforward $p_*$ establishes an equivalence between the following two categories:
\begin{enumerate}[(A)]
    \item The category of $Y$-pure sheaves of dimension $d$ on $Z$.
    \item The category of $p_*(\mathcal{O}_{Z})$-modules on $P$ that are $Y$-pure of dimension $d$.
\end{enumerate}
Moreover, if the schemes $Z$ and $P$ admit compatible $\mathbb{G}_m$-actions, then $p_*$ also induces an equivalence of the $\mathbb{G}_m$-equivariant versions of the categories (A) and (B).
\end{lemma}
\begin{proof}
The only thing to check is that a $Y$-flat sheaf $\mathcal{G}$ on $Z$ has $Y$-fibers that are pure of dimension $d$ if and only if the pushforward $p_*(\mathcal{G})$ on $P$ has $Y$-fibers that are pure of dimension $d$. Since the morphism $p$ is affine, the formation of the pushforward $p_*(-)$ commutes with passing to $Y$-fibers. Hence, we can reduce to the case when $p_{y}: Z_{y} \to P_{y}$ is a finite morphism of projective schemes of dimension $d$ over a field $k(y)$.

Let $\mathcal{G}$ be a coherent sheaf on $Z$. Suppose that $p_*(\mathcal{G})$ has an associated point of dimension $\leq d-1$. Let $\mathcal{T}\subset p_*(\mathcal{G})$ be the maximal subsheaf of $p_*(\mathcal{G})$ that is supported on a closed subscheme $H \subset P$ of dimension $\leq d-1$. The $p_*(\mathcal{O}_{Z})$-module generated by $\mathcal{T}$ is also supported on $H$, and so it must coincide with $\mathcal{T}$. Therefore, $\mathcal{T}$ is a $p_*(\mathcal{O}_{Z})$-submodule of $p_*(\mathcal{G})$. The subsheaf of $\mathcal{G}$ corresponding to $\mathcal{T}$ is also supported on dimension $\leq d-1$, because the morphism $p$ is finite. We conclude that $\mathcal{G}$ is not pure of dimension $d$. 

Conversely, if $\mathcal{G}$ is not pure of dimension $d$, then there exists some nontrivial subsheaf $\mathcal{T} \subset \mathcal{G}$ supported in dimension $\leq d-1$. The pushforward $p_*(\mathcal{T}) \subset p_*(\mathcal{G})$ is then supported in dimension $\leq d-1$, and so $p_*(\mathcal{G})$ is not pure of dimension $d$.
\end{proof}

\begin{thm} \label{thm: theta monotonicity of stacks pure sheaves}
The invariant $\nu$ is strictly $\Theta$-monotone and strictly $S$-monotone on the stacks $\Coh^{d}(X)$, $\Pair_{\mathcal{A}}^{d}(X)$, and $\Lambda\Coh^{d}(X)$.
\end{thm}
\begin{proof}
Let $R$ be a complete discrete valuation ring over $S$. Let $\kappa$ be the residue field of $R$ and let $K$ be the fraction field. Let $\mathfrak{X}$ denote either $\Theta_{R}$ or $\overline{ST}_{R}$. Set $Y \vcentcolon = Y_{\mathfrak{X}}$. Let $X_{R}$ denote the pullback under the morphism $\Spec(R) \rightarrow S$. We can further pullback $X_{Y}$ using the structure morphism $Y \rightarrow \Spec(R)$. Let $j : W \hookrightarrow Y$ denote the open complement of $0$ in $Y$. Let $j_{X}: X_{W} \hookrightarrow X_{Y}$ be the base-change. 

\smallskip
\noindent{\textit{Proof for $\Coh^d(X)$:}}
\smallskip
Suppose that we are given $f: \mathfrak{X}\setminus 0 \rightarrow \Coh^{d}(X)$. The morphism $f$ amounts to a $\mathbb{G}_m$-equivariant $W$-pure sheaf $\mathcal{F}$ of dimension $d$ on $X_{W}$. The pushforward $\widetilde{\mathcal{F}} \vcentcolon = (j_X)_{*} \, \mathcal{F}$ is a $\mathbb{G}_m$-equivariant $Y$-flat coherent sheaf by \cite[Lemmas 7.16, 7.17]{alper2019existence}. Let $Z \subset X_{Y}$ denote the $\mathbb{G}_m$-equivariant subscheme of $X_{Y}$ cut out by the $0$th Fitting ideal of $\widetilde{\mathcal{F}}$. Note that the $Y$-fibers of $Z$ have dimension $d$. We can view $\mathcal{F}$ as a $\mathbb{G}_m$-equivariant $W$-pure sheaf of dimension $d$ on $Z_{W}$.

Embed $X_{R}$ into some projective space $\mathbb{P}^N_{R}$ by using a multiple $\mathcal{O}_{X_R}(M):= \mathcal{O}_{X}(M)|_{X_R} $ of the ample line bundle $\mathcal{O}_{X_R}(1)$. Consider the $d$-dimensional closed subvariety $Z_{0} \subset \mathbb{P}^N_{\kappa}$. After replacing $R$ with a finite extension, we can assume that there exists a linear subspace $L_{\kappa} \subset \mathbb{P}^N_{\kappa}$ of dimension $N-d-1$ such that $L_{\kappa} \cap Z_{0}$ is empty. We lift $L_\kappa$ to a linear subspace $L_R \subset \bP^N_R$, and denote by $L_{Y} \subset \mathbb{P}^N_{Y}$ the base-change. The scheme theoretic intersection $L_{Y} \cap Z$ is a $\mathbb{G}_m$-equivariant proper $Y$-scheme. The image of the projection $L_{Y} \cap Z \to Y$ is a $\mathbb{G}_m$-stable closed subset of $Y$ that does not contain $0$, since $Z_{0} \cap L_{Y} = \emptyset$ by construction. We conclude that the image is empty, and therefore $Z$ does not intersect $L_{Y}$.

Consider the affine projection $\mathbb{P}^N_{R} \setminus L_{R} \to \mathbb{P}^{d}_{R}$. We base-change to $Y$ and consider the composition
\[ p: Z \hookrightarrow \mathbb{P}^N_{Y} \setminus L_{Y} \to  \mathbb{P}^{d}_{Y} \]
The morphism $p$ is finite, because it is proper and affine. Let $\Phi$ be the $\mathbb{G}_m$-equivariant sheaf of algebras $p_*(\mathcal{O}_{Z})$ on $\mathbb{P}^d_{Y}$. We regard $\Phi$ as a ring of differential operators by setting $\Phi_0 = \text{Im}\left(\mathcal{O}_{\mathbb{P}^d_{Y}} \to p_*(\mathcal{O}_{Z})\right)$ and $\Phi_1 = \Phi$.

By Lemma \ref{lemma: pureness under pushforward finite map}, $(p|_{Z_{W}})_*(\mathcal{F})$ is a $\mathbb{G}_m$-equivariant $W$-pure $\Phi|_{\mathbb{P}^{d}_{W}}$-module of dimension $d$ on $\mathbb{P}^{d}_{W}$. Lemma \ref{lemma: rational filling for projective space} applied to $X = \mathbb{P}^{d}_{S}$ shows that there exists
\begin{enumerate}[(1)]
    \item A $\mathbb{G}_m$-equivariant $Y$-pure sheaf $\mathcal{E}$ of dimension $d$ on $\mathbb{P}^d_{Y}$.
    \item A $\mathbb{G}_m$-equivariant relative Cartier divisor $D \hookrightarrow \mathbb{P}^d_{Y}$.
    \item A $\mathbb{G}_m$-equivariant morphism $\beta: \Phi_1|_{\mathbb{P}^d_{Y} \setminus D} \otimes \mathcal{E}|_{\mathbb{P}^d_{Y}\setminus D} \to \mathcal{E}|_{\mathbb{P}^d_{Y}\setminus D}$ that equips $\mathcal{E}|_{\mathbb{P}^d_{Y} \setminus D}$ with the structure of a $\Phi|_{\mathbb{P}^d_{Y} \setminus D}$-module.
    \item A $\mathbb{G}_m$-equivariant morphism $\psi: \mathcal{E}|_{\mathbb{P}^{d}_{W}} \to (p|_{Z_{W}})_*\mathcal{F}$ such that $\psi|_{\mathbb{P}^d_{W} \setminus D_{W}}$ is an isomorphism of $\Phi|_{\mathbb{P}^d_{W} \setminus D_{W}}$-modules.
\end{enumerate}
Observe that the $\mathbb{G}_m$-equivariant  $\Phi|_{\mathbb{P}_{W}^d}$-module $(p|_{Z_{W}})_*\mathcal{F}$ can be interpreted as a $\mathbb{G}_m$-equivariant morphism $\tau: W \to \text{Gr}_{\Phi \Coh^{d}(\mathbb{P}^d_{Y}), \mathbb{P}_Y^d, D, \mathcal{E}, \beta}$ over $Y$.

By Lemma \ref{lemma: pureness under pushforward finite map}, there is a well-defined $\mathbb{G}_m$-equivariant morphism $\text{Gr}_{\Phi \Coh^{d}(\mathbb{P}^d_{Y}), \mathbb{P}_Y^d, D, \mathcal{E}, \beta} \to \Coh^{d}(X_{Y})$ over $Y$ that sends a $T$-point $(\mathcal{G}, \varphi, \widetilde{\beta})$ in $\text{Gr}_{\Phi \Coh^{d}(\mathbb{P}^d_{Y}), \mathbb{P}_Y^d, D, \mathcal{E}, \beta}(T)$ to the $T$-pure sheaf on $Z_{T} \subset X_{T}$ corresponding to the $\Phi|_{\mathbb{P}^d_{T}}$-module $(\mathcal{G}, \widetilde{\beta})$. We summarize all of the data we have obtained so far in the following $\mathbb{G}_m$-equivariant commutative diagram:
\begin{figure}[H]
\centering
\begin{tikzcd}
   Y \setminus 0 \,  \ar[r, "\tau"] \ar[rr, bend left, "f"] \ar[dr] & \text{Gr}_{\Phi \Coh^{d}(\mathbb{P}^d_{Y}), \mathbb{P}_Y^d, D, \mathcal{E}, \beta} \ar[r] \ar[d] & \Coh^d(X_{Y}) \\  \,  &  Y  & 
\end{tikzcd}
\end{figure}
Now we can take the scheme closure $\Sigma$ of $Y\setminus 0$ inside $\text{Gr}_{\Phi \Coh^{d}(\mathbb{P}^d_{Y}), \mathbb{P}_Y^d, D, \mathcal{E}, \beta}$. Recall that $M$ denotes the multiple $\mathcal{O}_{X_R}(N)$ that we used to embed $X_R$ into a projective space. The rational line bundles $L_n$ on $\text{Gr}_{\Phi \Coh^{d}(\mathbb{P}^d_{Y}), \mathbb{P}_Y^d, D, \mathcal{E}, \beta}$ coming from $\Coh^{d}(X_{Y})$ agree with the ones pulled-back from $\Coh^{d}(\mathbb{P}_{S}^d)$, where $\mathbb{P}^d_{S}$ is equipped with the $\mathbb{Q}$-ample polarization $\frac{1}{M}\mathcal{O}_{\mathbb{P}^d_{S}}(1)$. It follows from Proposition \ref{prop: line bundle affine grassmannian torsion-free is asymptotically ample} that $L^{\vee}_n$ is eventually relatively ample on each projective stratum of $\text{Gr}_{\Phi \Coh^{d}(\mathbb{P}^d_{Y}), \mathbb{P}_Y^d, D, \mathcal{E}, \beta}$. We can therefore apply the same argument as in Theorem \ref{thm: theta monotonicity of stacks torsion-free sheaves} to conclude in this case.

\medskip
\noindent{\textit{Proof for either $\mathcal{M} = \Pair_{\mathcal{A}}^{d}(X)$ or $\mathcal{M} = \Lambda \Coh^{d}(X)$:}}
\medskip

Choose a $\mathbb{G}_m$-equivariant morphism $W= Y \setminus 0 \rightarrow \mathcal{M}$, which amounts to a $\mathbb{G}_m$-equivariant $W$-pure sheaf $\mathcal{F}$ of dimension $d$ on $X_{W}$ and some extra structure $\alpha$ in the form a morphism of sheaves with target $\mathcal{F}$. We can ignore the extra data $\alpha$ at first and apply the same argument as in the previous case to obtain the data (1)-(4) above.

Let $C$ denote the principal subscheme of $Z$ given by the inverse image of $p^{-1}(D)$. $C$ is given by a $\mathbb{G}_m$-equivariant global section $s$ of $\mathcal{O}_{Z}(n)$. After replacing $D$ and $p^{-1}(D)$ with some multiple and scaling $n$ accordingly, we can lift this to a $\mathbb{G}_m$-equivariant global section of $\mathcal{O}_{X_{Y}}(n)$. We denote by $\widetilde{D}$ the $\mathbb{G}_m$-invariant principal subscheme of $X_{Y}$ cut out by $\widetilde{s}$. By construction the intersection $\widetilde{D} \cap Z$ has dimension $\leq d-1$, and hence $\widetilde{D}$ is regular with respect to any $Y$-pure sheaf of dimension $d$ supported inside $Z$. By Lemma \ref{lemma: pureness under pushforward finite map}, we can view the $\Phi_{\mathbb{P}^d_{Y} \setminus D}$-module $\mathcal{E}|_{\mathbb{P}^d_{Y} \setminus D}$ as a $\mathbb{G}_m$-equivariant $Y$-pure sheaf on $Z \setminus p^{-1}(D)$. Since $Z \setminus p^{-1}(D) = Z \setminus \widetilde{D}$ is a closed subscheme of $X_{Y} \setminus \widetilde{D}$, we can also view $\mathcal{E}|_{\mathbb{P}^d_{Y} \setminus D}$ as a $\mathbb{G}_m$-equivariant $Y$-pure sheaf on $X_{Y} \setminus \widetilde{D}$.

By Lemma \ref{lemma: pureness under pushforward finite map}, for any $T$-point $(\mathcal{G}, \varphi, \widetilde{\beta})$ in $\text{Gr}_{\Phi \Coh^{d}(\mathbb{P}^d_{Y}), \mathbb{P}_Y^d, D, \mathcal{E}, \beta}$ we can view $\mathcal{G}$ as a $T$-pure sheaf on $X_{T}$ supported inside $Z_{T}$. The isomorphism $\varphi|_{X_{T} \setminus \widetilde{D}_{T}}$ can then be used to transport the extra structure $\widetilde{\alpha}$ to the restriction $\mathcal{G}|_{X_{T} \setminus \widetilde{D}_{T}}$. Let $H$ denote the subfunctor of $\text{Gr}_{\Phi \Coh^{d}(\mathbb{P}^d_{Y}), \mathbb{P}_Y^d, D, \mathcal{E}, \beta}$ consisting of $T$-points such that this extra structure $\widetilde{\alpha}$ on $\mathcal{G}|_{X_{T} \setminus \widetilde{D}_{T}}$ extends (uniquely) to $X_{T}$. We claim that $H$ is a closed strict-ind subscheme of $\text{Gr}_{\Phi \Coh^{d}(\mathbb{P}^d_{Y}), \mathbb{P}_Y^d, D, \mathcal{E}, \beta}$. Indeed, for any such $T$-point the structure $\widetilde{\alpha}$ is given by a morphism
\[ \widetilde{\alpha} : \mathcal{B} |_{X_{T} \setminus \widetilde{D}_{T}}\to \mathcal{G}|_{X_{T} \setminus \widetilde{D}_{T}}\]
where $\mathcal{B}$ is a finitely presented coherent sheaf on $X_{T}$ (either $\mathcal{B} = \mathcal{A}$ or $\mathcal{B} = \Lambda_1|_{X_{T}} \otimes \mathcal{G}$). There exists some $N \gg 0$ such that the morphism $\widetilde{\alpha}$ extends uniquely to a morphism $\xi: \mathcal{B} \to \mathcal{G}(N\widetilde{D}_{T})$. Consider the composition
\[ \gamma: \mathcal{B} \xrightarrow{\xi} \mathcal{G}(N \widetilde{D}_{T}) \to \mathcal{G}(N \widetilde{D}_{T}) / \mathcal{G}\]
The fiber product $H_{T}$ represents the locus where this $T$-section $\gamma$ of the Hom functor $\text{Hom}(\mathcal{B}, \, \mathcal{G}(N \widetilde{D}_{T}) / \mathcal{G})$ is equal to the $0$ morphism. By Lemma \ref{lemma: representability of hom sheaves}, the Hom functor $\text{Hom}(\mathcal{B}, \, \mathcal{G}(N \widetilde{D}_{T}) / \mathcal{G})$ is represented by a separated scheme over $T$. Therefore the locus $H_{T}$ where $\gamma =0$ is a closed subscheme of $T$, as claimed.
    
By definition, $H$ admits a morphism to $\mathcal{M}$. The triple $(\mathcal{F},\psi, \alpha)$ represents a $\mathbb{G}_m$-equivariant morphism $\tau: Y \setminus 0 \to H$. We get a similar $\mathbb{G}_m$-equivariant commutative diagram as in the proof for $\Coh^{d}(X)$ above:
\begin{figure}[H]
\centering
\begin{tikzcd}
   Y \setminus 0 \,  \ar[r, "\tau"] \ar[rr, bend left, "f"] \ar[dr] & H \ar[r] \ar[d] & \mathcal{M} \\  \,  &  Y  & 
\end{tikzcd}
\end{figure}
By the same reasoning as before, the pullback of $(L_n)^{\vee}$ to each projective stratum of the closed ind-subscheme $H \subset \text{Gr}_{\Phi \Coh^{d}(\mathbb{P}^d_{Y}), \mathbb{P}_Y^d, D, \mathcal{E}, \beta}$ is $Y$-ample for $n \gg 0$. Therefore we can apply the same argument as in Theorem \ref{thm: theta monotonicity of stacks torsion-free sheaves} to conclude.
\end{proof}
\end{subsection}
\end{section}

\begin{section}{Applications to the moduli of $\Lambda$-modules} \label{section: applications lambda modules}
In this section we will derive the main structural results for the stack of $\Lambda$-modules $\Lambda \Coh^d(X)$ using strict monotonicity (\Cref{thm: theta monotonicity of stacks pure sheaves}). We will need the following facts.
\begin{thm} \label{thm: companion paper}
Let $\nu$ denote the numerical invariant on the stack $\Lambda\Coh^d(X)$ pulled back from $\Coh^d(X)$ (Definition \ref{defn: poly numerical invariant on coh}). We have the following:
\begin{enumerate}[(a)]
    \item The numerical invariant $\nu$ satisfies HN boundedness.
    \item The semistable locus $\Lambda\Coh^d(X)^{\nu \dash {\rm ss}}$ consists of p-semistable $\Lambda$-modules (as in \cite[\S3]{simpson-repnfundamental-I}). For each $P \in \mathbb{Q}[n]$, the substack $\Lambda\Coh^d(X)^{\nu \dash {\rm ss}}_{P}$ of semistable $\Lambda$-modules with Hilbert polynomial $P$ is bounded.
\end{enumerate}
\end{thm}
\begin{proof}
This will be shown in the companion paper \cite{companion-paper}.
\end{proof}

\begin{thm} \label{thm: theta stratification pure sheaves}
The numerical invariant $\nu$ on the stack $\Lambda\Coh^d(X)$ pulled back from $\Coh^d(X)$ (Definition \ref{defn: poly numerical invariant on coh}) defines a weak $\Theta$-stratification on $\Lambda\Coh^{d}(X)$.

The semistable locus $\Lambda\Coh^d(X)^{\nu \dash {\rm ss}}$ consists of p-semistable $\Lambda$-modules. For each $P \in \mathbb{Q}[n]$, the substack $\Lambda\Coh^d(X)^{\nu \dash {\rm ss}}_{P}$ of semistable $\Lambda$-modules with Hilbert polynomial $P$ is bounded. If $S$ is defined over $\mathbb{Q}$, then $\Lambda\Coh^d(X)^{\nu \dash {\rm ss}}_{P}$ admits a separated good moduli space.
\end{thm}
\begin{proof}
We need to check the hypotheses of \Cref{thm: theta stability paper theorem}. The numerical invariant $\nu$ is strictly $\Theta$-monotone and strictly $S$-monotone by \Cref{thm: theta monotonicity of stacks pure sheaves}. Moreover, $\nu$ satisfies HN boundedness by \Cref{thm: companion paper}(a). By \Cref{thm: theta stability paper theorem}(1), this implies that $\nu$ defines a weak $\Theta$-stratification on $\Lambda\Coh^d(X)$.

On the other hand, we know that each open and closed substack $\Lambda\Coh^d(X)^{\nu \dash {\rm ss}}_{P}$ of the $\nu$-semistable locus is bounded. Therefore, if $S$ is defined over $\mathbb{Q}$, if follows by \Cref{thm: theta stability paper theorem}(2) that $\Lambda\Coh^d(X)^{\nu \dash {\rm ss}}_{P}$ admits a separated good moduli space.
\end{proof}
\end{section}

In \cite{companion-paper} we will provide the necessary results (Theorem \ref{thm: companion paper}) for the construction of the $\Theta$-stratification on $\Lambda\Coh^d(X)$ using the numerical invariant $\nu$. We shall also describe this stratification in terms of the Harder-Narasimhan stratification in the context of p-stability. Our canonical filtrations are coarser than the Gieseker-Harder-Narasimhan filtration described in \cite[\S3]{simpson-repnfundamental-I}. We call our filtration of an unstable sheaf the \emph{leading term filtration}. Theorem \ref{thm: theta stratification pure sheaves} gives an alternative proof of the existence of these filtrations that does not use Harder-Narasimhan theory. The Gieseker-Harder-Narasimhan filtration can be recovered from the leading term filtration by iterating the construction for the associated graded sheaves. This is explained in the paper \cite{rho-sheaves-paper} in the more general context of $\rho$-sheaves.

\begin{section}{Applications to moduli of pairs} \label{section: pairs}
In this section we define a family of Laurent polynomial numerical invariants $\nu^{(\delta)}$ on $\Pair_{\mathcal{A}}^{d}(X)$. They are indexed by a choice of rational Laurent polynomial $\delta \in \mathbb{Q}[n, n^{-1}]$. We show that each numerical invariant $\nu^{(\delta)}$ induces a $\Theta$-stratification on $\Pair_{\mathcal{A}}^{d}(X)$. If $\delta \geq 0$, then the corresponding notion of stability agrees with the ones considered in \cite{le-potier-coherent-systems}, \cite{wandel-moduli-pairs}.
\begin{subsection}{Numerical invariants on $\Pair_{\mathcal{A}}^{d}(X)$}
We start by describing filtrations of objects in $\Pair_{\mathcal{A}}^{d}(X)$. Let $k$ be a field over $S$. Let $(\mathcal{F}, \beta)$ be a pair in $\Pair_{\mathcal{A}}^{d}(X)(k)$, where $\mathcal{F}$ is a pure sheaf  of dimension $d$ on $X_{k}$ and $\beta$ is a morphism $\beta: \mathcal{A}|_{X_{k}} \rightarrow \mathcal{F}$. A filtration $f$ of $(\mathcal{F}, \beta)$ consists of a $\mathbb{Z}$-indexed filtration $(\mathcal{F}_{m})_{m\in \mathbb{Z}}$ of the pure sheaf $\mathcal{F}$ (as in Subsection \ref{section: numerical invariant}) such that $\beta: \mathcal{A}|_{X_{k}} \rightarrow \mathcal{F}$ factors through $\mathcal{F}_m$ for all $m\leq 0$. In other words, $\cF_0$ must contain the image of $\beta$. For any such filtration, the associated graded point $f|_{0}$ is the pair $(\text{gr}(\mathcal{F}), \text{gr}(\beta))$ on $X_{k}$ consisting of the graded sheaf $\text{gr}(\mathcal{F}) = \bigoplus_{m \in \mathbb{Z}} \mathcal{F}_m / \mathcal{F}_{m+1}$ and the homogeneous morphism $\text{gr}(\beta) : \mathcal{A}|_{X_{k}} \rightarrow \text{gr}(\mathcal{F})$ given by the composition
\[ \text{gr}(\beta): \mathcal{A}|_{X_{k}} \xrightarrow{\beta} \mathcal{F}_0 \rightarrow \mathcal{F}_0/ \mathcal{F}_1 \hookrightarrow \bigoplus_{m \in \mathbb{Z}} \mathcal{F}_m / \mathcal{F}_{m+1} \]

In the case of the moduli of pairs, it is useful to have some variations of the line bundle $L_n$.
\begin{defn} \label{defn: family line bundles pairs}
Let $\delta \in \mathbb{Q}[n,n^{-1}]$ be a Laurent polynomial in $n$ with rational coefficients. For any $n \in \mathbb{Z}$, we define the line bundle $L_n^{(\delta)}$ on each open and closed substack  $\Pair^d_{\mathcal{A}}(X)_{P} \subset \Pair^{d}_{\mathcal{A}}(X)$ given as follows.
\[ L_n^{(\delta)} \vcentcolon = M_n|_{\Pair^d_{\mathcal{A}}(X)_{P}} \otimes \left(b_d|_{\Pair^d_{\mathcal{A}}(X)_{P}}\right)^{\otimes \,\left(-\frac{\delta(n)}{\rk} - \overline{p}(n)\right)} \]
Here we are using the forgetful morphism $\Pair^{d}_{\mathcal{A}}(X)_{P} \to \Coh^{d}(X)_{P}$ to pull back the line bundles, and the locally constant functions $\overline{p}(n)$ and $\rk$ on $|\Coh^d(X)|$ are as defined in Subsection \ref{subsection: hilbert polynomials background}. %Recall that the reduced Hilbert polynomial $\overline{p}$ on $\Coh^{d}(X)_{P}$ is defined to be the unique scalar multiple of $P$ with leading coefficient $\frac{1}{d!}$, and we define the rank $\rk$ by the equality $P = \rk \cdot \overline{p}$.
\end{defn}
\begin{remark}
Note that $L_n^{(0)} = L_n$ for all $n$.
\end{remark}
For any fixed $\delta$, we denote by $\nu^{(\delta)}$ the Laurent polynomial numerical invariant on $\Pair_{\mathcal{A}}^{d}(X)$ determined by the pullback of the family $L_n^{(\delta)}$ and the rational quadratic norm $b$. This numerical invariant takes values in the group $\mathbb{R}[n,n^{-1}]$, which can be equipped with the structure of a totally ordered $\mathbb{R}$-vector space by defining $p_1 \geq p_2$ if and only if $p_1(n) \geq p_2(n)$ for all $n\gg0$.

We record here the value of $\nu^{(\delta)}$ for filtrations.
\begin{defn} \label{defn: numerical invariant for pairs}
Choose $\delta \in \mathbb{Q}[n, n^{-1}]$. Let $f: \Theta_{k} \to \Pair_{\mathcal{A}}^{d}(X)$ be a filtration given by $f = (\mathcal{F}_m)_{m \in \mathbb{Z}}$. Then $\nu^{(\delta)}(f)$ is the Laurent polynomial given by
\[ \nu^{(\delta)}(f) = \frac{ \sum_{m \in \mathbb{Z}} m \cdot \left(\overline{p}_{\mathcal{F}_m/ \mathcal{F}_{m+1}} - \frac{\delta}{\rk_{\mathcal{F}}} - \overline{p}_{\mathcal{F}}\right) \cdot \rk_{\mathcal{F}_m/\mathcal{F}_{m+1}}}{\sqrt{\left(\sum_{m\in \mathbb{Z}} \rk_{\mathcal{F}_m/\mathcal{F}_{m+1}} \cdot m^2\right)}}\]
\end{defn}
For each fixed $\delta$, the arguments in Subsection \ref{subsection: asymptotic ampleness line bundle pure sheaves} apply without change to the family $L_n^{(\delta)}$. Indeed, the argument in the proof of Lemma \ref{lemma: ampleness for line bundle with big enough index torsion-free grassmannian} shows that the line bundle $b_d$ is torsion on each projective stratum $\Gr_{X,D,\cE}^{\leq N, P} \subset \Gr_{X,D,\cE}$ up to a line bundle coming from the base. Therefore, the family $(L_n^{(\delta)})^\dual \cong L_n^\dual \otimes b_d^{\otimes (\delta(n)/\rk)}$ is eventually relatively ample on each projective stratum of the affine Grassmannian $\text{Gr}_{\Pair_{\mathcal{A}}^{d}(X)}$ associated to some data $\mathcal{E}, D, \beta$. We can use the same argument as in Theorem \ref{thm: theta monotonicity of stacks pure sheaves} to conclude the following.
\begin{thm} \label{thm: monotonicity of family of numerical invariants on pairs}
Fix $\delta \in \mathbb{Q}[n,n^{-1}]$. The numerical invariant $\nu^{(\delta)}$ is strictly $\Theta$-monotone and strictly $S$-monotone on the stack $\Pair_{\mathcal{A}}^{d}(X)$. \qed
\end{thm}
\end{subsection}

\begin{subsection}{Semistable locus and canonical filtrations when $\text{deg}(\delta) \geq d$}
We define the degree $\text{deg}(\delta)$ of a Laurent polynomial $\delta$ to be the maximum among the integers $h$ such that $n^{h}$ has nonzero coefficient in $\delta$.
\begin{prop} \label{prop: description semistable locus pairs big degree}
Fix $\delta \in \mathbb{Q}[n,n^{-1}]$ with $\text{deg}(\delta) \geq d$. Let $(\mathcal{F}, \beta) \in \Pair_{\mathcal{A}}^{d}(X)(k)$ be a field-valued point.
\begin{enumerate}[(i)]
    \item If $\delta <0$, then $(\mathcal{F}, \beta)$ is always $\nu^{(\delta)}$ unstable. There is a unique (up to scaling) canonical filtration $f = (\mathcal{F}_m)_{m \in \mathbb{Z}}$ that maximizes the numerical invariant $\nu^{(\delta)}$. It is given by
    \[\mathcal{F}_m = \begin{cases*} 0 & if $m \geq 2$, \\
    \mathcal{F} & if $m < 2$.
    \end{cases*}\]
    
    \item If $\delta > 0$, then $(\mathcal{F}, \beta)$ is $\nu^{(\delta)}$ semistable if and only if the cokernel of $\beta: \mathcal{A}|_{X_{k}} \to \mathcal{F}$ is supported in dimension $\leq d-1$. If $(\mathcal{F}, \beta)$ is unstable, then up to scaling there is a unique filtration $f = (\mathcal{F}_m)_{m \in \mathbb{Z}}$ maximizing $\nu^{(\delta)}$. It is given by
    \[\mathcal{F}_m = \begin{cases*} 0 & if $m \geq 1$, \\
    \text{Im}(\beta)^{sat} & if $m =0$, \\
    \mathcal{F} & if $0 > m$.
    \end{cases*}\]
    Here the saturation $\text{Im}(\mathcal{\beta})^{sat}$ denotes the smallest subsheaf $\text{Im}(\mathcal{\beta}) \subset \text{Im}(\mathcal{\beta})^{sat} \subset \mathcal{F}$ such that $\mathcal{F}/ \text{Im}(\beta)^{sat}$ is pure of dimension $d$.  
\end{enumerate}
\end{prop}
\begin{proof}
For this proof we set $D \vcentcolon = \text{deg}(\delta)$. 

\noindent{\textit{Proof of (i):}}
\medskip

Let $f$ denote the filtration described in the proposition. If $\delta_{D}$ denotes the leading coefficient of $\delta$, then the leading coefficient of $\nu^{(\delta)}(f)$ is given by $\frac{-\delta_{D}}{\sqrt{\rk_{\mathcal{F}}}}$. This is positive by assumption, and hence $f$ is destabilizing. For any other filtration $f' = (\mathcal{F}'_{m})_{m \in \mathbb{Z}}$ with $\nu^{(\delta)}(f') \geq \nu^{(\delta)}(f)$ we have $\text{deg}(\nu^{(\delta)}(f')) = D$. The leading coefficient of $\nu^{(\delta)}(f')$ is given by 
\[\nu^{(\delta)}(f')_{D} = \frac{ \sum\limits_{m \in \mathbb{Z}} m \cdot \frac{-\delta_{D}}{\rk_{\mathcal{F}}}\cdot \rk_{\mathcal{F}'_m / \mathcal{F}'_{m+1}}}{\sqrt{\left(\sum\limits_{m\in \mathbb{Z}} \rk_{\mathcal{F}'_m/\mathcal{F}'_{m+1}} \cdot m^2\right)}} = \frac{-\delta_D}{\rk_\cF} \cdot \frac{\sum\limits_{m\in \bZ} \left(m \sqrt{\rk_{\cF'_m/\cF'_{m+1}}}\right) \cdot \left( \sqrt{\rk_{\cF'_m/\cF'_{m+1}}}\right)}{\sqrt{\left(\sum\limits_{m\in \mathbb{Z}} \rk_{\mathcal{F}'_m/\mathcal{F}'_{m+1}} \cdot m^2\right)}}\]
To understand this formula in the case $D=d$, we note that the degree $d$ terms in $\overline{p}_{\mathcal{F}'_m / \mathcal{F}'_{m+1}}$ and $\overline{p}_{\mathcal{F}}$ cancel out.

The Cauchy-Schwartz inequality implies that $\nu^{(\delta)}(f')_{D} \leq \frac{-\delta_{D}}{\sqrt{\rk_{\mathcal{F}}}}$, with equality if and only if there is single $m$ for which $\cF'_m/\cF'_{m+1}\neq 0$. In this latter case we have that $f'$ is equal to $f$ up to scaling. Therefore, $f$ is the unique maximizing filtration for the pair $(\mathcal{F}, \beta)$.

%    \item
\medskip
\noindent{\textit{Proof of (ii):}}
\medskip

Suppose that the cokernel of $\beta$ is supported in dimension $\leq d-1$. Let $f = (\mathcal{F}_m)_{m \in \mathbb{Z}}$ be a filtration of $(\mathcal{F}, \beta)$. We shall show that $f$ is not destabilizing. By assumption, $\mathcal{F}$ itself is the only saturated subsheaf of $\mathcal{F}$ containing $\text{Im}(\beta)$. It follows that $\mathcal{F}_m = \mathcal{F}$ for $m \leq 0$. Hence
\[ \nu^{(\delta)}(f) = \frac{ \sum_{m > 0} m \cdot \left(\overline{p}_{\mathcal{F}_m/ \mathcal{F}_{m+1}} - \frac{\delta}{\rk_{\mathcal{F}}} - \overline{p}_{\mathcal{F}}\right) \cdot \rk_{\mathcal{F}_m/\mathcal{F}_{m+1}}}{\sqrt{\left(\sum_{m\in \mathbb{Z}} \rk_{\mathcal{F}_m/\mathcal{F}_{m+1}} \cdot m^2\right)}}\]
Observe that $m \cdot (\overline{p}_{\mathcal{F}_m / \mathcal{F}_{m+1}} - \frac{\delta}{\rk_{\mathcal{F}}} - \overline{p}_{\mathcal{F}}) \cdot \rk_{\mathcal{F}_{m}/\mathcal{F}_{m+1}} \leq 0$ for all $m > 0$, because  $\text{deg}(\delta) \geq d$, $\delta>0$, and the degree $d$ terms of $\overline{p}_\cF$ and $\overline{p}_{\cF_m/\cF_{m+1}}$ cancel whenever $\rk_{\cF_m/\cF_{m+1}} \neq 0$. Therefore, every summand in the numerator is $\leq 0$, and hence $\nu^{(\delta)}(f) \leq 0$.

On the other hand, suppose that the cokernel of $\beta$ is not supported in dimension $\leq d-1$. Then $\text{Im}(\beta)^{sat}$ is a proper subsheaf of $\mathcal{F}$. Let $f$ be the filtration defined in the statement of the proposition. We have
\[\nu^{(\delta)}(f) = \left(\frac{\delta}{\rk_{\mathcal{F}}} + \overline{p}_{\mathcal{F}} - \overline{p}_{\mathcal{F}/ \text{Im}(\beta)^{sat}}\right) \cdot \sqrt{\rk_{\mathcal{F}/ \text{Im}(\beta)^{sat}}}.\]
The leading coefficient $\nu^{(\delta)}(f)_{D} = \delta_{D} \cdot \sqrt{\rk_{\mathcal{F}/ \text{Im}(\beta)^{sat}}} / {\rk_{\mathcal{F}}}$ is positive, so this filtration is destabilizing. We end by showing the following claim: any other filtration $f' = (\mathcal{F}'_m)_{m \in \mathbb{Z}}$ with $\nu^{(\delta)}(f') \geq \nu^{(\delta)}(f)$ must coincide with $f$ up to scaling. This claim implies that $f$ is the unique maximizing filtration of $(\mathcal{F}, \beta)$ up to scaling.

For any such filtration $f'$, we have seen that the leading term of $\nu^{(\delta)}(f')_{D}$ is given by
\[\nu^{(\delta)}(f')_{D} = \frac{ \sum\limits_{m \in \mathbb{Z}} m \cdot \frac{-\delta_{D}}{\rk_{\mathcal{F}}}\cdot \rk_{\mathcal{F}'_m / \mathcal{F}'_{m+1}}}{\sqrt{\left(\sum\limits_{m\in \mathbb{Z}} \rk_{\mathcal{F}'_m/\mathcal{F}'_{m+1}} \cdot m^2\right)}}\]
If some of the $\cF_m$ with positive weight $m>0$ are nonzero, then setting them all to $0$ increases the numerator and decreases the denominator in the formula above. Hence, setting $\cF_m =0$ for all $m>0$ defines another filtration $f''$ with bigger numerical invariant, and equality $\nu^{(\delta)}(f') = \nu^{(\delta)}(f'')$ holds if and only if all nonzero weights $m$ with $\cF_m \neq 0$ were negative to begin with. Therefore, we can reduce to showing the claim in the case when all nonzero weights $m$ of $f'$ are negative. The formula above then reads
\[\nu^{(\delta)}(f')_{D} = \frac{\delta_D}{\rk_\cF} \cdot \frac{\sum\limits_{m<0} \left((-m)\cdot \sqrt{\rk_{\cF'_m/\cF'_{m+1}}}\right) \cdot \left( \sqrt{\rk_{\cF'_m/\cF'_{m+1}}}\right)}{\sqrt{\left(\sum\limits_{m<0} \rk_{\mathcal{F}'_m/\mathcal{F}'_{m+1}} \cdot (-m)^2\right)}}\]
An application of the Cauchy-Schwartz inequality shows that $\nu^{(\delta)}(f')_{D} \leq \frac{\delta_{D}}{\rk_{\mathcal{F}}} \cdot\sqrt{\rk_{\mathcal{F}/\mathcal{F}'_{0}}}$ with equality if and only if there is single $m$ for which $\cF'_m/\cF'_{m+1}\neq 0$ and $m < 0$. Therefore, up to scaling we can reduce to the case when
    \[\mathcal{F}'_m = \begin{cases*} 0 & if $m \geq 1$, \\
    \mathcal{F}'_0 & if $m =0$, \\
    \mathcal{F} & if $0 > m$.
    \end{cases*}\]
and so $\nu^{(\delta)}(f')_{D} = \frac{\delta_{D}}{\rk_{\mathcal{F}}} \cdot\sqrt{\rk_{\mathcal{F}/\mathcal{F}'_{0}}}$. Since $\text{Im}(\beta)^{sat}$ is the smallest saturated subsheaf containing $\text{Im}(\beta)$, we have $\text{Im}(\beta)^{sat} \subset \mathcal{F}'_{0}$. It follows that the inequality
\[\frac{\delta_{D}}{\rk_{\mathcal{F}}} \cdot\sqrt{\rk_{\mathcal{F}/\text{Im}(\beta)^{sat}}} = \nu^{(\delta)}(f)_D \leq \nu^{(\delta)}(f')_D = \frac{\delta_{D}}{\rk_{\mathcal{F}}} \cdot\sqrt{\rk_{\mathcal{F}/\mathcal{F}'_{0}}} \]
implies $\mathcal{F}'_0 = \text{Im}(\beta)^{sat}$, and so $f' = f$.
\end{proof}

\begin{prop}
Let $\delta \in \mathbb{Q}[n, n^{-1}]$ with $\text{deg}(\delta) \geq d$. Suppose that $\delta > 0$. For any $P \in \mathbb{Q}[x]$, let $\Pair_{\mathcal{A}}^{d}(X)^{\nu^{(\delta)} \dash \mathrm{ss}}_{P}$ denote the open substack of $\nu^{(\delta)}$-semistable pairs with Hilbert polynomial $P$. Then, 
\begin{enumerate}[(1)]
    \item $\Pair_{\mathcal{A}}^{d}(X)^{\nu^{(\delta)} \dash \mathrm{ss}}_{P}$ does not depend on the choice of $\delta$.
    \item $\Pair_{\mathcal{A}}^{d}(X)^{\nu^{(\delta)} \dash \mathrm{ss}}_{P}$ is represented by an algebraic space that is proper and of finite presentation over $S$.
\end{enumerate}
\end{prop}
\begin{proof}
Proposition \ref{prop: description semistable locus pairs big degree} (ii) shows that $\Pair_{\mathcal{A}}^{d}(X)^{\nu^{(\delta)} \dash \mathrm{ss}}_{P}$ does not depend on $\delta$. In fact, the explicit description of $\Pair_{\mathcal{A}}^{d}(X)^{\nu^{(\delta)} \dash \mathrm{ss}}_{P}$ shows that it coincides with the quotient husk functor $\text{QHusk}_{P}(\mathcal{A})$ defined by Koll\'ar in \cite{kollar2009hulls}. Part (2) follows from \cite[Thm. 10]{kollar2009hulls}.
\end{proof}
\begin{remark}
When $\mathcal{A} = \mathcal{O}_{X}$ and $d=1$, the space $\Pair_{\mathcal{O}_{X}}^{1}(X)^{\nu^{(\delta)} \dash \mathrm{ss}}_{P}$ described above recovers the moduli space of stable pairs in the sense of Pandharipande-Thomas \cite{pt-stable-pairs}.
\end{remark}

For the rest of this section we will focus on the case when $\text{deg}(\delta) \leq d-1$. These cases yield more interesting $\Theta$-stratifications.
\end{subsection}

\begin{subsection}{HN boundedness for pairs}
As preparation for the proof of HN boundedness, we prove the following lemmas. Recall that $\widehat{\mu}$ is the generalized slope defined in Subsection \ref{subsection: hilbert polynomials background}.

\begin{lemma} \label{lemma: basic boundedness lemma on filtrations}
For any bounded subset $B$ of geometric points in $\Coh^d(X)$ and real number $c \in \bR$, the subset
\[ \mathfrak{S}_{B, c} \vcentcolon = \left\{ \bigoplus_{m \in \mathbb{Z}} \mathcal{F}_m/ \mathcal{F}_{m+1} \; \left| \; \begin{matrix} \exists \cF \in B \text{ such that $(\mathcal{F}_m)_{m \in \mathbb{Z}}$ is a filtration of $\cF$} \; \; \; \\ \text{and $\widehat{\mu}(\mathcal{F}_m/\mathcal{F}_{m+1}) \geq c$ for all $m \in \mathbb{Z}$} \end{matrix} \right. \right\} \]
is bounded.
\end{lemma}
\begin{proof}
Let $T$ be a quasi-compact $S$-scheme and morphism $T \to \Coh^{d}(X)$ such that $B$ is contained in the image. The morphism $T \to \Coh^{d}(X)$ is represented by a $T$-pure sheaf $\mathcal{F}$ of dimension $d$ on $X_{T}$. By base-changing to $X_{T}$ and applying Noetherian approximation \cite[App. C]{thomason-trobaugh}, we can assume without loss of generality that the base scheme $T$ is Noetherian. Note that the rank $\rk_{\mathcal{F}_{t}}$ of the fibers takes only finitely-many values on $T$. This provides a uniform bound for the number of nonzero graded pieces in each graded sheaf $\bigoplus_{m \in \mathbb{Z}} \mathcal{F}_m/\mathcal{F}_{m+1} \in \mathfrak{S}_{B, c}$. We will induct on the maximum $N(B)$ of the number of nonzero graded pieces in elements of $\mathfrak{S}_{B, c}$. The base case $N(B)=1$ is clear, since then $\mathfrak{S}_{B, c}  \subset \{ \mathcal{F}_t \; \mid \; \widehat{\mu}(\mathcal{F}_t) \geq c \; \}$.

We proceed with the induction step. For each element $\bigoplus \mathcal{F}_m/\mathcal{F}_{m+1} \in \mathfrak{S}_{B, c}$, we denote by $\mathcal{F}_{m_{max}}/\mathcal{F}_{m_{max}+1}$ the nonzero graded piece with largest weight $m_{max}$. Notice that then $\mathcal{F}_{m_{max}+1} = 0$, and so $\mathcal{F}_{m_{max}}/\mathcal{F}_{m_{max}+1} = \mathcal{F}_{m_{max}}$ is a subsheaf of a fiber $\mathcal{F}_t$ satisfying $\widehat{\mu}(\mathcal{F}_{m_{max}}) \geq c$. It suffices to show that the set
\[ \left\{ \mathcal{G} \; \mid \; \exists t \in T \text{\text{ with $\mathcal{G} \subset \mathcal{F}_t$}} \; \text{and $\widehat{\mu}(\mathcal{G}) \geq c$} \right\} \]
is bounded, because then we can use the induction hypothesis on the bounded collection $B'$ consisting of all quotients $\mathcal{F}_{t}/\mathcal{F}_{m_{max}}$ (note that by construction $N(B') = N(\mathcal{F})-1$).

After passing to a (finite) Zariski cover of $T$, and replacing $\cO_{X_T}(1)$ with $\cO_{X_T}(M)$ and $c$ with $c/M$ for some sufficiently large $M>0$, we can assume without loss of generality that the ample line bundle $\mathcal{O}_{X_T}(1)$ induces a closed immersion $i : X_T \hookrightarrow \mathbb{P}^n_T$. For any point $t \in T$ and any pure sheaf $\mathcal{G}$ on $X_t$, the pushforward $(i_t)_* \mathcal{G}$ will be a pure sheaf of dimension $d$ on $\mathbb{P}^n_t$ with $\widehat{\mu}((i_t)_* \mathcal{G}) = \widehat{\mu}(\mathcal{G})$. Therefore, we can reduce to the case $X= \mathbb{P}^n$ for the bounded set of sheaves $\{\ (i_t)_{*} \mathcal{F}_{t} \; \mid \; t \in T\}$. Now we can conclude by \cite[Lemma 1.7.9]{huybrechts.lehn}, because the set of Castelnuovo-Mumford regularities $\{ \text{reg}((i_t)_{*} \mathcal{F}_{t}) \; \mid \; t \in T\}$ and the set of Hilbert polynomials $\{ P_{(i_t)_{*} \mathcal{F}_{t}} \; \mid \; t \in T\}$ are both finite by \cite[Lemma 1.7.6]{huybrechts.lehn} (note that the bound on \cite[Lemma 1.7.9]{huybrechts.lehn} only depends on the Hilbert polynomial and regularity of the sheaf).
\end{proof}

\begin{lemma} \label{lemma: near convexity of maximizing slopes}
Assume that $\text{deg}(\delta) \leq d-1$. Let $(\mathcal{F}, \beta)$ be a field-valued point of $\Pair^{d}_{\mathcal{A}}(X)$, and let $f$ be a nondegenerate filtration such that $\nu^{(\delta)}(f)>0$. Then there is another nondegenerate filtration $f'$ with $\nu^{(\delta)}(f') \geq \nu^{(\delta)}(f)$ that corresponds to an unweighted filtration $0 \subsetneq \cG_{(q)} \subsetneq \cG_{(q-1)} \subsetneq \cdots \subsetneq \cG_{(0)} = \cF$ with associated graded sheaves $\overline{\cG}_i := \cG_{(i)}/\cG_{(i+1)}$ such that either
\begin{enumerate}
    \item $\widehat{\mu}(\overline{\cG}_i) = \widehat{\mu}(\cF) + (d-1)!\cdot \delta_{d-1}/\rk_\cF$ for all $i$, or \\
    \item $\widehat{\mu}_{max}(\mathcal{F}) \geq \widehat{\mu}(\overline{\mathcal{G}}_q)>
\cdots > \widehat{\mu}(\overline{\mathcal{G}}_{j+1}) > \widehat{\mu}(\overline{\mathcal{G}}_{j-1}) > \cdots > \widehat{\mu}(\overline{\mathcal{G}}_{0}) \geq \widehat{\mu}_{min}(\mathcal{F})$.
\end{enumerate}
In the second case (2), $j$ denotes the largest index such that $\rm{Im}(\beta) \subset \cG_{(j)}$, and $\widehat{\mu}_{max}(\mathcal{F})$ (resp.  $\widehat{\mu}_{min}(\mathcal{F})$) denotes the maximum (resp. minimum) slope among the graded pieces of the Gieseker-Harder-Narasimhan filtration of $\mathcal{F}$. Moreover,
\[
\widehat{\mu}(\overline{\cG}_j) = \frac{1}{\rk_{\overline{\mathcal{G}}_j}} \cdot \left( \rk_{\mathcal{F}} \cdot \widehat{\mu}(\mathcal{F}) - \sum_{i \neq j}\rk_{\overline{\mathcal{G}}_i} \cdot \widehat{\mu}(\overline{\mathcal{G}}_i) \right).
\]
\end{lemma}

\begin{proof}
Any filtration $f_{\vec{w}}$ of the pair $(\cF,\beta)$ is given by an unweighted filtration 
\[0 \subsetneq \cG_{(q)} \subsetneq \cG_{(q-1)} \subsetneq \cdots \subsetneq \cG_{(0)} = \cF\]
along with a choice of integer weights $\vec{w} := (w_0)_{i=0}^q$ satisfying $w_0 < \cdots < w_q$ and $w_j \geq 0$, where $j$ is the largest index such that $\rm{Im}(\beta) \subset \cG_{(j)}$. By the definition of the numerical invariant $\nu^{(\delta)}$, we have
\[\nu^{(\delta)}(f_{\vec{w}}) = \frac{ \sum_{i=0}^q w_i \cdot \left(\overline{p}_{\overline{\mathcal{G}}_i} - \frac{\delta}{\rk_{\mathcal{F}}} - \overline{p}_{\mathcal{F}}\right) \cdot \rk_{\overline{\mathcal{G}}_i}}{\sqrt{\left(\sum_{i=0}^q \rk_{\overline{\mathcal{G}}_i} \cdot w_i^2\right)}}\]
The Laurent polynomial $\nu^{(\delta)}(f_{\vec{w}})$ has degree at most $d-1$. Let us denote by $\nu_{d-1}(f_{\vec{w}})$ and $\delta_{d-1}$ the coefficient of degree $d-1$ in $\nu^{(\delta)}(f_{\vec{w}})$ and $\delta$ respectively. By the formula above, we have
\begin{equation} \label{eqn: slope invariant}
\nu_{d-1}(f_{\vec{w}}) = \frac{ \sum_{i=0}^q w_i \cdot \left(\frac{\widehat{\mu}(\overline{\mathcal{G}}_i) -\widehat{\mu}(\mathcal{F})}{(d-1)!}- \frac{\delta_{d-1}}{\rk_{\mathcal{F}}} \right) \cdot \rk_{\overline{\mathcal{G}}_i}}{\sqrt{\left(\sum_{i=0}^q \rk_{\overline{\mathcal{G}}_i} \cdot w_i^2\right)}}
\end{equation}

This formula makes sense if we take $\vec{w}$ to be a tuple of real numbers, so we regard it as a function $\nu_{d-1} : \bR^{q+1} \setminus 0 \to \bR$. Note that $\nu_{d-1}$ is continuous and scale invariant. We will use the general fact that for any vector $\vec{u} \in \mathbb{R}^{q+1}$ and positive definite matrix $B$, the function $\nu(\vec{w}) = (\vec{w} \cdot \vec{u}) / \sqrt{\vec{w} \cdot B\vec{w}}$ is strictly quasi-concave on the subset of $\bR^{q+1} \setminus 0$ where $\nu>0$ \cite[Lemma 4.1.15]{halpernleistner2018structure}. By this we mean that for any linearly independent vectors $\vec{w}_1$ and $\vec{w}_2$ such that $\nu(\vec{w}_i)>0$ for $i=1,2$, one has $\nu(t\vec{w}_1+(1-t) \vec{w}_2) > \min(\nu(\vec{w}_1),\nu(\vec{w}_2))$ for $t \in (0,1)$.

Suppose that the filtration $f$ given in the statement has underlying unweighted filtration $\cG_{\bullet}$. We shall show that one can produce the required filtration $f'$ satisfying (1) or (2) by a process deleting steps in the underlying unweighted filtration $\cG_{\bullet}$ and changing the weights $\vec{w}$. Because we are only interested in weights up to scaling by a positive integer, and $\nu^{(\delta)}$ is invariant under scaling, it suffices to construct $f'$ as an $\bR$-weighted filtration then observe a posteriori that $f'$ has rational weights.

If $\nu_{d-1}(f_{\vec{w}}) \leq 0$ for all choices of $\vec{w}$ in the cone of weights
\[
C_{\cG_\bullet} := \{\vec{w} \in \bR^{q+1}| w_0 < \cdots < w_q \text{ and } w_j \geq 0\},
\]
then $\nu^{(\delta)}(f_{\vec{w}})>0$ implies that $\nu_{d-1}(f_{\vec{w}})=0$. In this case, you must have the coefficient of each $w_i$ vanish in \eqref{eqn: slope invariant}, or else you could perturb the $w_i$ so that $\nu_{d-1}(f)>0$. Thus $f$ itself satisfies condition (1) of the lemma.

%$w_j^{*} = Q \cdot (\widehat{\mu}(\overline{\cG}_j) - ({\delta_{d-1}}/{\rk_{\cF}}) - \widehat{\mu}(\cF))$

Otherwise, assume for the moment that there exists a $\vec{w}^{*}$ that maximizes $\nu_{d-1}$ on $C_{\cG_\bullet}$ and such that $\nu_{d-1}(f_{\vec{w}^{*}})>0$. This means that either: (A) $w_j^{*}>0$ and the weights are an unconstrained local max of $\nu_{d-1}$ on all of $\bR^{q+1} \setminus 0$; or (B) the weights are an unconstrained local max of $\nu_{d-1}$ on the subspace $\{w_j=0\} \subset \bR^{q+1}\setminus 0$. Either way, the critical point equations imply that
\begin{equation} \label{eqn: maximizing weights}
w_i^{*} = Q \cdot \left(\frac{\widehat{\mu}(\overline{\cG}_i) - \widehat{\mu}(\cF)}{(d-1)!} - \frac{\delta_{d-1}}{\rk_{\cF}} \right)\end{equation}
for some positive rational number $Q$ and all $i \neq j$. The weight $w_j^{*}$ is either $0$ in case (B) or is given by the same formula \eqref{eqn: maximizing weights} in case (A). Therefore this maximizer occurs at a rational point. It is unique because $\nu_{d-1}$ is strictly quasi-concave. Thus if the original weights of $f$ are not a multiple of this maximizer $\vec{w}^{*}$, then $\nu_{d-1}(f)<\nu_{d-1}(f_{\vec{w}^{*}})$, so we can set $f' = f_{\vec{w}^{*}}$. The inequalities in condition (2) of the lemma now follow from combining the explicit formula for $w_i^{*}$ in terms of $\widehat{\mu}(\overline{G}_i)$ and the hypothesis that $w_0^{*}<\cdots<w_q^{*}$.

In general, $\nu_{d-1}$ does not achieve a maximum on $C_{\cG_\bullet}$. However, if we denote by $\overline{C}_{\cG_\bullet}$ the closure of $C_{\cG_{\bullet}}$ in $\mathbb{R}^{q+1}$, then the continuous function $\nu_{d-1}$ must achieve a maximum on the compact space $(\overline{C}_{\cG_\bullet} \setminus 0) / \bR^\times_{>0}$. In particular, $\nu_{d-1}$ must attain a global maximum at some point of $\overline{C}_{\cG_\bullet} \setminus 0$, which is unique up to scaling because of the strict quasi-concavity of $\nu_{d-1}$. If this maximizer does not lie in $C_{\cG_\bullet}$, then it must lie on one of the boundary components where $w_i = w_{i+1}$. The restriction of $\nu_{d-1}$ to this boundary component can be identified with the formula \eqref{eqn: slope invariant} on the cone of weights $C_{\cG'_\bullet}$, where $\cG'_\bullet$ is the unweighted filtration obtained from $\cG_\bullet$ by deleting the $i^{th}$ step. Hence, after deleting finitely many steps we get an unweighted filtration $\cG'_\bullet$ such that $\nu_{d-1}$ does admit a maximum $\vec{w}^{*}$ on $C_{\cG'_\bullet}$. We have already shown that the resulting filtration $f'_{\vec{w}^{*}}$ satisfies the condition (2) of the lemma. Furthermore, $f'_{\vec{w}^{*}}$ maximizes $\nu_{d-1}$ on $\overline{C}_{\cG_\bullet} \setminus 0$ by construction, so $\nu_{d-1}(f)<\nu_{d-1}(f'_{\vec{w}^{*}})$, which in turn implies $\nu^{(\delta)}(f) \leq \nu^{(\delta)}(f'_{\vec{w}^{*}})$.

The inequalities with $\widehat{\mu}_{min}(\cF)$ and $\widehat{\mu}_{max}(\cF)$ follow from the fact that $\overline{\cG}_0$ and $\overline{\cG_q}$ are pure quotients and subobjects of $\cF$ respectively. Finally, the formula for $\widehat{\mu}(\cG_j)$ follows from the additivity of Hilbert polynomials.
\end{proof}

\begin{prop}[HN-boundedness] \label{prop: HN boundedness of pairs}
Fix $\delta \in \mathbb{Q}[n,n^{-1}]$ with $\text{deg}(\delta) \leq d-1$. Let $T$ be an affine Noetherian scheme. Let $g: T \rightarrow \Pair_{\mathcal{A}}^{d}(X)$. Then there exists a quasi-compact open substack $\mathcal{U}_T$ of $\Pair_{\mathcal{A}}^{d}(X)$ satisfying the following: For all geometric points $t \in T$ and all $\nu^{(\delta)}$-destabilizing filtrations $f$ of $g(t)$, there exists a nondegenerate filtration $f'$ of $g(t)$ with $f'|_{0} \in \mathcal{U}_T$ and $\nu^{(\delta)}(f') \geq \nu^{(\delta)}(f)$.
\end{prop}
\begin{proof}
The natural forgetful map $\text{Forget}: \Pair_{\mathcal{A}}^{d}(X) \rightarrow \Coh^{d}(X)$ is quasi-compact by Proposition \ref{prop: stack pairs algebraic}. Hence, it suffices to show that there is a quasi-compact open substack $\mathcal{W}_{T}$ of $\Coh^{d}(X)$ such that for all destabilizing filtrations $f$ for a pair $g(t)$ as in the statement of this proposition we can find a filtration $f'$ such that $\text{Forget}(f'|_{0}) \in \mathcal{W}_{T}$ and $\nu^{(\delta)}(f') \geq \nu^{(\delta)}(f)$.

The morphism $T \rightarrow \Pair_{\mathcal{A}}^{d}(X)$ is represented by a pair $(\mathcal{F}, \beta)$ consisting of a $T$-pure sheaf $\mathcal{F}$ of dimension $d$ on $X_{T}$ and a morphism $\beta: \mathcal{A}|_{X_{T}} \rightarrow \mathcal{F}$. Lemma \ref{lemma: basic boundedness lemma on filtrations} shows that it is sufficient to find a uniform lower bound $C$ such that for all $t \in T$ and all destabilizing filtrations $f$ of the pair $(\mathcal{F}_{X_t}, \beta_{X_t})$, there exists another filtration $f'= (\mathcal{F}_m)_{m \in \mathbb{Z}}$ satisfying $\nu^{(\delta)}(f') \geq \nu^{(\delta)}$ and $\widehat{\mu}(\mathcal{F}_m/\mathcal{F}_{m+1}) \geq C$ for all $m \in \mathbb{Z}$.

By Lemma \ref{lemma: near convexity of maximizing slopes}, we can always find a filtration $f'$ with $\nu^{(\delta)}(f') \geq \nu^{(\delta)}(f)$ and such that for all $i \neq j$ the associated graded piece $\overline{\mathcal{G}}_i$ satisfies either
\begin{enumerate}[(1)]
    \item $\widehat{\mu}(\overline{\mathcal{G}}_i) = \widehat{\mu}(\mathcal{F}_t) + (d-1)!\cdot \delta_{d-1}/\rk_{\mathcal{F}_t}$, or
    \item $\widehat{\mu}_{max}(\mathcal{F}_t) \geq \widehat{\mu}(\overline{\mathcal{G}}_i) \geq \widehat{\mu}_{min}(\mathcal{F}_t)$.
\end{enumerate}
Since the family $\mathcal{F}_t$ is bounded, either case yields uniform upper and lower bounds for $\widehat{\mu}(\overline{\mathcal{G}}_i)$ if $i \neq j$. The slope of the remaining associated graded piece $\overline{\mathcal{G}}_j$ can be bounded using the formula
\[
\widehat{\mu}(\overline{\cG}_j) = \frac{1}{\rk_{\overline{\mathcal{G}}_j}} \cdot \left( \rk_{\mathcal{F}_t} \cdot \widehat{\mu}(\mathcal{F}_t) - \sum_{i \neq j}\rk_{\overline{\mathcal{G}}_i} \cdot \widehat{\mu}(\overline{\mathcal{G}}_i) \right)
\]
from Lemma \ref{lemma: near convexity of maximizing slopes}, because we have established that all the terms in the right-hand side admit uniform upper and lower bounds.
\end{proof}
\end{subsection}

\begin{subsection}{$\Theta$-stratifications on $\Pair_{\mathcal{A}}^{d}(X)$}
\begin{thm} \label{thm: existence of weak theta stratification for pairs}
Let $\delta \in \mathbb{Q}[n,n^{-1}]$ with $\text{deg}(\delta) \leq d-1$. The invariant $\nu^{(\delta)}$ defines a weak $\Theta$-stratification of the stack $\Pair_{\mathcal{A}}^{d}(X)$. If the scheme $S$ is defined over $\mathbb{Q}$, then $\nu^{(\delta)}$ defines a $\Theta$-stratification.
\end{thm}
\begin{proof}
We use Theorem \ref{thm: theta stability paper theorem} (1). By Theorem \ref{thm: monotonicity of family of numerical invariants on pairs}, the invariant $\nu^{(\delta)}$ is strictly $\Theta$-monotone. On the other hand Proposition \ref{prop: HN boundedness of pairs} implies that $\nu^{(\delta)}$ satisfies the HN-boundedness condition.
\end{proof}

In particular, every unstable pair $(\mathcal{F}, \beta)$ defined over a field admits (after maybe passing to a purely inseparable field extension) a canonical filtration $(\mathcal{F}_m)_{m \in \mathbb{Z}}$ that maximizes the numerical invariant $\nu^{(\delta)}$. Such canonical filtrations are unique up to scaling the weights.

 One interesting feature of the ``non-abelian" moduli problem $\Pair_{\mathcal{A}}^{d}(X)$ is that the canonical filtration $(\mathcal{F}_m)_{m\in \mathbb{Z}}$ of a pair $(\mathcal{F}, \beta)$ need not be convex with respect to the numerical invariant. This is illustrated by the following example.
\begin{example} \label{example: filtration nonconvex pairs}
Set $\delta = 0$. Let $X = \mathbb{P}^1$ and $\mathcal{A} = \mathcal{O}$. Set $\mathcal{F} = \mathcal{O}(5) \oplus \mathcal{O}(1) \oplus \mathcal{O}$. We define $\beta: \mathcal{O} \rightarrow \mathcal{F}$ to be the inclusion into the last component. Then, the canonical filtration $(\mathcal{F}_m)_{m \in \mathbb{Z}}$ of the pair $(\mathcal{F}, s)$ is given up to scaling by
\[ \mathcal{F}_m =  \begin{cases} 0 \; \; \; \;  \; \; \; \;  \;\; \; \; \; \; \; \; \; \; \text{for $m >3$, } \\
\mathcal{O}(5) \; \; \; \; \; \; \; \; \; \; \; \; \text{for $3 \geq m >0$, } \\
\mathcal{O}(5) \oplus \mathcal{O}  \; \; \; \; \; \text{for $0 \geq m >-1$, } \\
\mathcal{F} \; \; \; \; \; \; \; \; \; \; \; \; \;\; \; \; \; \text{for $-1 > m$. }
\end{cases}\]
The associated graded pieces are $\mathcal{O}(5)$ in weight $3$, $\mathcal{O}$ in weight $0$, and $\mathcal{O}(1)$ in weight $-1$.
\end{example}
This phenomenon (i.e. nonconvexity of canonical filtrations) does not arise for the moduli stacks $\Coh^{d}(X)$ and $\Lambda \Coh^{d}(X)$.

\begin{remark}
We note, however, that the canonical filtration of a pair is always ``nearly convex". More precisely, one can modify the proof of Lemma \ref{lemma: near convexity of maximizing slopes} (possibly replacing $\nu_{d-1}$ with a lower order term) to see that the numerical invariant $\nu^{(\delta)}$ of the graded pieces $\mathcal{F}_m/\mathcal{F}_{m+1}$ will form a convex sequence except possibly at $\mathcal{F}_0/\mathcal{F}_1$.
\end{remark}

\begin{remark}
The canonical filtration coming from the $\Theta$-stratification agrees with the Harder-Narasimhan filtration for Bradlow pairs in the case of rank $2$, as defined in \cite{thaddeus-verlinde}. 

In higher rank, we find a definition of Harder-Narasimhan filtrations for Bradlow pairs in \cite[3.3.2]{mochizuki-invariants-surfaces-book}. It is interesting to note that our canonical filtration does not necessarily agree with Mochizuki's Harder-Narasimhan filtration when the rank is bigger than $2$. As an example, take Example \ref{example: filtration nonconvex pairs} for some constant $0 < \delta \ll1$.
\end{remark}

\end{subsection}
\begin{subsection}{Moduli spaces for pairs}
Let $\delta \in \mathbb{Q}[n, n^{-1}]$ with $\text{deg}(\delta) \leq d-1$. Since the invariant $\nu^{(\delta)}$ induces a weak $\Theta$-stratification, it follows that the locus of $\nu^{(\delta)}$-semistable pairs is an open substack of $\Pair_{\mathcal{A}}^{d}(X)$. We denote this open substack by $\Pair_{\mathcal{A}}^{d}(X)^{\nu^{(\delta)} \dash \mathrm{ss}}$. For each rational polynomial $P \in \mathbb{Q}[x]$, we set $\Pair_{\mathcal{A}}^{d}(X)^{\nu^{(\delta)} \dash \mathrm{ss}}_{P}$ to be the open and closed substack of $\Pair_{\mathcal{A}}^{d}(X)^{\nu^{(\delta)} \dash \mathrm{ss}}$ parametrizing $\nu^{(\delta)}$-semistable points such that the underlying pure sheaf has Hilbert polynomial $P$.
\begin{prop} \label{prop: boundedness semistable locus pairs}
Assume that $\text{deg}(\delta) \leq d-1$. The stack $\Pair_{\mathcal{A}}^{d}(X)^{\nu^{(\delta)} \dash \mathrm{ss}}_{P}$ is quasi-compact.
\end{prop}
\begin{proof}
By Proposition \ref{prop: stack pairs algebraic}, the forgetful morphism $\Pair_{\mathcal{A}}^{d}(X) \to \Coh^{d}(X)$ is quasi-compact. Hence, it suffices to show that the image of $\Pair_{\mathcal{A}}^{d}(X)^{\nu^{(\delta)} \dash \mathrm{ss}}_{P} \to \Coh^{d}(X)$ is quasi-compact. Let $k$ be a field and let $(\mathcal{F}, \beta) \in \Pair_{\mathcal{A}}^{d}(X)_{P}^{\nu^{(\delta)} \dash \mathrm{ss}}(k)$. We want to show that $\mathcal{F}$ belongs to a fixed bounded family relative to $S$. By Noetherian approximation and \cite[Thm. 4.4]{langer-boundedness-general}, it suffices to show that there exists a uniform upper bound $\mu_0$ such that every proper saturated subsheaf $\mathcal{E} \subset \mathcal{F}$ satisfies $\widehat{\mu}(\mathcal{E}) \leq \mu_0$.

Let $\delta_{d-1}$ denote the coefficient of $n^{d-1}$ in the Laurent polynomial $\delta$. Note that the slope $\widehat{\mu}(\mathcal{F})$ and rank $\rk_{\mathcal{F}}$ are fixed, since they only depend on $P$. We claim that $\mu_0 = \widehat{\mu}(\mathcal{F}) + \frac{\delta_{d-1}}{\rk_{\mathcal{F}}}$ is a valid upper bound. To see this, let $\mathcal{E} \subset \mathcal{F}$ be a proper saturated subsheaf. Consider the filtration $f= (\mathcal{F}_m)_{m \in \mathbb{Z}}$ of the pair given by
\[    \mathcal{F}_m \vcentcolon = \begin{cases*}
      0 & if $m > 1$, \\
      \mathcal{E}        & if $ m=1$, \\
      \mathcal{F} & if $0 \leq m$.
    \end{cases*}
    \]
Then we have $\nu^{(\delta)}(f) = \sqrt{\rk_{\mathcal{E}}} \cdot (\overline{p}_{\mathcal{E}} - \frac{\delta}{\rk_{\mathcal{F}}} - \overline{p}_{\mathcal{F}})$. The leading coefficient in degree $d-1$ is given by
\[ \nu^{(\delta)}(f)_{d-1} =  \sqrt{\rk_{\mathcal{E}}} \cdot \left(\widehat{\mu}(\mathcal{E}) - \frac{\delta_{d-1}}{\rk_{\mathcal{F}}}- \widehat{\mu}(\mathcal{F})\right)\]
Since $(\mathcal{F}, \beta)$ is semistable, we must have $\nu^{(\delta)}(f)_{d-1} \leq 0$. This implies that $\widehat{\mu}(\mathcal{E}) \leq \widehat{\mu}(\mathcal{F}) + \frac{\delta_{d-1}}{\rk_{\mathcal{F}}}$, as desired.
\end{proof}

In order to check the existence part of the valuative criterion for properness for pairs, we will use the following.
\begin{prop} \label{prop: valuative criterion}
Let $R$ be a complete discrete valuation ring with fraction field $K$. Suppose that we are given a $2$-commutative diagram as follows.
\begin{figure}[H]
\centering
\begin{tikzcd}
  \Spec(K) \ar[r] \ar[d] &  \Coh^{d}(X) \ar[d] \\   \Spec(R) \ar[r] &  S
\end{tikzcd}
\end{figure}
Then there exist a morphism $\Spec(R) \rightarrow \Coh^{d}(X)$ such that the following diagram is $2$-commutative
    \begin{figure}[H]
\centering
\begin{tikzcd}
  \Spec(K) \ar[r] \ar[d] &  \Coh^{d}(X) \ar[d] \\   \Spec(R) \ar[ur] \ar[r] &  S
\end{tikzcd}
\end{figure}
\end{prop}
\begin{proof}
The map $\Spec(K) \rightarrow \Coh^{d}(X)$ amounts to a pure sheaf $\mathcal{F}$ of dimension $d$ on $X_{K}$. By \cite[\href{https://stacks.math.columbia.edu/tag/01PF}{Tag 01PF}]{stacks-project} we can extend it to a coherent sheaf $\widetilde{\mathcal{F}}$ on $X_{R}$. Let $\varpi$ be a uniformizer of $R$ and let $\kappa$ denote the residue field. Let $j: \Spec(K) \hookrightarrow \Spec(R)$ denote the open immersion. Then we can kill the $\varpi$-torsion by taking the image of the unit $\widetilde{\mathcal{F}} \rightarrow j_*j^*\widetilde{\mathcal{F}}$, so we can assume without loss of generality that $\widetilde{\mathcal{F}}$ is $R$-flat. Let $B \subset \widetilde{\mathcal{F}}_\kappa$ be the maximal subsheaf of dimension $<d$, and form $\cE = \ker(\widetilde{\mathcal{F}} \to \widetilde{\mathcal{F}}_\kappa / B)$. We can again replace $\cE$ with $\ker(\cE \to \cE_k / B')$, where $B' \subset \cE_k$ is the maximal subsheaf of dimension $<d$. The proof of \cite[Theorem 2.B.1]{huybrechts.lehn} \footnote{The smoothness hypothesis is unnecessary for the argument. Using the notation in \cite{huybrechts.lehn}, a sheaf is semistable in $\Coh_{d, d-1}$ if and only if it is pure of dimension $d$.} applies verbatim to show that iterating this procedure results in a sheaf $\cE$ with $\cE|_{X_{\kappa}}$ pure and $\mathcal{E}|_{X_{K}} \cong \widetilde{\mathcal{F}}|_{X_{K}}$. This is the morphism $\Spec(R) \rightarrow \Coh^{d}(X)$ that we were looking for.
\end{proof}

\begin{thm} \label{thm: moduli space pairs}
Suppose that $\text{deg}(\delta) \leq d-1$. Choose a Hilbert polynomial $P \in \mathbb{Q}[x]$. Suppose that the scheme $S$ is defined over $\mathbb{Q}$. Then the stack $\Pair_{\mathcal{A}}^{d}(X)^{\nu^{(\delta)} \dash \mathrm{ss}}_{P}$ admits a proper good moduli space over $S$.
\end{thm}
\begin{proof}
We use Theorem \ref{thm: theta stability paper theorem} (2). First, $\nu^{(\delta)}$ is strictly $\Theta$-monotone and strictly $S$-monotone by Theorem \ref{thm: monotonicity of family of numerical invariants on pairs}. HN boundedness follows from Proposition \ref{prop: HN boundedness of pairs}. On the other hand, Proposition \ref{prop: boundedness semistable locus pairs} implies that the stack $\Pair_{\cA}^d(X)_{P}^{\nu^{(\delta)} \dash \mathrm{ss}}$ is quasi-compact.

We are only left to check that $\Pair_{\mathcal{A}}^{d}(X)$ satisfies the existence part of the valuative criterion for properness. Let $R$ be a complete discrete valuation ring over $S$ with fraction field $K$ and uniformizer $\varpi \in R$. Let $(\mathcal{F}, \beta)$ be a pair on $X_{K}$. By Proposition \ref{prop: valuative criterion}, we can extend $\mathcal{F}$ to a $R$-pure sheaf $\widetilde{\mathcal{F}}$ of dimension $d$ on $X_{R}$. Since $\widetilde{\mathcal{F}}$ is $R$-flat, for any $n \geq 0$ we have an inclusion $\widetilde{\mathcal{F}} \subset \varpi^{-n} \cdot \widetilde{\mathcal{F}}$ that restricts to an isomorphism on the generic fiber $X_{K}$. The morphism $\beta: \mathcal{A}|_{X_{K}} \to \mathcal{F}$ extends to a morphism $\widetilde{\beta}: \mathcal{A}_{X_{R}} \to \varpi^{-n} \cdot \widetilde{\mathcal{F}}$ for some $n\gg0$. The resulting pair $(\varpi^{-n} \cdot \widetilde{\mathcal{F}}, \widetilde{\beta})$ on $X_{R}$ extends $(\mathcal{F}, \beta)$, as desired.
\end{proof}

\begin{remark}
The same proof as in Proposition \ref{prop: description semistable locus pairs big degree} (i) shows that the $\nu^{(\delta)}$-semistable locus is empty whenever $\delta <0$. Hence, we can restrict our attention to $\delta \geq 0$ for the purposes of studying the moduli space.
\end{remark}

\begin{example} If $\delta =0$, then a pair $(\mathcal{F}, \beta)$ is $\nu^{(\delta)}$-semistable if and only if the sheaf $\mathcal{F}$ is Gieseker semistable. Indeed, given any destabilizing filtration $f$ of $\mathcal{F}$ we can obtain a filtration $f'$ of the pair $(\mathcal{F}, \beta)$ by shifting the weights so that $\mathcal{F}_0$ contains the image of $\beta$. Since $\nu^{(0)}$ remains unchanged after shifting weights, $f'$ is a destabilizing filtration for the pair $(\mathcal{F}, \beta)$.

Using the scaling action on $\cF$, one can show that every such semistable pair $(\mathcal{F}, \beta)$ contains the semistable pair $(\mathcal{F}, 0)$ in its closure. This shows that the moduli space of $\Pair_\cA^d(X)^{\nu^{(\delta)} \dash \mathrm{ss}}$ agrees with the moduli space of Gieseker semistable pure sheaves on $X$.
\end{example}

Next we give an alternative description of $\nu^{(\delta)}$-semistability in this case when $\delta>0$. This shows that our stability condition is analogous to the notion of $\delta$-stability for coherent systems considered by Le Potier \cite{le-potier-coherent-systems}. See also \cite{wandel-moduli-pairs} for a formulation that is closer to ours. 
\begin{prop} \label{prop: comparison of stability pairs}
Suppose that $\text{deg}(\delta) \leq d-1$. Let $(\mathcal{F}, \beta)$ be a field-valued point of $\Pair_{\mathcal{A}}^{d}(X)$. Let $\delta \in \mathbb{Q}[n,n^{-1}]$ with $\delta > 0$. Then $(\mathcal{F}, \beta)$ is $\nu^{(\delta)}(f)$-semistable if and only if the following two conditions are satisfied
\begin{enumerate}[(a)]
    \item $\beta \neq 0$.
    
    \item All proper saturated subsheaves $\mathcal{E} \subset \mathcal{F}$ satisfy
\[   \begin{cases*}
      \overline{p}_{\mathcal{F}} + \frac{\delta}{\rk_{\mathcal{F}}} \geq \overline{p}_{\mathcal{E}} & if $\text{Im}(\beta) \not\subset \mathcal{E}$,  \\
      \overline{p}_{\mathcal{F}} + \frac{\delta}{\rk_{\mathcal{F}}} \geq \overline{p}_{\mathcal{E}} + \frac{\delta}{\rk_{\mathcal{E}}}  & if $\text{Im}(\beta) \subset \mathcal{E}$. 
    \end{cases*}\]
\end{enumerate}
\end{prop}
\begin{proof}
\textit{$\big[$semistability $\Rightarrow$ (a) $+$ (b)$\big]$:}
Suppose that (a) is not satisfied, so $\beta =0$. Then the filtration $f = (\mathcal{F}_m)_{m \in \mathbb{Z}}$ given by
\[    \mathcal{F}_m \vcentcolon = \begin{cases*}
      0 & if $m > -1$, \\
      \mathcal{F}        & if $ -1 \geq m$.
    \end{cases*}
    \]
is a destabilizing filtration for the pair $(\mathcal{F}, \beta)$.

On the other hand, suppose that (b) is not satisfied. Let $\mathcal{E} \subset \mathcal{F}$ be a proper saturated sheaf violating condition (b) above. We show that $(\mathcal{F}, \beta)$ is unstable in each case.
\begin{enumerate}[(C1)]
    \item Assume that $\text{Im}(\beta) \not\subset \mathcal{E}$. Let $f= (\mathcal{F}_m)_{m \in \mathbb{Z}}$ denote the filtration of $(\mathcal{F}, \beta)$ defined by
    \[\mathcal{F}_m = \begin{cases*} 0 & if $m > 1$, \\
    \mathcal{E} & if $m =1$, \\
    \mathcal{F} & if $0 \geq m$.
    \end{cases*}\]
    Then we have $\nu^{(\delta)}(f) = \frac{(\overline{p}_{\mathcal{E}} - \frac{\delta}{\rk_{\mathcal{F}}} - \overline{p}_{\mathcal{F}}) \cdot \rk_{\mathcal{E}}}{\sqrt{\rk_{\mathcal{E}}}}$. By assumption this is strictly positive, and hence $f$ is a destabilizing filtration.
    
    \item Assume that $\text{Im}(\beta) \subset \mathcal{E}$. Let $f= (\mathcal{F}_m)_{m \in \mathbb{Z}}$ denote the filtration of $(\mathcal{F}, \beta)$ defined by
    \[\mathcal{F}_m = \begin{cases*} 0 & if $m > 0$, \\
    \mathcal{E} & if $m =0$, \\
    \mathcal{F} & if $0 > m$.
    \end{cases*}\]
    Then we have 
        \[  \nu^{(\delta)}(f) = \frac{(-\overline{p}_{\mathcal{F}/\mathcal{E}} + \frac{\delta}{\rk_{\mathcal{F}}} + \overline{p}_{\mathcal{F}}) \cdot \rk_{\mathcal{F}/\mathcal{E}}}{\sqrt{\rk_{\mathcal{F}/\mathcal{E}}}}\]
    Note that this can be rewritten as follows.
        \[  \nu^{(\delta)}(f) = \frac{\rk_{\mathcal{E}}}{\sqrt{\rk_{\mathcal{F}/\mathcal{E}}}} \cdot \left( \overline{p}_{\mathcal{E}} + \frac{\delta}{\rk_{\mathcal{E}}} - \overline{p}_{\mathcal{F}} - \frac{\delta}{\rk_{\mathcal{F}}}\right)\]
        By assumption this is strictly positive. Therefore $f$ is a destabilizing filtration. This concludes the ``only if" direction.
\end{enumerate}
\medskip
\noindent\textit{$\big[$(a) $+$ (b) $\Rightarrow$ semistability)$\big]$:}
\medskip

Suppose that $(\mathcal{F}, \beta)$ satisfies conditions (a) and (b) in the statement of the proposition. Let $f = (\mathcal{F}_{m})_{m \in \mathbb{Z}}$ be a filtration of $(\mathcal{F}, \beta)$. We shall show that $f$ is not destabilizing. By definition, the numerical invariant is given by
\begin{equation} \label{eqn 3}
\nu^{(\delta)}(f)  = \frac{1}{\sqrt{\left(\sum_{m \in \mathbb{Z}} \rk_{\mathcal{F}_m/\mathcal{F}_{m+1}} \cdot m^2\right)}} \cdot \left( \sum_{m \geq 1} m \cdot \left(\overline{p}_{\mathcal{F}_m/ \mathcal{F}_{m+1}} - \frac{\delta}{\rk_{\mathcal{F}}} - \overline{p}_{\mathcal{F}}\right) \cdot \rk_{\mathcal{F}_m/\mathcal{F}_{m+1}} \right.
\end{equation}
\[\hfill \left. + \sum_{m \leq -1, \, \mathcal{F}_{m+1} \neq \mathcal{F}} m \cdot \left(\overline{p}_{\mathcal{F}_m/ \mathcal{F}_{m+1}} - \frac{\delta}{\rk_{\mathcal{F}}} - \overline{p}_{\mathcal{F}}\right) \cdot \rk_{\mathcal{F}_m/\mathcal{F}_{m+1}} \right)
\]
Two applications of summation by parts produce the following two equalities.
\begin{gather} \label{eqn 1} \sum_{m \geq 1} m \cdot \left(\overline{p}_{\mathcal{F}_m/ \mathcal{F}_{m+1}} - \frac{\delta}{\rk_{\mathcal{F}}} - \overline{p}_{\mathcal{F}}\right) \cdot \rk_{\mathcal{F}_m/\mathcal{F}_{m+1}} = \sum_{m \geq 1} \left(\overline{p}_{\mathcal{F}_m} - \frac{\delta}{\rk_{\mathcal{F}}} - \overline{p}_{\mathcal{F}}\right)\cdot \rk_{\mathcal{F}_m}
\end{gather}
\begin{gather*} \sum_{m \leq -1, \, \mathcal{F}_{m+1} \neq \mathcal{F}} m \cdot \left(\overline{p}_{\mathcal{F}_m/ \mathcal{F}_{m+1}} - \frac{\delta}{\rk_{\mathcal{F}}} - \overline{p}_{\mathcal{F}}\right) \cdot \rk_{\mathcal{F}_m/\mathcal{F}_{m+1}} = -\sum_{m \leq -1, \, \mathcal{F}_{m+1} \neq \mathcal{F}} \left(\overline{p}_{\mathcal{F}/\mathcal{F}_{m+1}} - \frac{\delta}{\rk_{\mathcal{F}}} - \overline{p}_{\mathcal{F}}\right)\cdot \rk_{\mathcal{F}/\mathcal{F}_{m+1}}
\end{gather*}
We can use (\ref{eqn 1}) to rewrite the summations appearing in formula (\ref{eqn 3}). This yields
\begin{equation} \label{eqn 4}
    \nu^{(\delta)}(f) = \frac{1}{\sqrt{\left(\sum_{m \in \mathbb{Z}} \rk_{\mathcal{F}_m/\mathcal{F}_{m+1}} \cdot m^2\right)}} \cdot \Biggl( \sum_{m \geq 1} \left(\overline{p}_{\mathcal{F}_m} - \frac{\delta}{\rk_{\mathcal{F}}} - \overline{p}_{\mathcal{F}}\right) \cdot \rk_{\mathcal{F}_m}  \\
    + \sum_{m \leq -1, \, \mathcal{F}_{m+1} \neq \mathcal{F}}  \left(\overline{p}_{\mathcal{F}/\mathcal{F}_{m+1}} - \frac{\delta}{\rk_{\mathcal{F}}} - \overline{p}_{\mathcal{F}}\right)\cdot \rk_{\mathcal{F}/\mathcal{F}_{m+1}} \Biggl) \notag
\end{equation}
Let us further rewrite the terms in the second summation appearing in (\ref{eqn 4}) above. For each $m \leq -1$, we have $\text{Im}(\beta) \subset \mathcal{F}_{m+1}$. Since $\beta \neq 0$, we have $\mathcal{F}_{m+1} \neq 0$. If in addition $\mathcal{F}_{m+1} \neq \mathcal{F}$, then we have
\begin{gather} \label{eqn 2} -\left(\overline{p}_{\mathcal{F}/\mathcal{F}_{m+1}} - \frac{\delta}{\rk_{\mathcal{F}}} - \overline{p}_{\mathcal{F}}\right)\cdot \rk_{\mathcal{F}/\mathcal{F}_{m+1}} =  \left(\overline{p}_{\mathcal{F}_{m+1}} + \frac{\delta}{\rk_{\mathcal{F}_{m+1}}} - \overline{p}_{\mathcal{F}}  - \frac{\delta}{\rk_{\mathcal{F}}} \right)\cdot \rk_{\mathcal{F}_{m+1}}
\end{gather}
Using (\ref{eqn 2}), we rewrite equation (\ref{eqn 4}) as
\begin{equation} \label{eqn 5}
    \nu^{(\delta)}(f) = \frac{1}{\sqrt{\left(\sum_{m \in \mathbb{Z}} \rk_{\mathcal{F}_m/\mathcal{F}_{m+1}} \cdot m^2\right)}} \cdot \Biggl( \sum_{m \geq 1} \left(\overline{p}_{\mathcal{F}_m} - \frac{\delta}{\rk_{\mathcal{F}}} - \overline{p}_{\mathcal{F}}\right) \cdot \rk_{\mathcal{F}_m} \\ \notag
    + \sum_{m \leq -1, \, \mathcal{F}_{m+1} \neq \mathcal{F}}  \left(\overline{p}_{\mathcal{F}_{m+1}} + \frac{\delta}{\rk_{\mathcal{F}_{m+1}}} - \overline{p}_{\mathcal{F}}  - \frac{\delta}{\rk_{\mathcal{F}}} \right)\cdot \rk_{\mathcal{F}_{m+1}} \Biggl)
\end{equation}
In order to see that $f$ is not destabilizing, we shall show that each term in the two summations appearing in (\ref{eqn 5}) is $\leq 0$. 
\begin{enumerate}[1.]

\item $\text{We show} \; \sum_{m \leq -1, \, \mathcal{F}_{m+1} \neq \mathcal{F}}  \left(\overline{p}_{\mathcal{F}_{m+1}} + \frac{\delta}{\rk_{\mathcal{F}_{m+1}}} - \overline{p}_{\mathcal{F}}  - \frac{\delta}{\rk_{\mathcal{F}}} \right)\cdot \rk_{\mathcal{F}_{m+1}} \leq 0$. For all $m \leq -1$ with $\mathcal{F}_{m+1} \neq \mathcal{F}$, we know that $\text{Im}(\beta) \subset \mathcal{F}_{m+1} \subsetneq \mathcal{F}$. Therefore, by condition (b), we have
\[\left(\overline{p}_{\mathcal{F}_{m+1}} + \frac{\delta}{\rk_{\mathcal{F}_{m+1}}} - \overline{p}_{\mathcal{F}}  - \frac{\delta}{\rk_{\mathcal{F}}} \right)\cdot \rk_{\mathcal{F}_{m+1}} \leq 0\]
for each term in the sum. 

\item $\text{We show} \;\sum_{m \geq 1} \left(\overline{p}_{\mathcal{F}_m} - \frac{\delta}{\rk_{\mathcal{F}}} - \overline{p}_{\mathcal{F}}\right) \cdot \rk_{\mathcal{F}_m} \leq 0$. Let $m \geq 1$. If $\text{Im}(\beta) \not\subset \mathcal{F}_m$, then condition (b) implies that
\[\left(\overline{p}_{\mathcal{F}_m} - \frac{\delta}{\rk_{\mathcal{F}}} - \overline{p}_{\mathcal{F}}\right) \cdot \rk_{\mathcal{F}_m} \leq 0,\]
On the other hand if $\text{Im}(\beta) \subset \mathcal{F}_m$, then we know that $\mathcal{F}_m \neq 0$ and
\[ \left(\overline{p}_{\mathcal{F}_{m}} + \frac{\delta}{\rk_{\mathcal{F}_{m}}} - \overline{p}_{\mathcal{F}}  - \frac{\delta}{\rk_{\mathcal{F}}} \right)\cdot \rk_{\mathcal{F}_{m}} \leq 0. \]
Using the fact that $\delta > 0$, we get the chain of inequalities
\begin{gather*}
\left(\overline{p}_{\mathcal{F}_m} - \frac{\delta}{\rk_{\mathcal{F}}} - \overline{p}_{\mathcal{F}}\right) \cdot \rk_{\mathcal{F}_m} \leq \left(\overline{p}_{\mathcal{F}_{m}} + \frac{\delta}{\rk_{\mathcal{F}_{m}}} - \overline{p}_{\mathcal{F}}  - \frac{\delta}{\rk_{\mathcal{F}}} \right)\cdot \rk_{\mathcal{F}_{m}} \leq 0,
\end{gather*}
thus concluding the proof.
\end{enumerate}
\end{proof}

\begin{example} \quad
\begin{enumerate}[(1)]
    \item Suppose that $\delta>0$ with $\text{deg}(\delta) \leq -1$. In this case, Proposition \ref{prop: comparison of stability pairs} implies that a pair $(\mathcal{F}, \beta)$ is $\nu^{(\delta)}$-semistable if and only if $\beta \neq 0$ and $\mathcal{F}$ is Gieseker semistable. The resulting moduli space of $\Pair_{\mathcal{A}}^{d}(X)^{\nu^{(\delta)} \dash \mathrm{ss}}$ is projective over the moduli space of Gieseker semistable pure sheaves on $X$.
    
    \item Fix the choice of a Hilbert polynomial $P$. We restrict to the moduli stack $\Pair_{\mathcal{A}}^{d}(X)_{P}$ parametrizing pairs with Hilbert polynomial $P$. Let $\delta$ be a positive (constant) rational number. If $\delta$ is small enough, then Proposition \ref{prop: comparison of stability pairs} implies that a pair $(\mathcal{F}, \beta) \in \Pair_{\mathcal{A}}^{d}(X)_{P}$ is semistable if and only if $\beta \neq 0$, the sheaf $\mathcal{F}$ is Gieseker semistable, and all proper subsheaves $\mathcal{E} \subset \mathcal{F}$ with $\overline{p}_{\mathcal{E}} = \overline{p}_{\mathcal{F}}$ satisfy $\text{Im}(\beta) \not\subset \mathcal{E}$. If $\mathcal{A}$ is Gieseker semistable with the same reduced Hilbert polynomial as $\mathcal{F}$, then this last condition is equivalent to requiring that $\mathcal{F}/ \text{Im}(\beta)^{sat}$ is Gieseker stable.
    
    \item Suppose that $S = Spec(k)$ for a field $k$. Let $X$ be a smooth projective geometrically connected curve over $k$, equipped with the choice of a polarization $\mathcal{O}_{X}(1)$. Set $\mathcal{A} = \mathcal{O}_{X}$ and take $d=1$. Then $\Pair_{\mathcal{O}_X}^{1}(X)$ is the stack classifying vector bundles on $X$ along with the choice of a section. We can take $\delta$ to be a constant in $\mathbb{Q}$. If $\delta >0$, then Proposition \ref{prop: comparison of stability pairs} shows that $\nu^{(\delta)}$-stability coincides with stability for Bradlow pairs as considered by Thaddeus in his work on the Verlinde formula \cite{thaddeus-verlinde} \footnote{Using the notation in \cite{thaddeus-verlinde}, the precise comparison is $\delta = 2 \sigma$.}. Therefore, we recover the moduli space of $\delta$-semistable Bradlow pairs.
\end{enumerate}
\end{example}
\end{subsection}
\end{section}

\bibliography{torsion_free_grassmannian.bib}
\bibliographystyle{alpha}
\end{document}